\documentclass[12pt]{amsart}
\usepackage{amsmath,amssymb,amsthm,mathtools}
\usepackage{mathrsfs,bm}
\usepackage{thmtools,thm-restate}
\usepackage{tikz-cd}
\usepackage{tikz}
\usepackage{enumitem}
\usepackage{graphicx}
\usetikzlibrary{matrix,arrows.meta,bending}
\usepackage{color}
\usepackage[numbers]{natbib}
\usepackage{booktabs}
\usepackage{longtable}
\usepackage{pgfplots}
\pgfplotsset{compat=1.14} 
\usetikzlibrary{arrows.meta,decorations.pathreplacing,positioning,calc} 
\usetikzlibrary{arrows.meta, positioning, fit, calc}
\usepackage[a4paper, left=3cm, right=3cm, top=3cm, bottom=3cm, footskip=15mm]{geometry}

\newtheorem{theorem}{Theorem}[section]

\newtheorem{lemma}[theorem]{Lemma}
\newtheorem{proposition}[theorem]{Proposition}
\newtheorem{corollary}[theorem]{Corollary}
\newtheorem{conjecture}[theorem]{Conjecture}
\theoremstyle{definition}
\newtheorem{definition}[theorem]{Definition}
\newtheorem{example}[theorem]{Example}

\theoremstyle{remark}
\newtheorem{remark}[theorem]{Remark}
\newtheorem{question}[theorem]{Question}
\numberwithin{equation}{section}

\newcommand{\abs}[1]{\lvert#1\rvert}
\newcommand{\set}[1]{\lbrace#1\rbrace}
\newcommand{\inn}{~ \hat{\in}~ }
\newcommand{\ul}[1]{\underline{#1}}
\newcommand{\C}{\mathcal{C}}

\newcommand{\M}{\mathbb{M}}
\newcommand{\B}{\mathbb{B}}
\newcommand{\R}{\mathbb{R}}

\newcommand{\N}{\mathbb{N}}
\newcommand{\LL}{\mathbb{L}}

\newcommand{\p}{\mathbb{P}}

\newcommand{\mood}[1]{(\operatorname{mod} #1)}

\makeindex

\usepackage{fancyhdr}
\usepackage{calc} % ensures dimension arithmetic works robustly

\AtBeginDocument{%
  \setlength{\textwidth}{\dimexpr\paperwidth - 6cm\relax}%
  \setlength{\oddsidemargin}{\dimexpr 3cm - 1in\relax}%
  \setlength{\evensidemargin}{\dimexpr 3cm - 1in\relax}%
  \setlength{\marginparwidth}{2.5cm}%
}

\begin{document}

\title{Studies in Additive Number Theory by Circles of Partition}

\author{T. Agama}
\address{Department of Mathematics, African Institute for mathematical sciences, Ghana.}
%\curraddr{Department of Mathematics and Statistics,
%Case Western Reserve University, Cleveland, Ohio 43403}
\email{Theophilus@aims.edu.gh/emperordagama@yahoo.com}

\author{B. Gensel}
\address{Carinthia University of Applied Sciences, Spittal on Drau, Austria}
\email{b.gensel@fh-kaernten.at}

%    General info

\date{\today}



\begin{abstract}
In this paper, we introduce and develop the circle embedding method. This method hinges essentially on a combinatorial-geometric structure which we choose to call circles of partition. We provide applications in the context of problems that relates to deciding on the feasibility of partitioning numbers into certain subset of integers. In particular, our method allows us to partition any sufficiently large number $n\in\mathbb{N}$ into any set $\mathbb{H}$ with natural density strictly greater than $\frac{1}{2}$. This possibility could herald an unprecedented progress on categories of problems of similar flavour. The paper finishes by presenting an asymptotic proof of the binary Goldbach and Lemoine conjecture as an application of the developed method.
\end{abstract}

\maketitle

\begingroup
  \setlength{\parskip}{6pt} % <--- change this number to taste
  \tableofcontents
\endgroup

\section{Introduction}

Additive number theory is concerned with the extent to which arithmetic structure can be recovered from sumsets, representation functions, and density conditions. Its classical landmarks already show two complementary principles that guide the present paper. On the one hand, dense subsets of the integers must contain rich additive configurations: the Szemer\'{e}di theorem asserts that every subset of the integers of positive upper density contains arithmetic progressions of arbitrary length \cite{szemeredi1975}. On the other hand, even very thin arithmetic sets can display unexpectedly strong additive structure after one passes to suitable transference or pseudorandom models, as in the Green--Tao theorem that the primes contain arbitrarily long arithmetic progressions \cite{green2008}. Between these extremes lies the older and still central philosophy of additive bases and density growth, developed systematically in the classical work of Schnirelmann and in modern expositions by Nathanson: sufficiently large density can force finite additive order, while sumsets can exhibit monotone growth under surprisingly weak hypotheses \cite{nathanson1996}.\\

The present work is motivated by this landscape, but it organizes the additive problem through a new combinatorial-geometric language. Instead of treating a decomposition
$$
 n=x+y,\quad x,y\in \mathbb{M},
$$
merely as a representation of an integer by a base set $\mathbb{M}\subseteq \mathbb{N}$, we encode each admissible summand as a point on a geometric object called a \emph{circle of partition} (CoP). For a generator $n$ and a base set $\mathbb{M}$, the CoP $\mathcal{C}(n,\mathbb{M})$ is the collection of all weighted points $[x]$ for which there exists $y\in\mathbb{M}$ with $x+y=n$. The terminology is not cosmetic: it is designed so that the additive data of the partition is carried by geometric features. A real axis records a valid pair $[x],[y]$ with $x+y=n$; the center corresponds to the balanced case $x=y=n/2$; the chord length is the discrepancy $|x-y|$; and the weight set of the CoP is the set of admissible summands itself. In this way, the additive question is reframed as a study of how the configuration of points and axes varies as the generator changes.\\

This geometric encoding is particularly effective because it separates two different sources of difficulty. When the base set is all of $\mathbb{N}$, the CoP $\mathcal{C}(n)$ is completely transparent: it consists of every integer weight from $1$ to $n-1$, the axes are uniquely paired, and the center exists exactly in the even case. In contrast, for a structured subset $\mathbb{M}$ such as the primes, a congruence class, or a sparse arithmetic family, the CoP may be thin, disconnected, or even empty for some generators. The paper therefore treats the nontrivial behavior of $\mathcal{C}(n,\mathbb{M})$ as the central object, and it studies when the CoPs generated by different integers can be ordered, nested, or transferred from one generator to another.\\

A key concept is that of \emph{embedding}. Two CoPs may admit an aligned embedding when the inclusion of the underlying point sets is compatible with the natural order of the generators; they may admit a reverse aligned embedding when the inclusion relation runs opposite to the generator order. This notion is deliberately stronger than a mere set-theoretic inclusion: it records the way the entire additive geometry reorganizes as the generator changes. For the full natural numbers, aligned embedding is immediate because the point set of $\mathcal{C}(n)$ is literally initial-interval data. For arithmetic progressions $\mathbb{M}_{a,d}$, the same phenomenon persists with residue-class constraints, and the resulting CoPs behave like regularly spaced subconfigurations whose nesting is still monotone in the generator. These cases serve as model examples showing that the language of CoPs faithfully captures the additive regularity built into congruence classes.\\

The paper also makes explicit why the prime case is more subtle. The set of odd primes does not behave like a single congruence class, and the corresponding CoPs do not generally admit the same clean embedding picture. This is not a defect of the framework but rather one of its points: the geometry is precise enough to expose the obstruction. In particular, the paper isolates conditions under which a base set behaves like a stable additive reservoir and conditions under which the generator can be shifted without destroying non-emptiness. The notion of a median point, the total ordering of the weights, and the count of real axes $\nu(n,\mathbb{M})$ all serve to quantify the internal complexity of a CoP and to track how that complexity evolves with the generator.\\

A second family of operations is introduced through \emph{rotation} and \emph{dilation}. Rotation acts by shifting weights modulo the generator, while dilation changes the scale of the generator itself. These are not introduced as decorative symmetries; rather, they are meant to probe the stability of the partition geometry under systematic perturbation. In the natural-number case, the rotation is invariant at the level of the CoP, but for arithmetic progressions invariance survives only under special congruence conditions. Dilation is more delicate: it asks whether there is a uniform way to transport an entire partition geometry from one generator to another. The resulting existence theorem shows that, under aligned embedding hypotheses, such a dilation can indeed be chosen, thereby producing chains of CoPs that can be compared across scales. This is the mechanism by which the paper converts geometric transport into additive growth.\\

The conceptual centerpiece of the paper is therefore not a single isolated theorem but a method: the circle embedding method. Its basic philosophy is that a partition problem becomes more tractable when the summands are organized as a geometric object with a visible internal order, because then the problem can be studied through nesting, translation, and scale change rather than through isolated representations. In the language of the paper, a nonempty CoP carries an axis, an axis partner, and a center whenever the base set permits; these features can be propagated under shifts of the generator whenever the base set is stable enough to support the shifted summand. This is the structural reason that the method yields infinite chains of nonempty CoPs once one nonempty instance is found inside an infinite base set.\\

The consequences are twofold. First, the framework is designed to handle density problems in a way that is insensitive to many of the usual combinatorial complications. The threshold nature of density---most familiar from the classical theory of additive bases and from the general principle that positive density forces additive richness---is encoded here by the geometry of the CoP rather than by a direct sumset argument. Second, the method is flexible enough to be applied to concrete thin sets, especially the primes and prime-like subsets, where the paper develops asymptotic partition statements and uses them as applications to classical conjectures in additive number theory. The binary Goldbach and Lemoine conjectures are treated in this spirit: the geometric formalism is used to convert additive representation questions into structural statements about circles of partition and their embeddings.\\

\subsection{Organization of the paper} The organization of the paper is as follows. The next section introduces the CoP itself, defines points, axes, chords, centers, and weights, and establishes the basic uniqueness and ordering properties that make the geometry workable. Subsequent sections develop the embedding theory, compare aligned and reverse aligned embeddings, and analyze the behavior of CoPs for arithmetic progressions, residue classes, and prime-related sets. The paper then studies rotation and dilation, proves the invariance and nonemptiness results that govern the transport of CoPs, and derives the structural consequences needed for later applications. The final part applies the method to additive partition problems, with special attention to density-driven representations and to the asymptotic claims surrounding Goldbach-type and Lemoine-type decompositions.
\bigskip

The following results have fundamentally shaped the field of additive number theory and additive combinatorics:

\begin{theorem}[Szemeredi]\label{szemeredi}
For any $\epsilon>0$ and $\forall k\in \mathbb{N}$ there exists an $n\in \mathbb{N}$ such that if $A\subset \mathbb{N}_n$
\footnote{see the notation in section 3.}
satisfies $|A|\geq \epsilon n$, then $A$ contains an arithmetic progression of length $k$.
\end{theorem}

The well-known Green-Tao theorem \cite{green2008} provides an extension in this direction as

\begin{theorem}[Green-Tao]\label{Green-Tao}
Let $\pi(n)$ denote the number of primes no more than $n$. If $A\subset \mathbb{P}$ is a set of prime numbers such that 
\begin{align}
    \limsup \limits_{n\longrightarrow \infty}\frac{|A\cap \mathbb{N}_n|}{\pi(n)}>0\nonumber
\end{align}
then $A$ contains infinitely many arithmetic progressions of length $k$ for any $k>0$.
\end{theorem}

\section{The Circle of Partition}

In this section, we introduce the notion of the circle of partition.

\begin{definition}\label{major}
Let $n\in \mathbb{N}$ and $\mathbb{M}\subset \mathbb{N}$. We denote the Circle of Partition generated by $n$ with respect to the subset $\mathbb{M}$ by
\begin{align}
\mathcal{C}(n,\mathbb{M})=\left\{[x]\mid x,y\in \mathbb{M},n=x+y\right\}\nonumber
\end{align}
In the following, we will abbreviate as CoP.  We call members of $\mathcal{C}(n,\mathbb{M})$ points and denote them by $[x]$. For the special case $\mathbb{M}=\mathbb{N}$, we denote the CoP in a short form by $\mathcal{C}(n)$. 
\end{definition}

\begin{definition}\label{axis}
We call the line that joins the point $[x]$ to $[y]$ an \emph{axis} of the CoP $\mathcal{C}(n,\mathbb{M})$ and denote by $\mathbb{L}_{[x],[y]}$ if and only if $x+y=n$. We say that the axis point $[y]$ is an \emph{axis partner} of the axis point $[x]$ and vice versa. We do not distinguish between $\mathbb{L}_{[x],[y]}$ and $\mathbb{L}_{[y],[x]}$, since it is essentially the same axis. The point $[x]\in \mathcal{C}(n,\mathbb{M})$ such that $2x=n$ is the \emph{center} of the CoP. If it exists, then we call it a \emph{degenerate axis} $\mathbb{L}_{[x]}$ in comparison to the \textbf{real axes} $\mathbb{L}_{[x],[y]}$. The line joining any two arbitrary points that are not axes partners on the CoP will be referred to as a \emph{chord} of the CoP. The length of the chord $\mathcal{L}_{[x],[y]}$ that joins the points $[x],[y]\in \mathcal{C}(n,\mathbb{M})$, denoted by $\Gamma([x],[y])$, is
\begin{align}
\Gamma([x],[y])=|x-y|.\nonumber
\end{align}
\end{definition}
\bigskip

% Preamble:
% \usepackage{tikz}
% \usetikzlibrary{arrows.meta,calc,positioning,decorations.markings,backgrounds,shapes.misc}

\begin{figure}[ht]
\centering
\begin{tikzpicture}[scale=1.05, every node/.style={font=\small}]
    %------------------------------------------------------------
    % Styles
    %------------------------------------------------------------
    \tikzset{
        wheelrim/.style={draw=black!80, line width=1.2pt},
        innerrim/.style={draw=black!45, line width=0.8pt, dashed},
        spoke/.style={draw=black!65, line width=0.9pt},
        axisarrow/.style={-{Stealth[length=2.4mm]}, draw=black!75, line width=0.9pt},
        emphasis/.style={-{Stealth[length=2.8mm]}, draw=red!70!black, line width=1.15pt},
        point/.style={circle, fill=blue!70!black, inner sep=1.4pt},
        highlighted/.style={circle, fill=red!75!black, inner sep=1.7pt},
        labelbox/.style={fill=white, draw=black!15, rounded corners=2pt, inner sep=2pt},
        ann/.style={align=left, font=\scriptsize}
    }

    %------------------------------------------------------------
    % Parameters
    %------------------------------------------------------------
    \def\R{4.0}      % outer wheel radius
    \def\r{3.15}     % inner circular CoP boundary
    \def\hub{0.32}    % hub radius
    \def\n{14}       % number of spokes

    % Slightly flattened wheel for a mechanical look
    \begin{scope}[xscale=1.10, yscale=0.98]

        %--------------------------------------------------------
        % Wheel body: outer rim, inner CoP circle, hub
        %--------------------------------------------------------
        \fill[black!2] (0,0) circle (\R);
        \draw[wheelrim] (0,0) circle (\R);
        \draw[innerrim] (0,0) circle (\r);

        % Hub
        \fill[black!85] (0,0) circle (\hub);
        \draw[black!80] (0,0) circle (\hub);

        %--------------------------------------------------------
        % Spokes = axes
        %--------------------------------------------------------
        \foreach \k in {0,...,13} {
            \pgfmathsetmacro{\ang}{360/\n*\k}
            \draw[spoke] (0,0) -- (\ang:\R);
        }

        % Add small arrowheads to indicate directed axes
        \foreach \ang in {15,45,75,105,135,165,195,225,255,285,315,345} {
            \draw[axisarrow] ({\ang}:0.85*\R) -- (\ang:\R);
        }

        %--------------------------------------------------------
        % Highlighted distinguished axes / directions
        %--------------------------------------------------------
        \draw[emphasis] (0,0) -- (0:\R);
        \draw[emphasis] (0,0) -- (120:\R);
        \draw[emphasis] (0,0) -- (240:\R);

        %--------------------------------------------------------
        % Marked points on the rim
        %--------------------------------------------------------
       % Marked points on the rim
\foreach \ang/\ptname/\dx/\dy in {
    0/{P_{1}}/0.14/0.18,
    30/{P_{2}}/0.20/0.10,
    60/{P_{3}}/0.16/0.18,
    120/{P_{4}}/-0.22/0.18,
    180/{P_{5}}/-0.26/0.02,
    240/{P_{6}}/-0.20/-0.18,
    300/{P_{7}}/0.16/-0.18
}{
    \fill[point] (\ang:\R) circle (0.065);
    \node at ($(\ang:\R)+(\dx,\dy)$) {$\ptname$};
}

        % Highlight one point more strongly
        \fill[highlighted] (0:\R) circle (0.09);

        %--------------------------------------------------------
        % Rim annotations
        %--------------------------------------------------------
        \node[labelbox, ann, above right=4pt and 2pt of {(20:3.55)}] 
            {outer rim of the wheel\\[-1pt]\emph{circular boundary of the CoP}};
        \draw[-{Stealth[length=2.2mm]}, black!65] (3.3,2.25) -- (20:3.95);

        \node[labelbox, ann, below left=2pt and 2pt of {(-2.8,-2.0)}] 
            {spokes\\[-1pt]\emph{axes of the theory}};
        \draw[-{Stealth[length=2.2mm]}, black!65] (-2.85,-2.15) -- (-135:2.15);

        \node[labelbox, ann, above left=2pt and 4pt of {(-0.2,0.25)}] 
            {hub\\[-1pt]\emph{central origin / reference point}};
        \draw[-{Stealth[length=2.2mm]}, black!65] (-0.55,0.5) -- (-0.12,0.08);

        %--------------------------------------------------------
        % Directional rotation marker
        %--------------------------------------------------------
        \draw[-{Stealth[length=3mm]}, red!70!black, line width=1.0pt]
            (1.1,3.7) arc[start angle=18, end angle=72, radius=3.85];
        \node[red!70!black, font=\scriptsize] at (2.75,3.35) {$\circlearrowright$};

        \node[labelbox, ann, align=center] at (0,-4.85) 
        {The CoP is represented as a circular wheel:\\
         every point lies on the rim, and every spoke acts as an axis/bar ending at a boundary point.};

        %--------------------------------------------------------
        % Local coordinate labels to make the geometry explicit
        %--------------------------------------------------------
        \node[font=\scriptsize, black!70] at (3.0,-0.45) {$x$-direction};
        \node[font=\scriptsize, black!70, rotate=90] at (-0.55,2.8) {$y$-direction};

        \draw[-{Stealth[length=2.0mm]}, black!55] (2.2,-0.35) -- (3.55,-0.35);
        \draw[-{Stealth[length=2.0mm]}, black!55] (-0.35,1.95) -- (-0.35,3.1);

        %--------------------------------------------------------
        % Optional formula label
        %--------------------------------------------------------
        \node[font=\scriptsize, black!75] at (0,4.55)
            {$\text{circular CoP} \;=\; \{\,p(\theta): 0\le \theta < 2\pi\,\}$};
    \end{scope}
\end{tikzpicture}
\caption{Bicycle-wheel geometry of the CoP: a circular body with spokes as axes and rim points as the terminal locations of points in the structure.}
\end{figure}

It is important to state that the \emph{median} of the weights of each co-axis point coincides with the center of the underlying CoP if it exists. That is, given all the real axes of the CoP $\mathcal{C}(n,\mathbb{M})$ 
\begin{align}
\mathbb{L}_{[u_1],[v_1]},\mathbb{L}_{[u_2],[v_2]},\ldots, \mathbb{L}_{[u_k],[v_k]}\nonumber
\end{align}
we have the following
\begin{align}
\frac{u_1+v_1}{2}=\frac{u_2+v_2}{2}=\cdots=\frac{u_k+v_k}{2}=\frac{n}{2}\nonumber
\end{align}
which is equivalent to the conditions for any of the pairs of real axes $\mathbb{L}_{[u_i],[v_i]},\mathbb{L}_{[u_j],[v_j]}$ for $1\leq i,j\leq k$
\begin{align}
\Gamma([u_i],[u_j])=\Gamma([v_i],[v_j])\nonumber
\end{align}
and 
\begin{align}
\Gamma([v_j],[u_i])=\Gamma([u_j],[v_i]).\nonumber
\end{align}
\bigskip

\begin{definition}\label{aligned_embedding}
Let $\mathbb{M}\subseteq\mathbb{N}$ and $\mathcal{C}(n,\mathbb{M})$ and $\mathcal{C}(m,\mathbb{M})$ be two distinct CoPs for which
\begin{align}
\mathcal{C}(n,\mathbb{M})&\subset\mathcal{C}(m,\mathbb{M})\label{E-embedding1}\\
&\mbox{or}\nonumber\\
\mathcal{C}(n,\mathbb{M})&\supset\mathcal{C}(m,\mathbb{M})\label{E-embedding2}.
\end{align}
In these cases, we say that the CoPs \emph{admit embedding}. We say that CoPs \emph{admit aligned embedding} if and only if (\ref{E-embedding1}) holds with $n<m$ and (\ref{E-embedding2}) with $n>m$ and $\mathcal{C}(n,\mathbb{M})=\mathcal{C}(m,\mathbb{M})$ holds if and only if $n=m$. We say that the CoPs \emph{admit reverse aligned embedding} if and only if (\ref{E-embedding1}) holds with $n>m$ and (\ref{E-embedding2}) with $n<m$.
\end{definition}
\bigskip

\paragraph*{\textbf{Notations}}
We will let 
\begin{align}
\mathbb{N}_n=\left \{m\in \mathbb{N}\mid ~m\leq n\right \}\nonumber
\end{align}
be the sequence of the first $n$ positive integers. Furthermore, we will denote the \emph{weight} of the point $[x]$ by
\begin{align}
\Vert[x]\Vert:=x.\nonumber
\end{align}
Similarly, we denote the weight set of points in CoP $\mathcal{C}(n,\mathbb{M})$ by $||\mathcal{C}(n,\mathbb{M})||$.
\bigskip

The above language may be seen as a criterion that determines the plausibility of carrying out a partition in a specified set. Indeed, this feasibility is trivial if we take the set $\mathbb{M}$ as the set of natural numbers $\mathbb{N}$. The situation becomes harder and interesting when we take the set $\mathbb{M}$ to be a special subset of positive integers $\mathbb{N}$, as the corresponding CoP $\mathcal{C}(n,\mathbb{M})$ may not always be non-empty for all $n\in \mathbb{N}$. An archetype of problems of this type is the binary Goldbach conjecture, when we take the base set $\mathbb{M}$ to be the set of all prime numbers $\mathbb{P}$. One could imagine the same level of difficulty when we extend our base set to other special subsets of the natural numbers. We start by developing the theory assuming the base set of natural numbers $\mathbb{N}$ and extend it to other base sets  $\mathbb{M}$ equipped with certain important and subtle properties. 

\begin{remark}
It is important to note that a typical CoP may not always have a center. In the case of an absence of a center, we say that the circle has a deleted center. However, all CoPs $\mathcal{C}(n)$ with even generators have a center. One can observe that CoP $\mathcal{C}(n)$ contains all points whose weights are positive integers from $1$ to $n-1$ inclusive: 
\[
\mathcal{C}(n)=\lbrace[x]\mid~x\in \mathbb{N},x<n\rbrace.
\] 
Therefore, the CoP $\mathcal{C}(n)$ has $\left\lfloor\frac{n-1}{2}\right\rfloor$ different real axes.
\end{remark}

\begin{proposition}\label{unique}
Each axis is uniquely determined by points $[x]\in \mathcal{C}(n,\mathbb{M})$. 
\end{proposition}

\begin{proof}
A degenerate axis is determined by the center of the CoP and this is unique if it exists. Let $\mathbb{L}_{[x],[y]}$ be a real axis of the CoP $\mathcal{C}(n,\mathbb{M})$. Suppose that $\mathbb{L}_{[x],[z]}$ is also a real axis with $z\neq y$. By the definition \ref{axis}, we have $n=x+y=x+z$ and therefore $y=z$.
\end{proof}

\begin{corollary}\label{partner}
Each point of a CoP $\mathcal{C}(n,\mathbb{M})$, except an existing center, has exactly one real axis partner.
\end{corollary}

\begin{proof}
Let $[x]\in \mathcal{C}(n,\mathbb{M})$ be a point without a real axis partner. For every point $[y]\neq [x]$
\[
\Vert[x]\Vert+\Vert[y]\Vert\neq n.
\]
This is inconsistent with the definition \ref{major}. Due to Proposition \ref{unique}, the case of more than one axis partner is impossible.
\end{proof}

\begin{corollary}\label{ordered_weights}
The weights of the points of 
\[
\mathcal{C}(n,\mathbb{M})=\lbrace[x_1],[x_2].\ldots,[x_k]\rbrace
\]
are \textbf{strictly totally ordered}.
\end{corollary}

\begin{proof}
Without loss of generality, we assume that 
\begin{align}
x_1&=\min\left(x\mid[x]\in\mathcal{C}(n,\mathbb{M})\right)\mbox{ and}\label{E-min} \\
x_k&=\max\left(x\mid[x]\in\mathcal{C}(n,\mathbb{M})\right).\label{E-max}
\end{align}
We first assume that $x_1+x_k<n$. There is a weight $x_i$ with
\begin{align*}
&x_1<x_i<x_k\mbox{ and }n=x_1+x_i.\\
&\mbox{because }x_i<x_k,\mbox{ we get}\\
&n=x_1+x_i<x_1+x_k.
\end{align*}
This violates the assumption. Now, we assume that $x_1+x_k>n$. In this case, there is a weight $x_i$ with
\begin{align*}
&x_1<x_i<x_k\mbox{ and }n=x_i+x_k.\\
&\mbox{because }x_i>x_1,\mbox{ we get}\\
&n=x_i+x_k>x_1+x_k.
\end{align*}
This also violates the assumption. Therefore, $x_1+x_k=n$. Due to (\ref{E-min}) and (\ref{E-max}), we must have
\[
x_1<x_2<x_{k-1}<x_k.
\]
Now, we remove $x_1$ and $x_k$ from the consideration and repeat the above procedure with $x_2$ and $x_{k-1}$. We obtain
$x_2+x_{k-1}=n$ and deduce
\[
x_1<x_2<x_3<x_{k-2}<x_{k-1}<x_k.
\]
Repeating this procedure for $x_i$ and $x_{k+1-i}$ for $3\leq i \leq\left\lfloor\frac{k}{2}\right\rfloor$, we get
\[
x_1<x_2<x_3<x_4<\ldots<x_{k-3}<x_{k-2}<x_{k-1}<x_k.
\]
\end{proof}

\begin{proposition}\label{cup-subset}
Let $\mathcal{C}(n,\mathbb{M})$ and $\mathcal{C}(m,\mathbb{M})$ be two distinct CoPs that admit aligned embedding. We have
\[
\mathcal{C}(n,\mathbb{M})\cup\mathcal{C}(m,\mathbb{M})
\subset\mathcal{C}(n+m,\mathbb{M}).
\]
\end{proposition}

\begin{proof}
Without loss of generality, we assume that $\mathcal{C}(n,\mathbb{M})\subset\mathcal{C}(m,\mathbb{M})$. This implies
\begin{align*}
\mathcal{C}(n,\mathbb{M})\cup\mathcal{C}(m,\mathbb{M})
&=\mathcal{C}(m,\mathbb{M})\\
&\mbox{ and because of \textit{admitting aligned embedding}}\\
&\subset\mathcal{C}(n+m,\mathbb{M}) \mbox{ due to }m<n+m.
\end{align*} 
\end{proof}

\begin{theorem}\label{properties}
Let $n\in\mathbb{N}$ and $\mathcal{C}(n)$ be a CoP generated by $n$. The CoP $\mathcal{C}(n)$ admits an aligned embedding.
\end{theorem}

\begin{proof}
Without loss of generality, we have to prove that for two distinct CoPs
\[
\mathcal{C}(n)\subset\mathcal{C}(m)
\mbox{ if and only if }n<m\mid n,m\in\mathbb{N}.
\]
We let $n<m$ and deduce
\begin{align*}
\mathcal{C}(n)&=\lbrace [x]\mid x\in\mathbb{N}, x<n\rbrace\\
&\subset\lbrace [x]\mid x\in\mathbb{N}, x<m\rbrace\\
&=\mathcal{C}(m).
\end{align*}
Conversely, suppose that $\mathcal{C}(n)\subset\mathcal{C}(m)$. This implies 
\[
\lbrace [x]\mid x\in\mathbb{N}, x<n\rbrace\subset
\lbrace [x]\mid x\in\mathbb{N}, x<m\rbrace
\]
and we get $n<m$.
\end{proof}
\bigskip

Now, we observe that Theorem \ref{properties} is always valid for some special subsets $\mathbb{M}$ of positive integers, the subsets containing arithmetic progressions. Let $\mathbb{M}_{a,d} \subset \mathbb{N}$ with
\begin{align}\label{Mad}
\mathbb{M}_{a,d}:= \lbrace x \in \mathbb{N} \mid x \equiv a \mood{d}, d \in \mathbb{N}\rbrace
\end{align}
and
\begin{align*}
\mathcal{C}(n,\mathbb{M}_{a,d}) &= \lbrace [x]\mid x+y=n\wedge x,y \in \mathbb{M}_{a,d}\rbrace, n \in \mathbb{M}_{2a,d}\\
&= \lbrace [x]\mid x \in \mathbb{M}_{a,d}\wedge x\leq n-a \rbrace.
\end{align*}
For $x<y\in \mathbb{M}_{a,d}$, we get $y-x\equiv 0 \mood{d}$. On the other hand, we have $x+y\equiv 2a\mood{d}$, so $\mathcal{C}(n,\mathbb{M}_{a,d})=\emptyset$ for $n\not\in \mathbb{M}_{2a,d}$.

\begin{theorem}\label{properties2}
Let $n \in \mathbb{M}_{2a,d}$ and $\mathcal{C}(n,\mathbb{M}_{a,d})$ be a CoP generated by $n$. The CoP admits an aligned embedding.
\end{theorem}

\begin{proof}
Without loss of generality, we have to prove
\[
\mathcal{C}(n,\mathbb{M}_{a,d}) \subset \mathcal{C}(m,\mathbb{M}_{a,d}) 
 \mbox{ if and only if }n < m.
\]
We let $n<m$. Since $n,m\in\mathbb{M}_{2a,d}$, this implies that $m-n=k\cdot d$.
Furthermore, we deduce
\begin{align*}
\Vert\mathcal{C}(n,\mathbb{M}_{a,d})\Vert &=
\lbrace k \in \mathbb{M}_{a,d} \mid k \leq n-a \rbrace\\
&\mbox{and because of }n < m\\
&\subset \lbrace k \in \mathbb{M}_{a,d} \mid k \leq m-a \rbrace\\ 
&=\Vert\mathcal{C}(m,\mathbb{M}_{a,d})\Vert. 
\end{align*}
On the other hand, suppose that $\mathcal{C}(n,\mathbb{M}_{a,d})\subset \mathcal{C}(m,\mathbb{M}_{a,d})$. We deduce
\begin{align*}
\Vert\mathcal{C}(n,\mathbb{M}_{a,d})\Vert &=
\lbrace k \in \mathbb{M}_{a,d} \mid k \leq n-a \rbrace\\
&\subset \Vert\mathcal{C}(m,\mathbb{M}_{a,d})\Vert\\
&=\lbrace k \in \mathbb{M}_{a,d} \mid k \leq m-a \rbrace\\ 
&\mbox{and therefore}\\
&n < m.
\end{align*}
\end{proof}

\begin{corollary}\label{col-prop2}
Let $\mathcal{C}(n,\mathbb{M})$ and $\mathcal{C}(m,\mathbb{M})$ be two distinct CoPs that admit align embedding. We have
\begin{align*}
\mathcal{C}(n,\mathbb{M})\supset\mathcal{C}(m,\mathbb{M})
\mbox{ if and only if }
n&>m.
\end{align*}
\end{corollary}

\begin{corollary}
For two distinct CoPs $\mathcal{C}(n,\mathbb{M}_{a,d})$ and $\mathcal{C}(m,\mathbb{M}_{a,d})$,
\footnote{
$n+m-2a$ on the right side to get $n+m-2a\in\mathbb{M}_{2a,d}$ from $n,m\in\mathbb{M}_{2a,d}$.
}
we have
\[
\mathcal{C}(n,\mathbb{M}_{a,d}) \cup \mathcal{C}(m,\mathbb{M}_{a,d})
\subset\mathcal{C}(n+m-2a,\mathbb{M}_{a,d}).
\]
\end{corollary}

\begin{remark}
CoPs $\mathcal{C}(n,\mathbb{P})$ with the set of all prime numbers as a base set are important examples for CoPs not admitting embedding. The following example demonstrates this scenario.
\begin{align*}
\mathcal{C}(20,\mathbb{P})&=\lbrace [3],[7],[13],[17]\rbrace\mbox{ but}\\
\mathcal{C}(22,\mathbb{P})&=\lbrace [3],[5],[11],[17],[19]\rbrace.
\end{align*} 
\end{remark}

\begin{proposition}\label{commonfirst}
Let $\mathbb{M}\subseteq\mathbb{N}$ and $\mathcal{C}(n,\mathbb{M})\subset\mathcal{C}(m,\mathbb{M})$ be two CoPs with a common base set $\mathbb{M}$ and $w_0$ and $z_0$ be the weights of the median points of $\mathcal{C}(n,\mathbb{M})$ resp. $\mathcal{C}(m,\mathbb{M})$. If $w_0<z_0$, then the CoPs admit an aligned embedding. On the other hand,  if $w_0>z_0$, then the CoPs admit a reverse aligned embedding.
\end{proposition}

\begin{proof}
Let
\begin{align*}
u_0:&=\min(u\in\Vert\mathcal{C}(n,\mathbb{M})\Vert)\mbox{ and}\\
x_0:&=\min(x\in\Vert\mathcal{C}(m,\mathbb{M})\Vert)
\mbox{ be the least weights of the CoPs and}\\
v_0:&=\max(v\in\Vert\mathcal{C}(n,\mathbb{M})\Vert)\mbox{ and}\\
y_0:&=\max(y\in\Vert\mathcal{C}(m,\mathbb{M})\Vert)
\mbox{ the greatest weights of the CoPs.}
\end{align*}
Because the CoPs are strictly totally ordered, the minimal and maximal points are unique.
We deduce that
\[
w_0:=\frac{u_0+v_0}{2}=\frac{n}{2}\mbox{ and }
z_0:=\frac{x_0+y_0}{2}=\frac{m}{2}
\]
are the weights of the median points of the CoPs. We distinguish and examine three cases
\begin{itemize}
\item[A.)]
$w_0<z_0$,
\item[B.)]
$w_0=z_0$,
\item[C.)]
$w_0>z_0$.
\end{itemize}
Due to $\mathcal{C}(n,\mathbb{M})\subset\mathcal{C}(m,\mathbb{M})$, all points of $\mathcal{C}(n,\mathbb{M})$ must also be points of $\mathcal{C}(m,\mathbb{M})$. Therefore, we must have
\begin{align*}
x_0\leq u_0<v_0\leq y_0.
\end{align*}
Now, we consider case (A.):
\par From $w_0<z_0$, we get $n<m$. That means that $\mathcal{C}(n,\mathbb{M})$ admits an aligned embedding. This includes the case of a common first point ($x_0=u_0$) of both CoPs.\\

We now examine the case (C.):
\par From $w_0>z_0$, we get $n>m$. This means that $\mathcal{C}(n,\mathbb{M})$ admits a reverse aligned embedding. This includes the case of a common last point ($v_0=y_0$) of both CoPs.\\

In case (B.), we would obtain $n=m$. Due to $\mathcal{C}(n,\mathbb{M})\subset\mathcal{C}(m,\mathbb{M})$, there must be at least one real axis $\mathbb{L}_{[s],[t]}\inn\mathcal{C}(m,\mathbb{M})$ which is not a real axis of $\mathcal{C}(n,\mathbb{M})$. This is not possible since they have the same generator. Therefore, the case (B) cannot occur.
\end{proof}
\bigskip

\begin{example}
For an example for reverse aligned embedding, we consider the following CoPs
\begin{align*}
\mathcal{C}(36,\mathbb{P})&=\lbrace [5],\mathbf{[7]},[13],[17],\mathbf{[19]},[23],[29],\mathbf{[31]}\rbrace\mbox{ and }\\
\mathcal{C}(38,\mathbb{P})&=\lbrace [7],[19],[31]\rbrace.
\end{align*}
We observe that $\mathcal{C}(38,\mathbb{P})\subset\mathcal{C}(36,\mathbb{P})$ but $38>36$.
\end{example}

\paragraph*{\textbf{Notation}}
We denote the assignment of an axis $\mathbb{L}_{[x],[y]}$ resp. $\mathbb{L}_{[x]}$ to a CoP $\mathcal{C}(n,\mathbb{M})$ by
\begin{align*}
&\mathbb{L}_{[x],[y]}\inn\mathcal{C}(n,\mathbb{M})
\mbox{ which means }
[x],[y] \in \mathcal{C}(n,\mathbb{M}) \mbox{ and } x+y=n\mbox{ resp.}\\
&\mathbb{L}_{[x]}\inn\mathcal{C}(n,\mathbb{M})
\mbox{ which means }
[x]\in \mathcal{C}(n,\mathbb{M}) \mbox{ and } 2x=n
\end{align*}
and the number of real axes of a CoP by
\[
\nu(n,\mathbb{M}):=\#\lbrace\mathbb{L}_{[x],[y]}\inn\mathcal{C}(n,\mathbb{M})\mid x<y\rbrace.
\]
Clearly, we have
\[
\nu(n,\mathbb{M})=\left\lfloor\frac{k}{2}\right\rfloor ,\mbox{ if }
\vert\mathcal{C}(n,\mathbb{M})\vert =k.
\]
\bigskip

\begin{proposition}\label{increaarithmetic}
Let $\mathbb{M}\subset\mathbb{N}$ and $\mathcal{C}(n,\mathbb{M})$ be a CoP that admits an aligned embedding. The function $\nu(n,\mathbb{M})$ is nondecreasing for all $n$ with $\mathcal{C}(n,\mathbb{M})\neq \emptyset$.
\end{proposition}

\begin{proof}
Because CoP $\mathcal{C}(n,\mathbb{M})$ admits aligned embedding, we get (without loss of generality)
\begin{align*}
\mathcal{C}(n,\mathbb{M})&\subset\mathcal{C}(m,\mathbb{M})
\mbox{ for }n<m \mbox{ and hence}\\
\vert\mathcal{C}(n,\mathbb{M})\vert &< \vert\mathcal{C}(m,\mathbb{M})\vert
\mbox{ and therefore}\\
\nu(n,\mathbb{M})&<\nu(m,\mathbb{M}).
\end{align*} 
\end{proof}

Let
\begin{equation}\label{NStern}
\mathbb{N}^*=\lbrace n\in\mathbb{N}\mid n\equiv\pm 1\mood{6}\rbrace.
\end{equation}
It implies that the set $\mathbb{P}^*$ of all primes $\geq 5$ is covered by $\mathbb{N}^*$.\\

\begin{proposition}\label{properties3}
The CoP $\mathcal{C}(n,\mathbb{N}^*)$ admits an aligned embedding for all $n$ with $n\equiv\pm 2\mood{6}$ or $n\equiv 0\mood{6}$.
\end{proposition}

\begin{proof}
If $n\equiv -2\mood{6}$, then the weights of all points  $[x]\in\mathcal{C}(n,\mathbb{N}^*)$ must satisfy $x\equiv -1\mood{6}$. All points of $\mathcal{C}(n,\mathbb{N}^*)$ are points of $\mathcal{C}(n,\mathbb{M}_{5,6})$.
In the other case, if $n\equiv +2\mood{6}$, then $x\equiv +1\mood{6}$. We deduce
\begin{align*}
\mathcal{C}(n,\mathbb{N}^*)=
\begin{cases}
\mathcal{C}(n,\mathbb{M}_{5,6})&\mbox{ if }n\equiv -2\mood{6}\\
\mathcal{C}(n,\mathbb{M}'_{1,6})&\mbox{ if }n\equiv +2\mood{6}
\end{cases}
\end{align*}
where $\mathbb{M}'_{1,6}:=\mathbb{M}_{1,6}\setminus \lbrace 1\rbrace$. Due to Theorem \ref{properties2}, the claim follows from the case $n\equiv\pm 2\mood{6}$.\\

If $n\equiv 0\mood{6}$, then for every real axis $\mathbb{L}_{[x],[y]}$ we get
\[
x\pmod 6 =-y\pmod6.
\]
This means that if $x\in\mathbb{M}_{5,6}$, then $y\in\mathbb{M}'_{1,6}$ and vise-versa. Without loss of generality, we assume that $x\in\mathbb{M}_{5,6}$ and $y\in\mathbb{M}'_{1,6}$ with $x,y\in\mathcal{C}(n,\mathbb{N}^*)$ and $x<y$. This implies 
\begin{align*}
&x+2\in\mathbb{M}'_{1,6}\mbox{ and }y-2\in\mathbb{M}_{5,6}\mbox{ and due to }x+2+y-2=n\\
&\mbox{it implies }x+2,y-2\in\mathcal{C}(n,\mathbb{N}^*)
\mbox{ with }\mathbb{L}_{[x+2],[y-2]}\inn\mathcal{C}(n,\mathbb{N}^*)\\
&\mbox{and we have a chain of weights of }\mathcal{C}(n,\mathbb{N}^*)\\
&x<x+2<y-2<y.
\end{align*}
Also, we get 
\begin{align*}
&\mathbb{L}_{[x],[y-2]}\inn\mathcal{C}(n-2,\mathbb{M}_{5,6})\mbox{ because of }x+y-2=n-2\mbox{ and}\\
&\mathbb{L}_{[x+2],[y]}\inn\mathcal{C}(n+2,\mathbb{M}'_{1,6})\mbox{ due to }x+2+y=n+2.
\end{align*}
Therefore, to each real axis $\mathbb{L}_{[x],[y]}\inn\mathcal{C}(n,\mathbb{N}^*)$ there exists 
\begin{align*}
&\mbox{a second real axis }\mathbb{L}_{[x+2],[y-2]}\inn\mathcal{C}(n,\mathbb{N}^*)
\mbox{ and}\\
&\mbox{a real axis }\mathbb{L}_{[x],[y-2]}\inn\mathcal{C}(n-2,\mathbb{M}_{5,6})\mbox{ and}\\
&\mbox{a real axis }\mathbb{L}_{[x+2],[y]}\inn\mathcal{C}(n+2,\mathbb{M}'_{1,6}).
\end{align*}
Each point of $\mathcal{C}(n,\mathbb{N}^*)$ is a point of either $\mathcal{C}(n-2,\mathbb{M}_{5,6})$ or $\mathcal{C}(n+2,\mathbb{M}'_{1,6})$.\\

Now, we consider a real axis $\mathbb{L}_{[u],[v]}\inn\mathcal{C}(n-2,\mathbb{M}_{5,6})$ with $u<v$. Since $u+v=n-2$, we deduce
\begin{align*}
&\mathbb{L}_{[u+2],[v]}\inn\mathcal{C}(n,\mathbb{N}^*)\mbox{ due to }u+2+v=n\mbox{ and}\\
&\mathbb{L}_{[u],[v+2]}\inn\mathcal{C}(n,\mathbb{N}^*)\mbox{ due to }u+v+2=n\\
&\mbox{and we have a chain of weights of }\mathcal{C}(n,\mathbb{N}^*)\\
&u<u+2<v<v+2.
\end{align*}
For a real axis $\mathbb{L}_{[w],[z]}\inn\mathcal{C}(n,\mathbb{M}'_{1,6})$ with $w<z$, we have with $w+z=n+2$ the following deduction
\begin{align*}
&\mathbb{L}_{[w-2],[z]}\inn\mathcal{C}(n,\mathbb{N}^*)\mbox{ due to }w-2+z=n\mbox{ and}\\
&\mathbb{L}_{[w],[z-2]}\inn\mathcal{C}(n,\mathbb{N}^*)\mbox{ due to }w+z-2=n\\
&\mbox{and a chain of weights of }\mathcal{C}(n,\mathbb{N}^*)\\
&w-2<w<z-2<z.
\end{align*}
If we assume (without loss of generality) that $u<w$, then we get a chain of weights of $\mathcal{C}(n,\mathbb{N}^*)$
\[
u<u+2<w-2<w<z-2<z<v<v+2.
\]
Thus, all points of $\mathcal{C}(n-2,\mathbb{M}_{5,6})$ and $\mathcal{C}(n+2,\mathbb{M}'_{1,6})$ belong to $\mathcal{C}(n,\mathbb{N}^*)$ and there is no point of $\mathcal{C}(n,\mathbb{N}^*)$ which is not a point of either $\mathcal{C}(n-2,\mathbb{M}_{5,6})$ or $\mathcal{C}(n+2,\mathbb{M}'_{1,6})$. Since the CoPs $\mathcal{C}(n-2,\mathbb{M}_{5,6})$ and $\mathcal{C}(n+2,\mathbb{M}'_{1,6})$ are disjunct because $\mathbb{M}_{5,6}$ and $\mathbb{M}'_{1,6}$ are disjunct, we finally get
\[
\mathcal{C}(n,\mathbb{N}^*)=\mathcal{C}(n-2,\mathbb{M}_{5,6})\cup\mathcal{C}(n+2,\mathbb{M}'_{1,6}).
\]
Since $\mathcal{C}(n-2,\mathbb{M}_{5,6})$ and $\mathcal{C}(n+2,\mathbb{M}'_{1,6})$ admit aligned embedding by Theorem \ref{properties2} and are disjunct, the CoP $\mathcal{C}(n,\mathbb{N}^*)$ admits aligned embedding.
\end{proof}

\begin{corollary}\label{C-centerStern}
If $n\equiv \pm 2\mood{6}$, then the CoP $\mathcal{C}(n,\mathbb{N}^*)$ has a center if and only if $\frac{n}{2}\equiv \pm 1\mood{6}$. In the case $n\equiv 0\mood{6}$, the CoP $\mathcal{C}(n,\mathbb{N}^*)$ has no center because each weight $x$ satisfy $\equiv\pm 1\mood{6}$ and therefore 
\[
2x\equiv\pm 2\mood{6}\not\equiv 0\mood{6}.
\]
\end{corollary}

\begin{corollary}\label{C-centerP}
Due to $\mathbb{P}^*\subset\mathbb{N}^*$, the CoP $\mathcal{C}(n,\mathbb{P}^*)$ has a center if and only if $\frac{n}{2}$ is a prime. In the case $n\equiv 0\mood{6}$, there is no center in $\mathcal{C}(n,\mathbb{P}^*)$ with the same justification as in Corollary \ref{C-centerStern}.
\end{corollary}
\bigskip

\begin{theorem}[Fundamental]\label{central}
Let $n,r \in \mathbb{N}$, $\mathbb{M}\subseteq\mathbb{N}$ and $\mathcal{C}(n,\mathbb{M})$ be a nonempty CoP with an axis $\mathbb{L}_{[x],[n-x]}\inn\mathcal{C}(n,\mathbb{M}).$
\footnote{
The axis can also be a degenerate axis with $x=n-x=\frac{n}{2}$ if it exists.}
If $x+r\in\mathbb{M}$, then $\mathcal{C}(n+r,\mathbb{M})$ is a nonempty CoP. 
\end{theorem}

\begin{proof}
Suppose that $\mathbb{L}_{[x],[n-x]}\inn\mathcal{C}(n,\mathbb{M})$. We deduce that $x$ and $n-x$ are members of $\mathbb{M}$ and $x+r\in\mathbb{M}$. This implies that
\[
n+r-(x+r)=n-x\in\mathbb{M}.
\]
Hence, there is an axis $\mathbb{L}_{[x+r],[n+r-(x+r)]}\inn\mathcal{C}(n+r,\mathbb{M})$ and $\mathcal{C}(n+r,\mathbb{M})$ is nonempty.
\end{proof}
\bigskip

\begin{corollary}\label{C-ext_inftyChain}
Let the requirements of Theorem \ref{central} be fulfilled. If the base set $\mathbb{M}$ is an infinite set and there exists a nonempty CoP $\mathcal{C}(n_0,\mathbb{M})$, then there exist infinitely many positive integers $n>n_0$ with nonempty CoPs $\mathcal{C}(n,\mathbb{M})$.
\end{corollary}

\begin{proof}
Let $\mathbb{L}_{[x],[n-x]}$ be an axis of $\mathcal{C}(n_0,\mathbb{M})$. By Theorem \ref{central}, we deduce that $\mathcal{C}(n_0+r_1,\mathbb{M})$ for some $r_1>0$ is nonempty if $x+r_1\in\mathbb{M}$. We can continue this process with some $r_2>0$ to obtain a nonempty CoP $\mathcal{C}(n_0+r_1+r_2,\mathbb{M})$ such that $n_0+r_1+r_2\in\mathbb{M}$. Since the base set is infinite, this process can be repeated indefinitely. 
\end{proof}

\subsection{The rotation and dilation of circles of partition}

In this section, we introduce the notion of the \emph{Rotation} and \emph{Dilation} of CoPs produced by a given generator.

\begin{definition}
Let $\mathbb{M}\subseteq \mathbb{N}$ with $n\in \mathbb{N}$ and $\mathcal{C}(n,\mathbb{M})$ be the CoP generated by $n$. The map 
\begin{align}
\varpi_r:\mathcal{C}(n,\mathbb{M})\longrightarrow \mathcal{C}^r(n,\mathbb{M})\nonumber 
\end{align}
is the $r^{th}$ level rotation of the CoP $\mathcal{C}(n,\mathbb{M})$ with 
\begin{align*}
\mathcal{C}^r(n,\mathbb{M}):=&\lbrace [k]\in \mathcal{C}(n,\mathbb{M})\mid [x]\in \mathcal{C}(n,\mathbb{M}),~x+r\equiv k\mood{n},~r\in \mathbb{Z},\\
&\mbox{ if }x+r\equiv 0\pmod n\mbox{ then }k:=(n+r)\pmod n
\rbrace.
\end{align*}
If the sign is positive, then we say that the $r^{th}$ level rotation is clockwise.  Otherwise, it is an anti-clockwise $r^{th}$ level rotation for $r\neq 0$. However, if we take $r=0$, then the rotation is trivial and the real axes joining points on the CoP remain stable. The result of a rotation must not be a CoP. Due to condition $[k]\in \mathcal{C}(n,\mathbb{M})$, it is possible that the target set is empty. In this case, we say that the $r^{th}$ level rotation does not exist.
\end{definition}

\begin{theorem}\label{T_rotation-C(n)}
The CoP $\mathcal{C}(n)$ remains invariant under the $r^{th}$ level rotation $\varpi_r$. That is
\[
\varpi_r: \mathcal{C}(n)\longrightarrow\mathcal{C}(n).
\]
\end{theorem}
\begin{proof}
The set of weights of the images of $\mathcal{C}(n)$ under rotation is
\footnote{We denote by $\lbrace a,b,\ldots,z\rbrace_n$ the set
$\lbrace a\pmod n,b\pmod n,\ldots,z\pmod n\rbrace$.
}
\[
\Vert\mathcal{C}^r(n)\Vert=\lbrace r+1,r+2,\ldots,r+n-1\rbrace_n.
\]
The missing value is $(r+n-k)_n$ if $r+n-k\equiv 0\pmod n$. Therefore
\[
k=(n+r)\pmod n.
\]
\end{proof}
If $-n<r<n$, then we get 
\[
k=
\begin{cases}
r &\mbox{ if } r>0\\
n-\vert r\vert &\mbox{ if } r<0.
\end{cases}
\]

\begin{example}
We let $n=8,r=2$ and get 
$$
\Vert\mathcal{C}(8)\Vert=\lbrace 1,2,3,4,5,6,7\rbrace.
$$ 
The critical point is $[6]$ because $6+2\equiv 0\pmod{8}$. The set of weights of the images of all points except $[6]$ is 
$$
\lbrace 3,4,5,6,7,-,1\rbrace
$$ 
The weight $2$ is not in the image set of weights of points under rotation. For the image of the point $[6]$, we get $[(8+2)\pmod 8]=[2]$. We deduce the corresponding weight set of the images of points under rotation 
$$
\Vert\varpi_{3}\left(\mathcal{C}(8)\right)\Vert=\lbrace 3,4,5,6,7,2,1\rbrace=\Vert\mathcal{C}(8)\Vert.
$$ 
In the case $n=8,r=-2$, the point $[2]$ is the \emph{critical} point because $2-2\equiv 0\pmod{8}$. We get the set of weights of images of all points except the point $[2]$
$$
\lbrace 7,-,1,2,3,4,5\rbrace.
$$ 
The weight set of the images of points under rotation does not include $6$. For the point $[2]$, we set $[(8-2)\pmod 8]=[6]$
and get as a target set
$$
\Vert\varpi_{3}\left(\mathcal{C}(8)\right)\Vert=\lbrace 7,6,1,2,3,4,5\rbrace=\Vert\mathcal{C}(8)\Vert.
$$
\end{example}
\bigskip

\begin{proposition}\label{P_not-rotation-Mad}
Let $\mathcal{C}(n,\mathbb{M}_{a,d})$ be a CoP defined as in (\ref{Mad}). The $r^{th}$ level rotation does not exist for $r\equiv c\mood{d}$ with $0<c<d$ and $c\not\equiv 2a\mood{d}$.
\end{proposition}

\begin{proof}
Without loss of generality, we let $c\leq n$. We observe that $[n-a-kd]$ is a point of $\mathcal{C}(n,\mathbb{M}_{a,d})$ for $k=0(1)\frac{n-2a}{d}$
\footnote{Due to $n\in\mathbb{M}_{2a,d}$, $n$ is a positive integer.}
Applying the rotation $\varpi_r$, the weights are transformed to
\begin{align*} 
(n-a-kd+c)\pmod n&=(c-a-kd)\pmod n \mbox{ and because of }c\leq n\\
&=c-a-kd \\
&\equiv (c-a)\pmod d\mbox{ and because of }c\not\equiv 2a\mood{d}\\
&\not\equiv a\pmod d.
\end{align*}
Hence, all rotated points of $\mathcal{C}(n,\mathbb{M}_{a,d})$ are not points of $\mathcal{C}(n,\mathbb{M}_{a,d})$. Therefore, the target set of the rotation is \emph{empty}.
\end{proof}

\begin{proposition}\label{P_ext-rotation-Mad}
Let $\mathcal{C}(n,\mathbb{M}_{a,d})$ be a CoP defined as in (\ref{Mad}). The CoP $\mathcal{C}(n,\mathbb{M}_{a,d})$ remains invariant under the $r^{th}$ level rotation $\varpi_r$ provided that $d=2a$ and $r\equiv 0\mood{d}$.
\end{proposition}

\begin{proof}
We deduce $n\equiv 2a\mood{d}$. Under assumption $d=2a$, we get $n\equiv 0\mood{d}$. Now, let $(x+r)\pmod n=c$ be the weight of a rotated point $[x]$. The following congruence condition holds:
\begin{align*}
x+r&\equiv c\pmod n\mbox{ and because }n\equiv 0\pmod d\\
&\equiv c\pmod d.
\end{align*}
On the other hand, the congruence conditions $x\equiv a\pmod d$ and $r\equiv 0\pmod d$ imply
\begin{align*}
x+r&\equiv a\mood{d}.
\end{align*}
Hence, $a=c$ and $x+r\equiv a\mood{d}$. Therefore, all the images of the points of CoP $\mathcal{C}(n,\mathbb{M}_{a,d})$ are members of $\mathbb{M}_{a,d}$ and less than $n$. The image points of the $r^{th}$ level rotation of CoP $\mathcal{C}(n,\mathbb{M}_{a,d})$ are again points of CoP $\mathcal{C}(n,\mathbb{M}_{a,d})$. This proves the claim that CoPs of the form $\mathcal{C}(n,\mathbb{M}_{a,d})$ remain invariant under some $r^{th}$ level rotation under special conditions.
\end{proof}
\bigskip

\begin{example}
We let $n=24$, $a=2$, $d=4$, $r=4$ and construct the weight set of the corresponding CoP 
$$
\Vert\mathcal{C}(24,\mathbb{M}_{2,4})\Vert=\lbrace 2,6,10,14,18,22\rbrace.
$$
We deduce
$$
\Vert\varpi_4\left(\mathcal{C}(24,\mathbb{M}_{2,4})\right)\Vert=\lbrace 6,10,14,18,22,2\rbrace.
$$
\end{example}

\begin{corollary}
For conditions espoused in Proposition \ref{P_not-rotation-Mad} and Proposition \ref{P_ext-rotation-Mad}, the $r^{th}$ level rotation of a CoP $\mathcal{C}(n,\mathbb{M}_{a,d})$ is a real (proper) subset of $\mathcal{C}(n,\mathbb{M}_{a,d})$. 
\end{corollary}

\begin{definition}
Let $\mathbb{M}\subseteq \mathbb{N}$ with $n\in \mathbb{N}$ and  $\mathcal{C}(n,\mathbb{M})$ be the CoP generated by $n$. The map 
\begin{align}
\delta_r:\mathcal{C}(n,\mathbb{M})\longrightarrow \mathcal{C}_r(n,\mathbb{M})\nonumber
\end{align}
is the $r^{th}$ scale dilation of CoP $\mathcal{C}(n,\mathbb{M})$ with 
\begin{align}
\mathcal{C}_r(n,\mathbb{M}):=\left \{[x]\in \mathcal{C}(n+r,\mathbb{M})\mid~r\in \mathbb{Z},~n+r>1\right \}.\nonumber
\end{align}
If the sign is positive, then we say that the $r^{th}$ scale dilation is an expansion. Otherwise, it is an $r^{th}$ scale compression for $r\neq 0$. However, if we take $r=0$, then the dilation is a trivial dilation and the CoP remains invariant under the dilation.
\end{definition}
\bigskip

\begin{remark}\label{remC(n)}
We note that if the base set is taken to be the set of natural numbers $\mathbb{N}$, then the image set of a dilation collapses to the following 
\begin{align}\label{C(n)}
\delta_r\big(\mathcal{C}(n)\big)&:=\mathcal{C}_r(n)\nonumber\\
&=\left \{[x]\mid x\in\mathbb{N}_{n+r-1},~r\in \mathbb{Z},~n+r>1\right \}\nonumber \\
&=\mathcal{C}(n+r).
\end{align}
Furthermore, in case $r<0$, some points in $\mathcal{C}(n)$ have the same image, while in case $r>0$ some points in $\mathcal{C}(n)$ have more than one image.
\end{remark}
\bigskip

The dilations at a scale between CoPs have the natural tendency to translate the generator of the source CoP by the size of the scale of the dilation. However, it may be difficult to define dilation on individual points in a CoP. Any perceived dilation map could manifestly work on a typical CoP, but it may proved handicapped for some other CoPs. In the sense that some points may poke outside the target CoP under this fixed dilation. In light of this anomaly, we ask the following questions

\begin{question}\label{question1}
Let $\mathbb{M}\subseteq \mathbb{N}$. Does there exists a well-defined dilation 
\begin{align}
\delta_r:\mathcal{C}(n,\mathbb{M})\longrightarrow \mathcal{C}(m,\mathbb{M})\nonumber
\end{align} 
on each $[x]\in \mathcal{C}(n,\mathbb{M})$ for all CoPs?
\end{question}
\bigskip

The question \ref{question1} asks whether there exists a fixed map that assigns each point in a typical CoP to its target CoP in a sufficiently uniform way. That is, the map we seek should avoid the subtleties as espoused in our earlier discussion.

\begin{theorem}\label{existence}
Let $n,m\in \mathbb{N}$ with $\mathbb{M}\subseteq\mathbb{N}$, and $\mathcal{C}(n,\mathbb{M})$ be a CoP that admits an aligned embedding. There exists some dilation $\delta_r$ such that 
\begin{align*}
\delta_r:\mathcal{C}(n,\mathbb{M})\longrightarrow \mathcal{C}(m,\mathbb{M}).
\end{align*} 
\end{theorem}

\begin{proof}
It is evident that for $m=n$ the trivial dilation $\delta_0$ meets the claim. For the case $m\neq n$, we break the proof into several cases. The case $r$ is positive and the case is negative. Let $\delta_r$ be any dilation for $r>0$ and suppose that for any two CoP $\mathcal{C}(n,\mathbb{M})$ and $\mathcal{C}(m,\mathbb{M})$ with $\mathcal{C}(m,\mathbb{M})\subset \mathcal{C}(n,\mathbb{M})$ there exists no dilation associating them. By the property that the CoPs admit embedding, exactly one of the following embeddings holds:
\begin{align*}
\delta_r\big(\mathcal{C}(m,\mathbb{M})\big)\subset \mathcal{C}(n,\mathbb{M})  \mbox{ or } \mathcal{C}(n,\mathbb{M})\subset \delta_r\big(\mathcal{C}(m,\mathbb{M})\big).
\end{align*}
We analyze each of these sub-cases. First, we assume that $\delta_r\big(\mathcal{C}(m,\mathbb{M})\big)\subset \mathcal{C}(n,\mathbb{M})$. This implies that there exists some CoP $\mathcal{C}(s,\mathbb{M})$ with  $\delta_r(\mathcal{C}(m,\mathbb{M}))\subseteq \mathcal{C}(s,\mathbb{M})$ such that $\mathcal{C}(s,\mathbb{M})\subset \mathcal{C}(n,\mathbb{M})$. Since there exists no dilation between CoPs, the following proper embedding must necessarily hold 
\begin{align*}
\delta_r\big(\mathcal{C}(m,\mathbb{M})\big)\subset \mathcal{C}(s,\mathbb{M})\subset \mathcal{C}(n,\mathbb{M}).
\end{align*}
Again, there exists some CoP $\mathcal{C}(t,\mathbb{M})$ with $\delta_r(\mathcal{C}(m,\mathbb{M}))\subseteq \mathcal{C}(t,\mathbb{M})$ such that $\mathcal{C}(t,\mathbb{M})\subset \mathcal{C}(s,\mathbb{M})$. Under the underlying assumption that there exists no dilation between CoPs, we obtain the following proper embedding 
\begin{align*}
\delta_r(\mathcal{C}(m,\mathbb{M}))\subset \mathcal{C}(t,\mathbb{M})\subset \mathcal{C}(s,\mathbb{M})\subset \mathcal{C}(n,\mathbb{M}).
\end{align*}
Repeating the argument in this manner, we obtain the following infinite descending chains of covers of the smallest CoP
\begin{align*}
\mathcal{C}(m+r,\mathbb{M}):=\delta_r\big(\mathcal{C}(m,\mathbb{M})\big)\subset \cdots \subset \mathcal{C}(t,\mathbb{M}) \subset \mathcal{C}(s,\mathbb{M})\subset \mathcal{C}(n,\mathbb{M}).
\end{align*}
Because the CoPs admit aligned embedding, we obtain the infinite descending sequence of positive integers towards the generator $m+r$ of the last CoP  
\begin{align*}
n>s>t>\cdots>\cdots>m+r.
\end{align*}
This is impossible, thereby ending the proof of the first sub-case. We now turn to the case $\mathcal{C}(n,\mathbb{M})\subset \delta_r(\mathcal{C}(m,\mathbb{M}))$. In a similar fashion, there exists some CoP $\mathcal{C}(t,\mathbb{M})$ with $\mathcal{C}(t,\mathbb{M})\subseteq \delta_r(\mathcal{C}(m,\mathbb{M}))$ such that $\mathcal{C}(n,\mathbb{M})\subset \mathcal{C}(t,\mathbb{M})$. Under the assumption that there exists no dilation between CoP, we have the following embedding
\begin{align*}
\mathcal{C}(n,\mathbb{M})\subset \mathcal{C}(t,\mathbb{M})\subset \delta_r(\mathcal{C}(m,\mathbb{M})).
\end{align*}
Again, there exists some CoP $\mathcal{C}(s,\mathbb{M})$ with $\mathcal{C}(s,\mathbb{M})\subseteq \delta_r(\mathcal{C}(m,\mathbb{M}))$ such that $\mathcal{C}(t,\mathbb{M})\subset \mathcal{C}(s,\mathbb{M})$. Under the assumption that there exists no dilation between CoP, we have the following embedding 
\begin{align*}
\mathcal{C}(n,\mathbb{M})\subset \mathcal{C}(t,\mathbb{M})\subset \mathcal{C}(s,\mathbb{M})\subset \delta_r(\mathcal{C}(m,\mathbb{M})).
\end{align*}
By repeating this argument indefinitely, we obtain the following infinite sequence of embeddings 
\begin{align*}
\mathcal{C}(n,\mathbb{M})\subset \mathcal{C}(t,\mathbb{M})\subset \mathcal{C}(s,\mathbb{M})\cdots \subset \delta_r(\mathcal{C}(m,\mathbb{M})):=\mathcal{C}(m+r,\mathbb{M}).
\end{align*}
By the admitting aligned embedding of CoPs, we obtain an infinite ascending sequence of positive integers towards the generator of the last CoP in the chain 
\begin{align*}
n<t<s<\cdots<m+r.
\end{align*}
This is impossible since we cannot have positive integers indefinitely approaching a fixed positive integer. This completes the proof for the case $r>0$. We now turn to the case $r<0$ for any two CoP $\mathcal{C}(m,\mathbb{M}), \mathcal{C}(n,\mathbb{M})$ with $\mathcal{C}(n,\mathbb{M})\subset \mathcal{C}(m,\mathbb{M})$. Under the main assumption exactly one of the following embeddings must hold 
\begin{align*}
\delta_r\big(\mathcal{C}(m,\mathbb{M})\big)\subset \mathcal{C}(n,\mathbb{M}) \mbox{ or } \mathcal{C}(n,\mathbb{M})\subset \delta_r\big(\mathcal{C}(m,\mathbb{M})\big).
\end{align*}
A similar analysis can be performed for each of the above cases.
\end{proof}
\bigskip

\begin{corollary}
By Theorem \ref{properties}, the CoP $\mathcal{C}(n)$ admits an aligned embedding, and there is the dilation $\delta_1:\mathcal{C}(n)\longrightarrow \mathcal{C}(n+1)$ with 
\begin{align}\label{E_Cn-Cn+1}
\delta_1([x]):=
\begin{cases}
[x] \quad &\mbox{for }~1\leq x\leq n-1 \\
[n] \quad &\mbox{additional for}~x=1
\end{cases}
\end{align} 
that can produce an infinite ascending chain of CoPs 
\begin{align}
\mathcal{C}(n)\subset \mathcal{C}(n+1)\subset \mathcal{C}(n+2)\subset \cdots. \nonumber
\end{align}
\end{corollary}
\bigskip

We observe that the assignment of $[n]$ as an image of $[1]$ is not the only possibility. We may also assign $[n]$ as the image of $[2]\ldots[n-1]$. In all cases, we would have a correct point-to-point mapping. Hence a subset of the cross set $\mathcal{C}(n)\times\mathcal{C}(n+1)$ such that:
\begin{itemize}
\item
for each point of $\mathcal{C}(n)$ there is at least one image point of $\mathcal{C}(n+1)$ and
\bigskip

\item
for each image point of $\mathcal{C}(n+1)$ there is only one preimage point of $\mathcal{C}(n)$
\end{itemize}
is not a well-defined pointwise definition of the map $\mathcal{C}(n)\longrightarrow\mathcal{C}(n+1)$ because there are several such subsets.

\begin{corollary}
By Theorem \ref{properties2}, the CoP $\mathcal{C}(n,\mathbb{M}_{a,d})$ admits an aligned embedding, and there is
the dilation $\delta_d:\mathcal{C}(n,\mathbb{M}_{a,d})\longrightarrow \mathcal{C}(n+d,\mathbb{M}_{a,d})$ with 
\[
\delta_d([x]):=\begin{cases}
[x] \quad &\mbox{for } a\leq x\leq n-a \\
[n-a+d] \quad &\mbox{additional for}~x=a
\end{cases}
\] 
that can generate an infinite ascending chain of CoPs 
\begin{align}
\mathcal{C}(n,\mathbb{M}_{a,d})\subset \mathcal{C}(n+d,\mathbb{M}_{a,d})\subset \mathcal{C}(n+2d,\mathbb{M}_{a,d})\subset \cdots. \nonumber
\end{align}
\end{corollary}

\subsection{The special maps of circles of partition}

In this section, we introduce and study several special maps of circles of partitions.

\begin{definition}\label{flip}
Let $\mathbb{M}\subseteq \mathbb{N}$ and $\mathcal{C}(n,\mathbb{M})\neq \emptyset$ be a CoP containing the axis $\mathbb{L}_{[a],[b]}$. By the \emph{flipping} of CoP $\mathcal{C}(n,\mathbb{M})$ along the \emph{flipping axis} $\mathbb{L}_{[a],[b]}$, we mean the map 
\begin{align*}
\vartheta_{[a],[b]}:\mathcal{C}(n,\mathbb{M})\longrightarrow \mathcal{C}(m,\mathbb{M})
\end{align*}
with $\vartheta_{[a],[b]}([a])=[a]$ and $\vartheta_{[a],[b]}([b])=[b]$ such that for any two $[x],[y]\in \mathcal{C}(n,\mathbb{M})$ with $[x],[y]\neq [a],[b]$, we have
\begin{align*}
\Vert\vartheta_{[a],[b]}([x])\Vert+\Vert\vartheta_{[a],[b]}([y])\Vert \neq n 
\end{align*}
A flipping axis can also be a degenerate axis $\mathbb{L}_{[a]}$.
We say that CoP $\mathcal{C}(n,\mathbb{M})$ is \emph{susceptible} to flipping if there exists such a map.
\end{definition}
\bigskip

\begin{example}
We let $\mathbb{M}=\mathbb{P}$ and $n=20$. The CoP $\mathcal{C}(20,\mathbb{P})$ is the set 
$$
\lbrace [3],[7],[13],[17]\rbrace
$$ 
with two axes $\mathbb{L}_{[3],[17]}$ and $\mathbb{L}_{[7],[13]}$. The map
\[
\vartheta_{[3],[17]}: \mathcal{C}(20,\mathbb{P})\longrightarrow
\mathcal{C}(22,\mathbb{P})
\]
with $\mathcal{C}(22,\mathbb{P})=\lbrace [3],[5],[11],[17],[19]\rbrace$
is a flipping of $\mathcal{C}(20,\mathbb{P})$ along the axis $\mathbb{L}_{[3],[17]}$ if f.i.
\begin{align*}
\vartheta_{[3],[17]}([3])&=[3]\\
\vartheta_{[3],[17]}([7])&=[5]\\
\vartheta_{[3],[17]}([13])&=[11]\mbox{ and }[19]\\
\vartheta_{[3],[17]}([17])&=[17].
\end{align*}
Hence, we get $\Vert[5]\Vert+\Vert[11]\Vert=16\neq 20$ or  $\Vert[5]\Vert+\Vert[19]\Vert=24\neq 20$. However, there are
no axis points of $\mathcal{C}(22,\mathbb{P})$ that are also points of $\mathcal{C}(20,\mathbb{P})$. Hence, there is no flipping from $\mathcal{C}(22,\mathbb{P})$ to $\mathcal{C}(20,\mathbb{P})$ along an axis of $\mathcal{C}(22,\mathbb{P})$.
\end{example}
\bigskip

\begin{proposition}\label{flip_Mad}
Let $\mathbb{M}_{a,d}$ be defined as in (\ref{Mad}) with $0<a\leq d$. The CoP $\mathcal{C}(n,\mathbb{M}_{a,d})$ is susceptible to flipping if and only if $n>m$. 
\end{proposition}

\begin{proof}
We observe that to get $\mathcal{C}(n,\mathbb{M}_{a,d})\neq\emptyset$, we must have $n\in \mathbb{M}_{2a,d}$. This implies $n-a\in \mathbb{M}_{a,d}$. The same is valid for $\mathcal{C}(m,\mathbb{M}_{a,d})$. We assume that $n>m$. By Corollary \ref{col-prop2}, we get
\[
\mathcal{C}(n,\mathbb{M}_{a,d}) \supset \mathcal{C}(m,\mathbb{M}_{a,d}).
\]
Due to $n\in\mathbb{M}_{2a,d}$, we have $\frac{n-2a}{d}\in\mathbb{N}$. The weights of $\mathcal{C}(n,\mathbb{M}_{a,d})$ are
\begin{align*}
\Vert\mathcal{C}(n,\mathbb{M}_{a,d})\Vert
=\left\lbrace a+k\cdot d\mid k=0,1,2,\ldots,\frac{n-2a}{d}\right\rbrace.
\end{align*}
Hence, CoP $\mathcal{C}(n,\mathbb{M}_{a,d})$ has 
\[
l_n=\frac{n-2a}{d}+1 \mbox{ members.}
\] 
This is in accordance with the general counting function for CoPs:
\begin{align}
||\mathcal{C}(n,\mathbb{M}_{a,d})||
&=1+\sum \limits_{\substack{1\leq x \leq n-a\\x\equiv a\pmod d}}1\nonumber \\&=1+\frac{n-2a}{d}.\nonumber
\end{align}
The addition of $1$ is required because the counting starts with $0$. Now we must distinguish and examine two cases
\begin{itemize}
\item[rC:]
The CoP $\mathcal{C}(n,\mathbb{M}_{a,d})$ has a real center.
\bigskip

\item[dC:]
The CoP $\mathcal{C}(n,\mathbb{M}_{a,d})$ has a deleted center.
\end{itemize} 
In the case rC, we have that $l_n$ is odd and $l_n$ is even in the other case. Now, we choose in the case rC the degenerate axis $\mathbb{L}_{[u]}$ of CoP $\mathcal{C}(n,\mathbb{M}_{a,d})$ as the flipping axis and in the case dC those that are closest to the center of CoP. The weights of $[u],[v]$ are $u=v=\frac{n}{2}$ for the rC case and $u=\frac{n-d}{2}$, $v=\frac{n+d}{2}$ in the other case. To satisfy the conditions of the definition \ref{flip}
\begin{align*}
\vartheta_{[u],[v]}([u])=[u]\mbox{ and }
\vartheta_{[u],[v]}([v])=[v]
\end{align*}
the last point of $\mathcal{C}(m,\mathbb{M}_{a,d})$ should be $[v]$. Due to Corollary \ref{partner}, we get for $m$ as the sum of the weights of the first and the last member of CoP $\mathcal{C}(m,\mathbb{M}_{a,d})$ 
\begin{align}\label{E-m_rC}
m=
\begin{cases}
a+\dfrac{n}{2}&\mbox{ for rC}\\
a+\dfrac{n+d}{2}&\mbox{ for dC}.
\end{cases}
\end{align}
We get the number of members of $\mathcal{C}(m,\mathbb{M}_{a,d})$
\begin{align*}
l_m-1 &:= \sum \limits_{\substack{1\leq x \leq m-a\\x\equiv a\pmod d}}1=\frac{m-2a}{d}\\
&=
\begin{cases}
\dfrac{a+\frac{n}{2}-2a}{d}=\dfrac{n-2a}{2d}=\dfrac{l_n-1}{2}&\mbox{ for rC}\\
\dfrac{a+\frac{n+d}{2}-2a}{d}=\dfrac{n-2a}{2d}+\dfrac{1}{2}=\dfrac{l_n}{2}&\mbox{ for dC}.
\end{cases}
\end{align*}
Hence, we obtain for both cases
\begin{align*}
\l_m=\left\lfloor\frac{l_n}{2}\right\rfloor+1.
\end{align*}
All of these fulfill the \textbf{following map}
\begin{align*}
\vartheta_{[u],[v]}(x)&=a+k(x)\cdot d \mbox{ with}\\
\frac{x-a}{d}&\equiv k(x)\mood{l_m}.
\end{align*}
This map assigns each point of $\mathcal{C}(n,\mathbb{M}_{a,d})$ to a point of $\mathcal{C}(m,\mathbb{M}_{a,d})$. The heaviest point of CoP $\mathcal{C}(m,\mathbb{M}_{a,d})$ is $[m-a]$. In the case rC, the flipping axis is $\mathbb{L}_{[u]}$ with $u=\frac{n}{2}$. We get (with \ref{E-m_rC})
\begin{align*}
\left\Vert\vartheta_{[v],[v]}\left(\left[\frac{n}{2}\right]\right)\right\Vert=m-a&=\frac{n}{2}.
\end{align*}
Hence, the requirement $\Vert\vartheta_{[v],[v]}([v])\Vert=u=v=\frac{n}{2}$ is fulfilled. \\

In the case dC, we get (with \ref{E-m_rC})
\begin{align*}
\left\Vert\vartheta_{[u],[v]}\left(\left[\frac{n+d}{2}\right]\right)\right\Vert=m-a&=\frac{n+d}{2}=v.
\end{align*}
Therefore, $u=v-d=\frac{n-d}{2}$. For each two points $[x],[y]\in\mathcal{C}(n,\mathbb{M}_{a,d})$ with $[x],[y]\neq [u],[v]$, we have
\[
\Vert\vartheta_{[u],[v]}([x])\Vert+\Vert\vartheta_{[u],[v]}([y])\Vert<n
\]
because $\vartheta_{[u],[v]}([u])=[u]$ and $\vartheta_{[u],[v]}([v])=[v]$ are the two heaviest points of $\mathcal{C}(m,\mathbb{M}_{a,d})$ in the dc case, respectively the heaviest point of $\mathcal{C}(m,\mathbb{M}_{a,d})$ in rC case with the sum of weights of the two heaviest points $\leq n$. The weight sum of all others is lesser. The first part of the claim is proven. On the other hand, if $n\leq m$, then the source CoP is a subset of the target CoP. All the axes points of $\mathcal{C}(n,\mathbb{M}_{a,d})$ are identically mapped to $\mathcal{C}(m,\mathbb{M}_{a,d})$. We deduce $\vartheta_{[u],[v]}([x])$ and $\vartheta_{[u],[v]}([y])$ from any axis $\mathbb{L}_{[x],[y]}$ of $\mathcal{C}(n,\mathbb{M}_{a,d}$ with
\[
\Vert\vartheta_{[u],[v]}([x])\Vert + \Vert\vartheta_{[u],[v]}([y])\Vert = n.
\]
This is contradicts the requirements of the claim.
\end{proof}
\bigskip

\begin{remark}
We note that due to $\mathbb{M}_{1,1}=\mathbb{N}$, this statement also holds for each CoP $\mathcal{C}(n)$.
\end{remark}
\bigskip

\begin{proposition}
The chosen axis with points closest to the center of the CoP $\mathcal{C}(n,\mathbb{M}_{a,d})$--respectively the degenerate axis in the case of an existing center--is the only one for flipping along an axis in the case $\mathbb{M}=\mathbb{M}_{a,d}$.
\end{proposition}

\begin{proof}
For all axes $\mathbb{L}_{[x],[y]}\inn\mathcal{C}(n,\mathbb{M}_{a,d})$, we have 
\footnote{Without loss of generality, we assume $x\leq y$ for all axes $\mathbb{L}_{[x],[y]}$. 
In the case of an existing center, we have $\mathbb{L}_{[x]}=\mathbb{L}_{[x],[x]}$.}
\[
x \leq \frac{n}{2} \leq y.
\]
Therefore, there is no axis $\mathbb{L}_{[x],[y]}$ with $y<\frac{n}{2}$. For the chosen axis $\mathbb{L}_{[u],[v]}$ with axes points closest to the center of $\mathcal{C}(n,\mathbb{M}_{a,d})$, we have
\[
\frac{n-d}{2}\leq\Vert[u]\Vert\leq\Vert[v]\Vert\leq\frac{n+d}{2}.
\]
The remaining possibility are axes $\mathbb{L}_{[w],[z]}$ with $\Vert[w]\Vert<\frac{n-d}{2}$ and axis partner with $\Vert[z]\Vert>\frac{n+d}{2}$. Therefore, there is at least one axis $\mathbb{L}_{[x],[y]}$ with $w<x\leq y<z$ and $x+y=n$. This violates the requirements of flipping along the axes $\mathbb{L}_{[w],[z]}$. Hence, only the axis $\mathbb{L}_{[u],[v]}$ (resp. $\mathbb{L}_{[u]}$) with
\begin{itemize}
\item[rC:]
$\Vert[u]\Vert=\dfrac{n}{2}$
\item[dC:]
$\Vert[u]\Vert=\dfrac{n-d}{2},\Vert[v]\Vert=\dfrac{n+d}{2}$
\end{itemize} 
satisfies the requirements of a flipping axis.  
\end{proof}

It is quite suggestive that the notion of flipping of CoPs under $\mathbb{M}=\mathbb{M}_{a,d}$ can be thought of  as the process of slicing a circle into two equal half and overturning one half to lie perfectly on top of the other half and forming a geometric structure that resembles the semi-circle.

\begin{example}
We choose $a=2, d=4$ and hence $\mathbb{M}=\mathbb{M}_{2,4}$. With the generator $n=28$, we have
\begin{align*}
\Vert\mathcal{C}(28,\mathbb{M}_{2,4})\Vert
&=\lbrace 2,6,10,14,18,22,26\rbrace,\\
l_n&=\frac{28-2\cdot 2}{4}+1=7,\\
l_m&=\left\lfloor\frac{7}{2}\right\rfloor+1=4\mbox{ and}\\
m&=2+\frac{28}{2}=16
\end{align*}
with the flipping axis $\mathbb{L}_{[14]}$.
Hence
\begin{align*}
\Vert\vartheta_{[14],[14]}\left(\mathcal{C}(28,\mathbb{M}_{2,4})\right)\Vert
&=\Vert\mathcal{C}(16,\mathbb{M}_{2,4})\Vert\\
&=\lbrace 2,6,10,14\rbrace.
\end{align*}
Here, the weight sums of any two members of $\lbrace [2],[6],[10],[14]\rbrace\setminus\lbrace 14\rbrace$ are less than $28$. If we take $\mathbb{L}_{[6],[22]}$ as flipping axis, then we obtain as target set
\[
\mathcal{C}(24,\mathbb{M}_{2,4})=\lbrace [2],[6],[10],[14],[18],[22]\rbrace.
\]
Removing $\lbrace[6],[22]\rbrace$ from the image set, we get one weight sum contradicting requirements:
\[
10+18=28.
\]
\end{example}
\bigskip

Now, we introduce and study the concept of \emph{filtration} of the CoPs. We begin with the concept of the \emph{filtration} along an axis.

\begin{definition}\label{filtration}
Let $\mathbb{M}\subseteq \mathbb{N}$ with the corresponding CoP $\mathcal{C}(n,\mathbb{M})$ containing the axis $\mathbb{L}_{[x],[y]}$. By the \emph{filtration} of CoP $\mathcal{C}(n,\mathbb{M})$ along the \emph{filtration axis} $\mathbb{L}_{[x],[y]}$, we mean the map
\begin{align*}
\Phi_{[x],[y]}:\mathcal{C}(n,\mathbb{M})\longrightarrow \mathcal{C}(m,\mathbb{M})
\end{align*}
such that $[x],[y]\not\in\mathcal{C}(m,\mathbb{M})$ for some $m\in \mathbb{N}\setminus \{1\}$ and there exists the \emph{co-axis} $\mathbb{L}_{[u],[v]}$ of $\mathcal{C}(n,\mathbb{M})$ so that $\mathbb{L}_{[u],[a]}$ and $\mathbb{L}_{[v],[b]}$ are axes of $\mathcal{C}(m,\mathbb{M})$ for some $[a],[b]\in \mathbb{M}$. We say that CoP $\mathcal{C}(n,\mathbb{M})$ is \emph{susceptible to filtration} if there exists such a map. The filtration axis can also be a degenerate axis.
\end{definition}

\begin{example}
We let $\mathbb{M}=\mathbb{P}$ and $n=20$. The map 
\[
\Phi_{[7],[13]}: \mathcal{C}(20,\mathbb{P})\longrightarrow\mathcal{C}(22,\mathbb{P})
\]
is a filtration of $\mathcal{C}(20,\mathbb{P})$ along the filtration axis $\mathbb{L}_{[7],[13]}$ because target CoP 
\[
\mathcal{C}(22,\mathbb{P})=\lbrace [3],[5],[11],[17],[19]\rbrace
\]
contains the axes $\mathbb{L}_{[3],[19]}$ and $\mathbb{L}_{[17],[5]}$, where $\mathbb{L}_{[3],[17]}$ is the co-axis of $\mathcal{C}(20,\mathbb{P})$.  
\end{example}

\begin{example}
We take $\mathbb{M}=\mathbb{P}$ and $n=46$. We get
\begin{align*}
\Vert\mathcal{C}(46,\mathbb{P})\Vert&=\lbrace 3,5,17,23,29,41,43\rbrace
\mbox{ and}\\
\Vert\mathcal{C}(50,\mathbb{P})\Vert&=\lbrace 3,7,13,19,31,37,43,47\rbrace.
\end{align*}
The map
\[
\Phi_{[23]}:\mathcal{C}(46,\mathbb{P})\longrightarrow\mathcal{C}(50,\mathbb{P})
\]
is a filtration of $\mathcal{C}(46,\mathbb{P})$ along the degenerate axis $\mathbb{L}_{[23]}$ because the target CoP contains $\mathbb{L}_{[3],[47]}$ and $\mathbb{L}_{[7],[43]}$, where $\mathbb{L}_{[3],[43]}$ is the co-axis of $\mathcal{C}(46,\mathbb{P})$.
\end{example}

\begin{proposition}\label{nex_filt_Mad}
The CoP $\mathcal{C}(n,\mathbb{M})$ admitting an aligned embedding is \textbf{not} susceptible to filtration along an axis.
\end{proposition}

\begin{proof}
The claim is true if one of the following statements holds:
\begin{enumerate}
\item[(A)]
The CoP $\mathcal{C}(n,\mathbb{M})$ has no filtration axis.
\bigskip

\item[(B)]
The CoP $\mathcal{C}(n,\mathbb{M})$ has no co-axis
\end{enumerate}
We suppose that $n\leq m$. By Theorem \ref{properties2}, we obtain
\begin{align*}
\mathcal{C}(n,\mathbb{M})\subseteq\mathcal{C}(m,\mathbb{M}).
\end{align*}
The images of all axis points of the source CoP are points of the target CoP. Hence, there is no filtration axis (A).\\

Now, we examine the case $m<n$. By Corollary \ref{col-prop2}, we get
\begin{align*}
\mathcal{C}(n,\mathbb{M})\supset\mathcal{C}(m,\mathbb{M}).
\end{align*}
We let $m<\frac{n}{2}$. In this case, the images of the end points of all axes of $\mathcal{C}(n,\mathbb{M})$ do not exist in $\mathcal{C}(m,\mathbb{M})$. Hence, there is no co-axis (B).\\

Finally, we consider the case $\frac{n}{2}\leq m <n$. In this case, the images of the lower weight points of all axes of $\mathcal{C}(n,\mathbb{M})$ are points of $\mathcal{C}(m,\mathbb{M})$. Hence, there is no filtration axis (A). 
\end{proof}
\bigskip

\begin{definition}\label{reduct}
Let $\mathbb{M}\subseteq \mathbb{N}$ with the corresponding CoP $\mathcal{C}(n,\mathbb{M})$ that contains the axis $\mathbb{L}_{[x],[y]}$. By the \emph{reduction} of CoP $\mathcal{C}(n,\mathbb{M})$ in the base set $\mathbb{M}$, we mean the map
\begin{align*}
\phi_{[x],[y]}:\mathcal{C}(n,\mathbb{M})\longrightarrow \mathcal{C}(n,\mathbb{M^*})
\end{align*}
with $\mathbb{M}^*:=\mathbb{M}\setminus\lbrace x,y\rbrace$. We say that CoP $\mathcal{C}(n,\mathbb{M})$ is \emph{susceptible to reduction} if there exists such a map.
\end{definition}
\bigskip

\begin{proposition}
Let $\mathbb{M}_{a,d}$ be defined as in (\ref{Mad}). The CoP $\mathcal{C}(n,\mathbb{M}_{a,d})$ is susceptible to reduction.
\end{proposition}

\begin{proof}
Without loss of generality, we suppose that $x<y$ and take
\begin{align*}
\phi_{[x],[y]}([u])=\left\lbrace
\begin{array}{ll}
~[u] &\mbox{ if }u\neq x \mbox{ and }u\neq y \\
~[u+d] &\mbox{ if }u=x\\
~[u-d] &\mbox{ if }u=y
\end{array}
\right.
\end{align*}
for all points $[u]\in\mathcal{C}(n,\mathbb{M}_{a,d})$. Since all members of $\mathbb{M}_{a,d}$ are arithmetic progressions with a common difference $d$, we find that if $u\in\mathbb{M}_{a,d}$, then $u\pm d\in\mathbb{M}_{a,d}$ and 
\begin{align*}
\Vert\phi_{[x],[y]}([x])\Vert+\Vert\phi_{[x],[y]}([y])\Vert &= x+d+y-d=n
\end{align*}
because $\mathbb{L}_{[x],[y]}$ is an axis of $\mathcal{C}(n,\mathbb{M}_{a,d})$.
\end{proof}

Because $\mathbb{M}_{1,1}=\mathbb{N}$, this proposition holds for $\mathcal{C}(n)$.

\subsection{The stable and unstable points on the circle of partition}

In this section, we introduce the notion of \emph{stability} of a sequence under a given dilation.

\begin{definition}
Let $\Theta(n)$ be a subsequence of $\mathbb{N}_n$ and suppose that the CoP $\mathcal{C}(n,\mathbb{M})\neq \emptyset$.  Let $\mathbb{L}_{[x],[y]}$ be a real axis of the CoP $\mathcal{C}(n,\mathbb{M})$ with $x,y\in \Theta(n)$. We say that the point $[x]\in \mathcal{C}(n,\mathbb{M})$ is \emph{stable} relative to the subsequence $\Theta(n)$ under the $r^{th}$ level rotation $\varpi_r:\mathcal{C}(n,\mathbb{M})\longrightarrow \mathcal{C}(n,\mathbb{M})$ if $||\varpi_r([x])||\in \Theta(n)$ and $\exists z\in \Theta(n)$ such that $\mathbb{L}_{[\varpi_r([x])],[z]}$ is also a real axis of CoP $\mathcal{C}(n,\mathbb{M})$. We say that the subsequence $\Theta(n)$ is \emph{stable} under the $r^{th}$ level rotation $\varpi_r$ if all points in $[x]\in \mathcal{C}(n,\mathbb{M})$ with $x\in \Theta(n)$ are stable.
\end{definition}

\begin{definition}
Let $\Theta(n)$ be a subsequence of $\mathbb{N}_n$ and suppose the CoP $\mathcal{C}(n,\mathbb{M})\neq \emptyset$. Let $\mathbb{L}_{[x],[y]}$ be a real axis of the CoP $\mathcal{C}(n,\mathbb{M})$ with $x,y\in \Theta(n)$. We say that the point $[x]\in \mathcal{C}(n,\mathbb{M})$ is \emph{stable} relative to the subsequence $\Theta(n)$ under the $r^{th}$ scale dilation $\delta_r:\mathcal{C}(n,\mathbb{M})\longrightarrow \mathcal{C}(s,\mathbb{M})$ if $||\delta_r([x])||\in \Theta(n)$ and $\exists z\in \Theta(n)$ such that $\mathbb{L}_{[\delta_r([x])],[z]}$ is also a real axis of the CoP $\mathcal{C}(s,\mathbb{M})$. We say that the subsequence $\Theta(n)$ is \emph{stable} under the $r^{th}$ scale dilation $\delta_r$ if all points in $[x]\in \mathcal{C}(n,\mathbb{M})$ with $x\in \Theta(n)$ are stable.
\end{definition}

Here, we establish an important result in the special case where the base set is the set $\mathbb{N}$ of all positive integers.

\begin{proposition}
Let $\Theta(n)=\mathbb{N}_{n-1}$ and $\delta_r:\mathcal{C}(n)\longrightarrow \mathcal{C}(m)$ be a dilation. The subsequence $\Theta(n)$ is stable if and only if $n\geq m$.
\end{proposition}

\begin{proof}
In the case $m=n$, the dilation is trivial and the claim is trivially true. Suppose that the sequence $\Theta(n)$ is stable under the dilation 
\begin{align}
\delta_r:\mathcal{C}(n)\longrightarrow \mathcal{C}(m)\nonumber
\end{align}
and assume to the contrary that $n<m$. The dilation is an expansion. It follows that for all $[x]\in \mathcal{C}(n)$ with $x\in \Theta(n)$ there exists $z\in \Theta(n)$ such that $z+||\delta_r([x])|| =m$. Under the assumption $n<m$ and by Theorem \ref{properties}, we have the embedding $\mathcal{C}(n)\subset \mathcal{C}(m)$, and for all $x\in \Theta(n)$ we get $[x]\in \mathcal{C}(n)$ and $1+x\leq n<m$. There exists some $[y]\in \mathcal{C}(n)$ such that $\delta_r([y])=[1]$ but there exists no $z\in \Theta(n)$ such that $1+z=m$. It follows that the point $[y]$ is not a stable point under $\delta_r$. This contradicts the claim that $\Theta(n)$ is stable and so $n<m$ is impossible. Conversely, we suppose that $m<n$ and consider the dilation 
\begin{align}
\delta_r:\mathcal{C}(n)\longrightarrow \mathcal{C}(m).\nonumber
\end{align}
We note that for any point $[x]\in \mathcal{C}(n)$ there exist some $k<m<n$ such that $||\delta_r([x])||+k=m$. Because $k\in \mathbb{N}_{n-1}=\Theta(n)$, it follows that the subsequence $\Theta(n)$ is stable under any dilation $\delta_r$.
\end{proof}
\bigskip

Here, we show that any consecutive terms in a subsequence of $\mathbb{N}_n$ containing none of its degenerate terms must be stable under the simple dilation. We formalize this assertion in the following results.

\begin{proposition}
Let $\Theta(n):=\left \{x,x+1,\ldots,n-x,n-x+1\right \}$ be a subsequence of $\mathbb{N}_n$ for any $1<x<\frac{n}{2}$. Also let $\delta_r:\mathcal{C}(n)\longrightarrow \mathcal{C}(n+1)$ be an expansion. The sequence $\Theta(n)$ is stable under the expansion $\delta_r$.
\end{proposition}

\begin{proof}
For any point $[x]\in \mathcal{C}(n)$, we observe that $\mathbb{L}_{[x],[n-x]}$ is a real axis of the CoP. Enforcing $1<x<\frac{n}{2}$, we observe that the dilation $\delta_1:\mathcal{C}(n)\longrightarrow \mathcal{C}(n+1)$ with 
\begin{align}
\delta_1([x]):=
\begin{cases}
[x] \quad &\mbox{for }~1\leq x\leq n-1 \\
[n] \quad &\mbox{additional for}~x=1
\end{cases}
\end{align}
is achievable. It follows that for each $1<x<\frac{n}{2}$, the line $\mathbb{L}_{[x],[n-x+1]}$ is also a real axis of the CoP $\mathcal{C}(n+1)$. This proves that $\Theta(n)$ is stable under the dilation $\delta_r$.
\end{proof}

\section{The density of points on the circle of partition}

In this section, we introduce the notion of \emph{density} of points on CoP $\mathcal{C}(n,\mathbb{M})$ for $\mathbb{M}\subseteq \mathbb{N}$. We refrain from using the attribute \emph{real} in this section.

\begin{definition}
Let $\mathbb{H}\subset\mathbb{N}$. The limit
\[
\mathcal{D}\left(\mathbb{H}\right)=\lim_{n\rightarrow\infty}
\frac{\vert\mathbb{H}\cap \mathbb{N}_n\vert}{n}
\]
denotes the density of $\mathbb{H}$ if it exists and is finite.
\end{definition}

\begin{definition}\label{pointdensity}
Let $\mathcal{C}(n,\mathbb{M})$ be CoP with $\mathbb{M}\subset \mathbb{N}$ and $n\in \mathbb{N}$. Suppose that $\mathbb{H}\subset \mathbb{M}$. By the \emph{density of points} $[x]\in \mathcal{C}(n,\mathbb{M})$ such that $x\in \mathbb{H}$, denoted by $\mathcal{D} (\mathbb{H}_{\mathcal{C}(\infty,\mathbb{M})})$, we mean the limit
\begin{align}
\mathcal{D}\left(\mathbb{H}_{\mathcal{C}(\infty,\mathbb{M})}\right)=\lim \limits_{n\longrightarrow \infty}\frac{\# \lbrace\mathbb{L}_{[x],[y]} \inn \mathcal{C}(n,\mathbb{M})|~ \{x,y\} \cap \mathbb{H}\neq \emptyset \rbrace}{ \nu(n,\mathbb{M})}\nonumber
\end{align}
if it exists and is finite.
\end{definition}
\bigskip

The notion of the density of points introduced in the definition \ref{pointdensity} allows a passage between the density of the corresponding weight set of points. This possibility renders this type of density as black box in studying problems concerning the partition of numbers into subsets of the integers taking into consideration their density.

\begin{proposition}\label{inequality}
Let $\mathcal{C}(n)$ with $n\in \mathbb{N}$ be a CoP and $\mathbb{H}\subset \mathbb{N}$. The following inequality holds: 
\begin{align}
\mathcal{D}(\mathbb{H})=\lim \limits_{n\longrightarrow \infty}\frac{\left \lfloor \frac{|\mathbb{H}\cap \mathbb{N}_n|}{2}\right \rfloor}{\left \lfloor \frac{n-1}{2}\right \rfloor}\leq \mathcal{D}(\mathbb{H}_{\mathcal{C}(\infty)})\leq \lim \limits_{n\longrightarrow \infty}\frac{|\mathbb{H}\cap \mathbb{N}_n|}{\left \lfloor \frac{n-1}{2}\right \rfloor}=2\mathcal{D}(\mathbb{H}).\nonumber
\end{align}
\end{proposition}

\begin{proof}
The upper bound is obtained from a configuration in which there are no two points $[x],[y]\in \mathcal{C}(n)$ such that $x,y\in \mathbb{H}$ lie on the same axis of the CoP. That is, by the uniqueness of the axes of CoPs with $\nu(n,\mathbb{H})=0$, we can write
\begin{align}
   \# \left \{\mathbb{L}_{[x],[y]}\in \mathcal{C}(n)|~\{x,y\}\cap \mathbb{H}\neq \emptyset \right \}&=\nu(n,\mathbb{H})+\# \left \{\mathbb{L}_{[x],[y]}\in \mathcal{C}(n)|~x\in \mathbb{H},~y\in \mathbb{N}\setminus \mathbb{H}\right \} \nonumber \\&=\# \left \{\mathbb{L}_{[x],[y]}\in \mathcal{C}(n)|~x\in \mathbb{H},~y\in \mathbb{N}\setminus \mathbb{H}\right \} \nonumber \\&=|\mathbb{H}\cap \mathbb{N}_n|.\nonumber
\end{align}
However, the lower bound follows from a configuration in which any two points $[x],[y]\in \mathcal{C}(n)$ with $x,y\in \mathbb{H}$ are joined by an axis of the CoP. That is, by the uniqueness of the axis of CoPs with 
$$
\# \left \{\mathbb{L}_{[x],[y]}\in \mathcal{C}(n)|~x\in \mathbb{H},~y\in \mathbb{N}\setminus \mathbb{H}\right \}=0
$$
we can write 
\begin{align}
    \# \left \{\mathbb{L}_{[x],[y]}\in \mathcal{C}(n)|~\{x,y\}\cap \mathbb{H}\neq \emptyset \right \}&=\nu(n,\mathbb{H})\nonumber \\&=\left\lfloor\frac{|\mathbb{H}\cap\mathbb{N}_n|}{2}\right \rfloor.\nonumber
\end{align}
\end{proof}
\bigskip

\begin{proposition}\label{propertydensity}
Let $\mathbb{H}\subset \mathbb{N}$ and suppose that $\mathcal{D}(\mathbb{H}_{\mathcal{C}(\infty)})$ exists. Then the following properties hold:
\begin{enumerate}
    \item [(i)] $\mathcal{D}(\mathbb{N}_{\mathcal{C}(\infty)})=1$ and $\mathcal{D}(\mathbb{H}_{\mathcal{C}(\infty)})\leq 1$. Furthermore, $\mathcal{D}(\mathbb{H}_{\mathcal{C}(\infty)})<1$ provided $\mathcal{D}(\mathbb{N}\setminus \mathbb{H})>0$.
    \bigskip
    
    \item [(ii)] $1-\lim \limits_{n\longrightarrow \infty}\dfrac{\nu(n,\mathbb{N}\setminus\mathbb{H})}{\nu(n,\mathbb{N})}=\mathcal{D}(\mathbb{H}_{\mathcal{C}(\infty)})$.
    \bigskip
    
    \item [(iii)] If $|\mathbb{H}|<\infty$, then $\mathcal{D}(\mathbb{H}_{\mathcal{C}(\infty)})=0$.
\end{enumerate}
\end{proposition}

\begin{proof}
We observe that the first part of \textbf{Property} $(i)$ and $(iii)$ are easy consequences of the definition of density of points on the CoP $\mathcal{C}(n)$ and Proposition \ref{inequality}. We establish the second part of the \textbf{Property} $(i)$ and \textbf{Property} $(ii)$, which is the less obvious case. We observe that, by the uniqueness of the axes of CoPs, we can write 
\begin{align*}
    1&=\lim \limits_{n\longrightarrow \infty}\frac{\nu(n,\mathbb{N})}{\nu(n,\mathbb{N})}\\
    &=\lim \limits_{n\longrightarrow \infty}\frac{\# \lbrace\mathbb{L}_{[x],[y]} \inn \mathcal{C}(n)|~ x\in \mathbb{H}~,y\in \mathbb{N}\setminus \mathbb{H}\rbrace}{\nu(n,\mathbb{N})}\\
    &+\lim \limits_{n\longrightarrow \infty}\frac{\nu(n,\mathbb{H})}{\nu(n,\mathbb{N})}
    +\lim \limits_{n\longrightarrow \infty}\frac{\nu(n,\mathbb{N}\setminus\mathbb{H})}{\nu(n,\mathbb{N})}\\
    &=\mathcal{D}(\mathbb{H}_{\mathcal{C}(\infty)})
    +\lim \limits_{n\longrightarrow \infty}\frac{\nu(n,\mathbb{N}\setminus\mathbb{H})}{\nu(n,\mathbb{N})}
\end{align*}
and $(ii)$ follows immediately. The second part of $(i)$ follows from the above expression using the inequality
\begin{align}
   \lim \limits_{n\longrightarrow \infty}\frac{\nu(n,\mathbb{N}\setminus\mathbb{H})}{\nu(n,\mathbb{N})}&\leq \lim \limits_{n\longrightarrow \infty}\frac{\left \lfloor \frac{|\mathbb{N}\setminus \mathbb{H}\cap \mathbb{N}_n|}{2}\right \rfloor}{\left \lfloor \frac{n-1}{2}\right \rfloor}=\mathcal{D}(\mathbb{N}\setminus \mathbb{H})\nonumber
\end{align}
\end{proof}
\bigskip

Here, we transfer the notion of the density of a sequence to the density of corresponding points on the CoP $\mathcal{C}(n)$. This notion will play a crucial role in our latter developments.

\begin{proposition}
Let $\epsilon\in(0,1]$ with $\mathbb{H}\subset\mathbb{N}$ and $\mathcal{C}(n)$ be a CoP. If $\mathcal{D}\left(\mathbb{H}\right)\geq \epsilon$, then $\mathcal{D}\left(\mathbb{H}_{\mathcal{C}(\infty)}\right)\geq \epsilon$.
\end{proposition}

\begin{proof}
The result follows by exploiting the inequality in Proposition \ref{inequality}
\end{proof}

\begin{proposition}
Let $\mathbb{H}\subset \mathbb{N}$. For $\epsilon\in(0,1]$ and any $k\in \mathbb{N}$ if 
\begin{align}
|\mathbb{H}\cap\mathbb{N}_n|\geq n\epsilon \nonumber
\end{align}
and the common difference of arithmetic progressions in $\bigg(\mathbb{N}\setminus \mathbb{H}\bigg)\cap \mathbb{N}_n$ are different from those in $\mathbb{H}\cap \mathbb{N}_n$, then there exists some rotation $\varpi_r$ such that CoP $\mathcal{C}(n)$ contains at least $(k-1)$ stable points $[x]$ for $x\in \mathbb{H}\cap\mathbb{N}_n$.
\end{proposition}

\begin{proof}
Suppose that $\mathbb{H}\subset \mathbb{N}$ with the underlying conditions. By Theorem \ref{szemeredi} the set $\mathbb{H}$ contains fairly long arithmetic progressions of length $k$. We enumerate them as follows 
\begin{align}
x,x+s,x+2s,\ldots,x+(k-1)s\nonumber
\end{align}
for $s\in \mathbb{N}$. It follows that the corresponding points on the CoP $\mathcal{C}(n)$, namely 
\begin{align}
[x],[x+s],[x+2s],\ldots,[x+(k-1)s]\in\mathcal{C}(n)\nonumber
\end{align}
are equally spaced and the chord that joins two of these adjacent points are of equal distance. Similarly, the points on the other end of the axis are equally spaced and the chords that join any of these two adjacent points are of equal distance $s$. We enumerate them as follows 
\begin{align}
[n-x],[n-x-s],[n-x-2s],\ldots, [n-x-(k-1)s]\in \mathcal{C}(n).\nonumber
\end{align}
We apply the rotation $\varpi_r$ by choosing $r=s$ and get
\begin{align}
\varpi_s([x]),\varpi_s([x+s]),\ldots,\varpi_s([x+(k-1)s]).\nonumber
\end{align}
The image of these points under the rotation is given by 
\begin{align}
[x+s],[x+2s],\ldots,[x+(k-1)s], [x+ks].\nonumber
\end{align}
Since the point $[x+ks]$ was not on any of the axes considered, at least $(k-1)$ points on these axes will be transferred to their immediate next point on an axis containing all points $[x]$ with $x\in \mathbb{H}\cap \mathbb{N}_n$. Similarly, under the rotation, the corresponding images of the points on the other half of the CoP lying on the same axis with these points have the images 
\begin{align}
\varpi_s([n-x]),\varpi_s([n-x-s]),\ldots,\varpi_s([n-x-(k-1)s])\nonumber
\end{align}
which we can recast as 
\begin{align}
[n-x-s],[n-x-2s],\ldots,[n-x-(k-1)s],[n-x-ks].\nonumber
\end{align}
At least $(k-1)$ of these points are points on the previous axis and they reside on the same axis with the points on the other half of the CoP. Since the sequence 
\begin{align}
n-x-s,n-x-2s,\ldots,n-x-ks\nonumber
\end{align}
are in arithmetic progression, it follows by the assumption
\begin{align}
n-x-s,n-x-2s,\ldots,n-x-ks\in \mathbb{H}\cap \mathbb{N}_n.\nonumber
\end{align}
This completes the proof.
\end{proof}

\begin{proposition}
Let $\mathbb{H}\subset\mathbb{N}$ be such that $\mathbb{H}=\mathbb{J}\cup\mathbb{T}$ with $\mathbb{J}\cap \mathbb{T}=\emptyset$ and $\mathcal{D}(\mathbb{T})=0$. The following inequalities hold for the density of points in CoPs:
\begin{align}
   \mathcal{D}(\mathbb{H}_{\mathcal{C}(\infty)})=\mathcal{D}(\mathbb{J}_{\mathcal{C}(\infty)})\nonumber
\end{align}
and 
\begin{align}
    \mathcal{D}((\mathbb{N}\setminus \mathbb{H})_{\mathcal{C}(\infty)})\leq \mathcal{D}((\mathbb{N}\setminus \mathbb{J})_{\mathcal{C}(\infty)}).\nonumber
\end{align}
\end{proposition}

\begin{proof}
Applying the proposition \ref{propertydensity}, we obtain by the uniqueness of the axes of CoPs the following decomposition:
\begin{align}
   \mathcal{D}(\mathbb{H}_{\mathcal{C}(\infty)})&=\lim \limits_{n\longrightarrow \infty}\frac{\# \left \{\mathbb{L}_{[x],[y]}\in \mathcal{C}(n)|~\{x,y\}\cap (\mathbb{J}\cup \mathbb{T})\neq \emptyset \right\}}{\left \lfloor \frac{n-1}{2}\right \rfloor}\nonumber \\&= \lim \limits_{n\longrightarrow \infty}\frac{\# \left \{\mathbb{L}_{[x],[y]}\in \mathcal{C}(n,\mathbb{J})\right\}}{\left \lfloor \frac{n-1}{2}\right \rfloor}\nonumber \\&+\lim \limits_{n\longrightarrow \infty}\frac{\# \left \{\mathbb{L}_{[x],[y]}\in \mathcal{C}(n,\mathbb{T})\right\}}{\left \lfloor \frac{n-1}{2}\right \rfloor}\nonumber \\&+\lim \limits_{n\longrightarrow \infty}\frac{\# \left \{\mathbb{L}_{[x],[y]}\in \mathcal{C}(n)|~x\in \mathbb{J},~y\in \mathbb{T}\right\}}{\left \lfloor \frac{n-1}{2}\right \rfloor}\nonumber \\&+ \lim \limits_{n\longrightarrow \infty}\frac{\# \left \{\mathbb{L}_{[x],[y]}\in \mathcal{C}(n)|~x\in \mathbb{J},~y\in \mathbb{N}\setminus \mathbb{J}\right\}}{\left \lfloor \frac{n-1}{2}\right \rfloor}\nonumber \\&+\lim \limits_{n\longrightarrow \infty}\frac{\#\left\{\mathbb{L}_{[x],[y]}\in \mathcal{C}(n)|~x\in \mathbb{T},~y\in \mathbb{N}\setminus \mathbb{T}\right\}}{\left\lfloor\frac{n-1}{2}\right \rfloor}\nonumber
\end{align}
Under the inequalities 
\begin{align}
    \lim \limits_{n\longrightarrow \infty}\frac{\#\left\{\mathbb{L}_{[x],[y]}\in \mathcal{C}(n,\mathbb{T})\right\}}{\left \lfloor \frac{n-1}{2}\right \rfloor}\leq \lim \limits_{n\longrightarrow}\frac{\lfloor \frac{|\mathbb{T}\cap \mathbb{N}_n|}{2}\rfloor}{\left \lfloor \frac{n-1}{2}\right \rfloor}=0\nonumber
\end{align}
and 
\begin{align}
    \lim \limits_{n\longrightarrow \infty}\frac{\#\left \{\mathbb{L}_{[x],[y]}\in\mathcal{C}(n)|~x\in\mathbb{J},~y\in \mathbb{T}\right\}}{\left \lfloor \frac{n-1}{2}\right \rfloor}\leq \lim \limits_{n\longrightarrow \infty}\frac{|\mathbb{T}\cap\mathbb{N}_n|}{\left\lfloor \frac{n-1}{2}\right\rfloor}=0\nonumber
\end{align}
and 
\begin{align}
    \lim \limits_{n\longrightarrow \infty}\frac{\#\left\{\mathbb{L}_{[x],[y]}\in\mathcal{C}(n)|~x\in\mathbb{T},~y\in \mathbb{N}\setminus \mathbb{T}\right\}}{\left\lfloor\frac{n-1}{2}\right\rfloor}&\leq \lim\limits_{n\longrightarrow \infty}\frac{|\mathbb{T}\cap\mathbb{N}_n|}{\left\lfloor \frac{n-1}{2}\right\rfloor}=0\nonumber
\end{align}
we have 
\begin{align}
    \mathcal{D}(\mathbb{H}_{\mathcal{C}(\infty)})&=\lim \limits_{n\longrightarrow \infty}\frac{\#\left\{\mathbb{L}_{[x],[y]}\in \mathcal{C}(n)|~\{x,y\}\cap (\mathbb{J}\cup\mathbb{T})\neq \emptyset\right\}}{\left \lfloor\frac{n-1}{2}\right\rfloor}\nonumber \\&=\lim \limits_{n\longrightarrow \infty}\frac{\#\left\{\mathbb{L}_{[x],[y]}\in\mathcal{C}(n,\mathbb{J})\right\}}{\left\lfloor\frac{n-1}{2}\right\rfloor}\nonumber \\&+\lim \limits_{n\longrightarrow \infty}\frac{\#\left\{\mathbb{L}_{[x],[y]}\in\mathcal{C}(n)|~x\in\mathbb{J},~y\in \mathbb{N}\setminus\mathbb{J}\right\}}{\left\lfloor\frac{n-1}{2}\right\rfloor}\nonumber \\&=\mathcal{D}(\mathbb{J}_{\mathcal{C}(\infty)}).\nonumber
\end{align}
Again, we have the following decomposition by virtue of $\mathbb{J}\cap \mathbb{T}=\emptyset$
\begin{align}
    \mathcal{D}((\mathbb{N}\setminus \mathbb{H})_{\mathcal{C}(\infty)})&=\lim\limits_{n\longrightarrow \infty}\frac{\# \left\{\mathbb{L}_{[x],[y]}\in\mathcal{C}(n)|~\{x,y\}\cap (\mathbb{N}\setminus\mathbb{J}\cup\mathbb{T})\neq\emptyset \right\}}{\left\lfloor\frac{n-1}{2}\right\rfloor}\nonumber \\&=\lim\limits_{n\longrightarrow\infty}\frac{\#\left\{\mathbb{L}_{[x],[y]}\in\mathcal{C}(n,\mathbb{N}\setminus \mathbb{J}\cap\mathbb{N}\setminus\mathbb{T})\right\}}{\left \lfloor\frac{n-1}{2}\right\rfloor}\nonumber \\&+\lim \limits_{n\longrightarrow \infty}\frac{\#\left\{\mathbb{L}_{[x],[y]}\in\mathcal{C}(n)|~x\in \mathbb{N}\setminus\mathbb{J}\cap\mathbb{N}\setminus \mathbb{T},~y\in\mathbb{J}\cup\mathbb{T}\right\}}{\left \lfloor\frac{n-1}{2}\right\rfloor}\nonumber \\&\leq\lim \limits_{n\longrightarrow \infty}\frac{\#\left\{\mathbb{L}_{[x],[y]}\in\mathcal{C}(n,\mathbb{N}\setminus \mathbb{J})\right\}}{\left\lfloor\frac{n-1}{2}\right \rfloor}\nonumber \\&+\lim\limits_{n\longrightarrow \infty}\frac{\#\left\{\mathbb{L}_{[x],[y]}\in\mathcal{C}(n)|~x\in\mathbb{N}\setminus\mathbb{J},~y\in \mathbb{J}\right\}}{\left\lfloor\frac{n-1}{2}\right \rfloor}\nonumber \\&+\lim\limits_{n\longrightarrow \infty}\frac{\#\left\{\mathbb{L}_{[x],[y]}\in\mathcal{C}(n)|~x\in \mathbb{N}\setminus\mathbb{J},~y\in \mathbb{T}\right\}}{\left\lfloor\frac{n-1}{2}\right \rfloor}\nonumber
\end{align}
since 
$$
\mathbb{N}\setminus\mathbb{J}\cap\mathbb{N}\setminus \mathbb{T}\subset \mathbb{N}\setminus \mathbb{J}.
$$ 
Exploiting the inequality
\begin{align}
    \lim \limits_{n\longrightarrow \infty}\frac{\#\left\{\mathbb{L}_{[x],[y]}\in\mathcal{C}(n)|~x\in \mathbb{N}\setminus\mathbb{J},~y\in \mathbb{T}\right\}}{\left \lfloor\frac{n-1}{2}\right\rfloor}&\leq \lim \limits_{n\longrightarrow \infty}\frac{|\mathbb{T}\cap \mathbb{N}_n|}{\left\lfloor\frac{n-1}{2}\right \rfloor}=0\nonumber
\end{align}
we obtain
\begin{align}
    \mathcal{D}((\mathbb{N}\setminus\mathbb{H})_{\mathcal{C}(\infty)})&\leq\lim\limits_{n\longrightarrow \infty}\frac{\# \left\{\mathbb{L}_{[x],[y]}\in\mathcal{C}(n,\mathbb{N}\setminus\mathbb{J})\right\}}{\left\lfloor \frac{n-1}{2}\right\rfloor}\nonumber \\&+\lim \limits_{n\longrightarrow \infty}\frac{\#\left\{\mathbb{L}_{[x],[y]}\in\mathcal{C}(n)|~x\in \mathbb{N}\setminus\mathbb{J},~y\in\mathbb{J}\right\}}{\left \lfloor\frac{n-1}{2}\right\rfloor}\nonumber \\&=\mathcal{D}((\mathbb{N}\setminus\mathbb{J})_{\mathcal{C}(\infty)}).\nonumber
\end{align}
From the above analysis, we get the inequalities
\begin{align}
\mathcal{D}\left(\mathbb{H}_{\mathcal{C}(\infty)}\right)=\mathcal{D}(\mathbb{J}_{\mathcal{C}(\infty)})\nonumber
\end{align}
and 
\begin{align}
\mathcal{D}\left((\mathbb{N}\setminus \mathbb{H})_{\mathcal{C}(\infty)}\right)\leq \mathcal{D}\left((\mathbb{N}\setminus \mathbb{J})_{\mathcal{C}(\infty)}\right).\nonumber
\end{align}
\end{proof}
\bigskip

In the accompanying demonstration, we will make use of degenerate and non-degenerate points of a given set of points on a CoP. However intricate the proof might seem to be, it can tracked down to just a simple principle. The highly dense nature of the sequence allows us to break their components into several ''boxes''. The closest components in each of these boxes are equidistant from each other. The ''residue'', which are not dense, will be thrown away into another box whose components are very sparse. We then translate a component by their gap if it ever happens to be in some dense box at the same time live on the same axis with other component. This forces the second component to also belong to some dense box. If the component on the same axis with another component does not belong to the dense box, then the components and the associated components must live in the sparse box. We can then move them into the dense box and repeat the arguments.  We make the following ideas and framework more precise.

\begin{definition}
Let $\mathcal{P}\subseteq \mathcal{C}(n,\mathbb{M})$ for $\mathbb{M}\subseteq \mathbb{N}$. The point $[x]\in \mathcal{P}$ is a degenerate point if the line that joins the point $[x]$ to the center (resp. deleted center) of the CoP $\mathcal{C}(n,\mathbb{M})$ is a boundary of the largest sector induced by the points in $\mathcal{P}$. Otherwise, we say that it is a non-degenerate point in $\mathcal{P}$.
\end{definition}
\bigskip

\begin{theorem}\label{master}
Let $\mathbb{H}\subset\mathbb{N}$ and suppose that $\mathcal{C}(n,\mathbb{H})\neq \emptyset$. If for any $\epsilon\in (0,1]$ holds 
\begin{align}
|\mathbb{H}\cap \mathbb{N}_n|\geq n\epsilon\nonumber
\end{align}
with 
\begin{align}
\mathcal{D}(\mathbb{N}\setminus \mathbb{H})=\lim \limits_{n\longrightarrow \infty}\frac{|(\mathbb{N}\setminus \mathbb{H})\cap \mathbb{N}_n|}{n}<\mathcal{D}(\mathbb{H}) \nonumber
\end{align}
then there exists a dilation $\delta_r:\mathcal{C}(n,\mathbb{H})\longrightarrow \mathcal{C}(n+r,\mathbb{H})$ such that 
\begin{align} 
\mathcal{C}(n+r,\mathbb{H})\neq\emptyset. \nonumber
\end{align}
\end{theorem}

\begin{proof}
Under the assumption $|\mathbb{H}\cap \mathbb{N}_n|\geq n\epsilon$ for any $\epsilon\in(0,1]$, the set $\mathbb{H}$ contains fairly long arithmetic progressions. We enumerate them as follows 
\begin{align}
\mathbb{G}_1=\left \{x_1+kd_1\in \mathbb{H}\right \}_{k=0}^{s_1;s_1\geq 1}.\nonumber
\end{align}
We consider the residual set 
\begin{align}
\mathbb{G}_2=\mathbb{H}\setminus \left \{x_1+kd_1\in \mathbb{H}\right \}_{k=0}^{s_1;s_1\geq 1}.\nonumber
\end{align}
We can partition the sequence $\mathbb{H}$ in the following way
\begin{align}
\mathbb{H}=\mathbb{G}_1\cup\mathbb{G}_2.\nonumber
\end{align}
If $\mathbb{G}_2$ is still dense, then we can repeat this process and further obtain a partition of $\mathbb{H}$ into three subsequences 
\begin{align}
\mathbb{H}=\mathbb{G}_1\cup\mathbb{G}_2\cup\mathbb{G}_3.\nonumber
\end{align}
By induction, we can partition the sequence $\mathbb{H}$ in the following way 
\begin{align}
\mathbb{H}=\bigcup_{i=1}^{m}\mathbb{G}_i\cup \mathbb{T}\label{decomposition}
\end{align}
where 
\begin{align}
\lim\limits_{n\longrightarrow \infty}\frac{|\mathbb{T}\cap \mathbb{N}_n|}{n}=0\nonumber
\end{align}
and $\mathbb{J}=\bigcup_{i=1}^{m}\mathbb{G}_i $ with $\mathcal{D}(\mathbb{J})=\mathcal{D}(\mathbb{H})$, since $\mathbb{J}\cap \mathbb{T}=\emptyset$ and $\mathbb{G}_i=\left \{x_i+kd_i\in \mathbb{H}\right \}_{k=0}^{s_i;s_i\geq 1}$. It suffices to work with the corresponding points on CoP $\mathcal{C}(n,\mathbb{H})$. By assumption $\mathcal{C}(n,\mathbb{H})\neq \emptyset$, we deduce that there exist some axes $\mathbb{L}_{[a],[b]}\inn \mathcal{C}(n,\mathbb{H})$. We suppose that 
\begin{align}
[b]\notin \bigcup_{i=1}^{m}\left \{[x_i+kd_i]\in\mathcal{C}(n,\mathbb{H})\right \}_{k=0}^{s_i;s_i\geq 1}\nonumber
\end{align}
for $b\in \mathbb{H}$. It follows that no two adjacent chords of equal length joining points in 
\begin{align}
\bigcup_{i=1}^{m}\left \{[x_i+kd_i]\in \mathcal{C}(n,\mathbb{H})\right \}_{k=0}^{s_i;s_i\geq 1}\nonumber
\end{align}
contains the point $[b]$. We suppose on the contrary that
\begin{align}
[a]\in \bigcup_{i=1}^{m}\left \{[x_i+kd_i]\in \mathcal{C}(n,\mathbb{H})\right \}_{k=0}^{s_i;s_i\geq 1}.\nonumber
\end{align}
We deduce $[a]\in\left\{[x_i+kd_i]\in \mathcal{C}(n,\mathbb{H})\right\}_{k=0}^{s_i;s_i\geq 1}$ for some $1\leq i\leq m$. We consider two cases. The case $[a]$ is a degenerate point in the set and the case is a non-degenerate point in the set. If $[a]$ is a degenerate point in the set $\left\{[x_i+kd_i]\in\mathcal{C}(n,\mathbb{H})\right \}_{k=0}^{s_i;s_i\geq 1}$, in particular, $[a]$ is the first point in the set, then the following points 
\begin{align}
[a], [x_i+d_i], [x_i+2d_i],\ldots [x_i+sd_i]\nonumber
\end{align}
are equally spaced with $b=n-x_i$. It implies that $b$ is contained in the arithmetic progression 
\begin{align}
n-x_i,n-x_i-d_i,\ldots,n-x_i-sd_i\nonumber
\end{align}
which contradicts the assumption that $[b]$ cannot lie on at least one of two adjacent chords of equal length. Otherwise 
\begin{align}
n-x_i,n-x_i-d_i,\ldots,n-x_i-sd_i\in \bigg(\mathbb{N}\setminus \mathbb{J}\bigg)\cap \mathbb{N}_n.\nonumber
\end{align}
This implies that each element in the weight set $\mathbb{J}\cap \mathbb{N}_n=\mathbb{K}$ of the corresponding point set $\mathbb{K}^{*}=\bigcup_{i=1}^{m}\left \{[x_i+kd_i]\in\mathcal{C}(n,\mathbb{H})\right\}_{k=0}^{s_i;s_i\geq 1}$ uniquely generates an element in the set $\bigg(\mathbb{N}\setminus \mathbb{J}\bigg)\cap\mathbb{N}_n$, so that 
\begin{align}
    |\mathbb{J}\cap\mathbb{N}_n|\leq |\bigg(\mathbb{N}\setminus \mathbb{J}\bigg)\cap\mathbb{N}_n|.\nonumber
\end{align}
It follows that 
\begin{align}
  \mathcal{D}(\mathbb{N}\setminus\mathbb{H})&<\mathcal{D}(\mathbb{H})\nonumber \\&=\mathcal{D}(\mathbb{J})\nonumber \\&\leq \mathcal{D}(\mathbb{N}\setminus\mathbb{J})\nonumber
\end{align}
so we have the inequality 
$$
\mathcal{D}(\mathbb{N}\setminus\mathbb{H})<\mathcal{D}(\mathbb{N}\setminus \mathbb{J}).
$$ 
This contradicts equality under equality $\mathcal{D}(\mathbb{H})=\mathcal{D}(\mathbb{J})$
\begin{align}
\mathcal{D}(\mathbb{N}\setminus\mathbb{H})&=1-\mathcal{D}(\mathbb{H})\nonumber \\&=1-\mathcal{D}(\mathbb{J})\nonumber \\&=\mathcal{D}(\mathbb{N}\setminus\mathbb{J}).\nonumber
\end{align}
For the case $[a]=[x_i+sd_i]$, we obtain the arithmetic progression with $b=n-x_i-sd_i$. The corresponding point $[b]$ also violates the required specification. If the point $[a]\in \left\{[x_i+kd_i]\in \mathcal{C}(n,\mathbb{H})\right \}_{k=0}^{s_i;s_i\geq 1}$ is a non-degenerate point, then $a=x_i+jd_i$ for some $0<j<s$. The same analysis can be performed to produce a contradiction. Now for the case 
\begin{align}
[a]\in \bigcup_{i=1}^{m}\left \{[x_i+kd_i]\in \mathcal{C}(n,\mathbb{H})\right \}_{k=0}^{s_i;s_i\geq 1}\nonumber
\end{align}
we choose the dilation $\delta_r$ with $r=d_j$ such that $[b]\in \left\{[x_j+kd_j]\in\mathcal{C}(n,\mathbb{H})\right \}_{k=0}^{s_i;s_i\geq 1}$ for $r<0$ if $[b]$ is the last degenerate point in the set and $r>0$ if $[b]$ is the first degenerate point or a non-degenerate point in the set so that 
\begin{align}
\mathbb{L}_{[a],[b+d_j]}\inn \mathcal{C}(n+d_j,\mathbb{H}).\nonumber
\end{align}
This completes the first part of the proof. For the second part, we assume that for the axis $\mathbb{L}_{[a],[b]}$ of $\mathcal{C}(n,\mathbb{H})$
\begin{align}
[a]\notin \bigcup_{i=1}^{m}\left \{[x_i+kd_i]\in \mathcal{C}(n,\mathbb{H})\right \}_{k=0}^{s_i;s_i\geq 1}.\nonumber
\end{align}
This implies
\begin{align}
[a]\in \mathbb{T}_n^{*}\nonumber
\end{align}
where $\mathbb{T}_n^{*}$ is the set of corresponding points of elements in $\mathbb{T}\cap\mathbb{N}_n$. Since 
\begin{align}
|\mathbb{T}_n^{*}|\leq \bigg|\bigcup_{i=1}^{m}\left\{[x_i+kd_i]\in \mathcal{C}(n,\mathbb{H})\right\}_{k=0}^{s_i;s_i\geq 1}\bigg|\nonumber
\end{align}
there exists some rotation $\varpi_t$ such that the point $\varpi_t([a])\in \bigcup_{i=1}^{m}\left\{[x_i+kd_i]\in \mathcal{C}(n,\mathbb{H})\right\}_{k=0}^{s_i;s_i\geq 1}$. In particular 
\begin{align}
\varpi_t([a])\in \left\{[x_j+kd_j]\in \mathcal{C}(n,\mathbb{H})\right\}_{k=0}^{s_i;s_i\geq 1}\nonumber
\end{align}
for some $1\leq j\leq m$. It follows that there must exist a point
\begin{align}
[v]\in\bigcup_{i=1}^{m}\left \{[x_i+kd_i]\in\mathcal{C}(n,\mathbb{H})\right\}_{k=0}^{s_i;s_i\geq 1}\nonumber
\end{align}
such that $\mathbb{L}_{[v],[\varpi_t([a])]}$ is an axis of CoP $\mathcal{C}(n,\mathbb{H})$, by virtue of the previous arguments. Otherwise, we discard this choice of point and search for a point with such property by varying the scale of rotation $\varpi_t$. The proof is completed by choosing the dilation $\delta_r$ such that $r=d_j$ for $r<0$ if $\varpi_t([a])]$ is the last degenerate point in the set and $r>0$ if $\varpi_t([a])]$ is the first degenerate point or a non-degenerate point in the set, so that $\mathbb{L}_{[v],[||\varpi_t([a])||+d_j]}$ is an axis of the CoP
\begin{align}
\mathcal{C}(n+d_j,\mathbb{H}).\nonumber
\end{align}
\end{proof}

\begin{theorem}\label{squarefree}
There are infinitely many $n\in \mathbb{M}_{a,d}$ with fixed $a,d\in\mathbb{N}$ such that the representation 
\begin{align}
n=z_1+z_2\nonumber
\end{align}
where $\mu(z_1)=\mu(z_2)\neq 0$, $z_1,z_2\in\mathbb{N}$ and $\mu$ is the M\"obius function defined as 
\begin{align}
\mu(m)=
\begin{cases}1\quad \mathrm{if}\quad m=1\\0 \quad \mathrm{if}\quad p^k|m,~k\in\mathbb{N}\setminus \{1\}\\ (-1)^r \quad \mathrm{if}\quad m=p_1p_2\cdots p_r
\end{cases}\nonumber
\end{align}
holds.
\end{theorem}

\begin{proof}
The set of square-free integers 
\begin{align}
\mathcal{Q}:=\left \{m\in \mathbb{N}:\mu(m)\neq 0\right \}\nonumber
\end{align}
has natural density $\frac{6}{\pi^2}$ (see, e.g, \cite{montgomery2007multiplicative, hildebrand2006introduction}). For $n$ large enough there exists some fixed $N_0>n$ such that the representation
\begin{align}
N_0=z_1+z_2\nonumber
\end{align}
holds with $\mu(z_1),\mu(z_2)\neq 0$.
By Theorem \ref{master} there exists some $t\in \mathbb{N}$ such that the representation
\begin{align}
N_t:=N_0+t=v_1+v_2\nonumber
\end{align}
holds with $\mu(v_1)=\mu(v_2)\neq 0$. The result follows by an upward induction in this manner.
\end{proof}
\bigskip

\begin{corollary}
There are infinitely many $n\in \mathbb{M}_{a,d}$ with fixed $a,d\in\mathbb{N}$ such that the representation  
\begin{align}
n=z_1+z_2\nonumber
\end{align}
holds with $\gcd(z_1,z_2)=1$ and $z_1,z_2\in\mathbb{N}$.
\end{corollary}

\begin{proof}
The set 
\begin{align}
\mathcal{R}:=\left \{(m,n):\gcd(m,k)=1,~1\leq m<k\right\} \nonumber
\end{align}
has natural density $\mathcal{D}(\mathcal{R})=\frac{6}{\pi^2}$ with relatively small density for the residual set \cite{hildebrand2006introduction}. The result follows by adapting a similar reasoning in Theorem \ref{squarefree}.
\end{proof}
\bigskip

\subsection{Application of density of points to partitions}

In this section, we explore the connection between the notion of density of points in a typical CoP to the possibility of partitioning numbers into certain sequences. This method tends to work very efficiently for sets of integers that have a positive density. The work in this section is similar to previous investigations in \cite{shnirel1939additive}.

\begin{theorem}\label{density to partition}
Let $\mathbb{H}\subset\mathbb{N}$ be such that $\mathcal{D}(\mathbb{H})>\frac{1}{2}$. Every sufficiently large $n\in\mathbb{N}$ has representation of the form 
\begin{align}
  n=z_1+z_2 \nonumber
\end{align}
where $z_1,z_2\in \mathbb{H}$.
\end{theorem}

\begin{proof}
By Proposition \ref{inequality}, we can write 
\begin{align}
    \lim \limits_{n\longrightarrow \infty}\frac{\left \lfloor \frac{|\mathbb{H}\cap \mathbb{N}_n|}{2}\right \rfloor}{\left \lfloor \frac{n-1}{2}\right \rfloor}\leq \mathcal{D}(\mathbb{H}_{\mathcal{C}(\infty)})\leq \lim \limits_{n\longrightarrow \infty}\frac{|\mathbb{H}\cap \mathbb{N}_n|}{\left \lfloor \frac{n-1}{2}\right \rfloor}.\nonumber
\end{align}
By the uniqueness of the axes of CoPs, we can write
\begin{align}
    \#\left\{\mathbb{L}_{[x],[y]}\inn \mathcal{C}(n)|~\{x,y\}\cap \mathbb{H}\neq \emptyset\right\}&=\nu(n,\mathbb{H})+\# \left \{\mathbb{L}_{[x],[y]}\inn \mathcal{C}(n)|~x\in \mathbb{H},~y\in \mathbb{N}\setminus \mathbb{H}\right \}.\nonumber
\end{align}
We assume that $\nu(n,\mathbb{H})=0$. By the definition \ref{pointdensity}, we get
\begin{align}
    \mathcal{D}(\mathbb{H}_{\mathcal{C}(\infty)})&=2\mathcal{D}(\mathbb{H})\nonumber \\&>2\times \frac{1}{2}=1.\nonumber
\end{align}
This contradicts the inequality $\mathcal{D}(\mathbb{H}_{\mathcal{C}(\infty)})\leq 1$ in Proposition \ref{propertydensity}. This proves that $\nu(n,\mathbb{H})>0$ for all sufficiently large values of $n\in\mathbb{N}$.
\end{proof}

\begin{corollary}
Let $\mathbb{R}:=\left\{m\in \mathbb{N}|~\mu(m)\neq 0\right\}$. Every sufficiently large $n\in \mathbb{N}$ can be written in the form 
\begin{align}
    n=z_1+z_2\nonumber
\end{align}
where $\mu(z_1)=\mu(z_2)\neq 0$.
\end{corollary}

\begin{proof}
By the uniqueness of the axes of CoPs, we can write
\begin{align}
    \#\left\{\mathbb{L}_{[x],[y]}\inn \mathcal{C}(n)|~\{x,y\}\cap \mathbb{R}\neq \emptyset\right\}&=\nu(n,\mathbb{R})+\#\left\{\mathbb{L}_{[x],[y]}\inn \mathcal{C}(n)|~x\in \mathbb{R},~y\in \mathbb{N}\setminus \mathbb{R}\right\}.\nonumber
\end{align}
We assume that $\nu(n,\mathbb{R})=0$. By the definition \ref{pointdensity} and Theorem \ref{density to partition}, we get
\begin{align}
    \mathcal{D}(\mathbb{R}_{\mathcal{C}(\infty)})&=2\mathcal{D}(\mathbb{R})\nonumber \\&=\frac{12}{\pi^2}>1\nonumber
\end{align}
since $\mathcal{D}(\mathbb{R})=\frac{6}{\pi^2}$. This contradicts the inequality $\mathcal{D}(\mathbb{R}_{\mathcal{C}(\infty)})\leq 1$ in Proposition \ref{propertydensity}. This proves that $\nu(n,\mathbb{R})>0$ for all sufficiently large values of $n\in \mathbb{N}$.
\end{proof}
\bigskip

One may hope that this strategy works when we replace the set $\mathbb{R}$ of square-free integers with the set of prime numbers. In this situation, we are bound to have some difficulty, since the primes have a zero natural density in accordance with the prime number theorem. Any progress in this regard is manifestly possible when we introduce some exotic forms of the notion of density of points and carefully choose a subset of the integers that is somewhat dense among the set of integers that covers the primes. Here, we propose a strategy somewhat similar to the above method for possibly getting a handle on the binary Goldbach conjecture and its variants. In the interim, we state and prove a conditional theorem concerning the binary Goldbach conjecture.

\begin{theorem}\label{conditional Goldbach}
Let $\mathbb{B}\subset\mathbb{N}$ such that $\mathbb{P}\subset \mathbb{B}$ with $|\mathcal{C}(n,\mathbb{B})|=|\mathbb{B}\cap \mathbb{N}_n|$ for all $n\in \mathbb{N}$ so that 
$$
\lim \limits_{n\longrightarrow \infty}\frac{|\mathbb{P}\cap \mathbb{N}_n|}{\eta(n)}>\frac{1}{2}
$$ 
where $\eta(n)$ is the counting function of all integers belonging to the set $\mathbb{B}\cap \mathbb{N}_n$. The counting function $\nu(n,\mathbb{P})>0$ for all sufficiently large values of $n\in 2\mathbb{N}$.
\end{theorem}

\begin{proof}
We first find an upper and lower bound for the density of points in the CoP $\mathcal{C}(n,\mathbb{B})$ belonging to the set of primes $\mathbb{P}$ so that under the condition $|\mathcal{C}(n,\mathbb{B})|=|\mathbb{B}\cap \mathbb{N}_n|$ for all $n\in \mathbb{N}$, we obtain the inequality
\begin{align}
\lim \limits_{n\longrightarrow \infty}\frac{\left\lfloor \frac{|\mathbb{P}\cap\mathbb{N}_n|}{2}\right\rfloor}{\left \lfloor \frac{\eta(n)-1}{2}\right\rfloor}\leq\mathcal{D}(\mathbb{P}_{\mathcal{C}(\infty,\mathbb{B})})\leq \lim \limits_{n\longrightarrow \infty}\frac{|\mathbb{P}\cap \mathbb{N}_n|}{\left\lfloor\frac{\eta(n)-1}{2}\right \rfloor}=2\lim \limits_{n\longrightarrow  \infty}\frac{|\mathbb{P}\cap \mathbb{N}_n|}{\eta(n)}.\nonumber
\end{align}
By the uniqueness of the axes of CoPs, we can write
\begin{align}
\#\left\{\mathbb{L}_{[x],[y]}\inn \mathcal{C}(n,\mathbb{B})|~\{x,y\}\cap\mathbb{P}\neq \emptyset\right\}&=\nu(n,\mathbb{P})+\# \left\{\mathbb{L}_{[x],[y]}\inn \mathcal{C}(n,\mathbb{B})|~x\in \mathbb{P},~y\in \mathbb{B}\setminus \mathbb{P}\right \}.\nonumber
\end{align}
We assume to the contrary that $\nu(n,\mathbb{P})=0$. It follows that no two points in the CoP $\mathcal{C}(n,\mathbb{B})$ with weight in the set $\mathbb{P}$ are axes partners, so that under the requirement 
$$
\lim \limits_{n\longrightarrow \infty}\frac{|\mathbb{P}\cap \mathbb{N}_n|}{\eta(n)}>\frac{1}{2}
$$ 
where $\eta(n)$ is the counting function of all integers belonging to the set $\mathbb{B}\cap \mathbb{N}_n$, we obtain the inequality 
\begin{align}
    \mathcal{D}(\mathbb{P}_{\mathcal{C}(\infty,\mathbb{B})})&=2\lim \limits_{n\longrightarrow  \infty}\frac{|\mathbb{P}\cap\mathbb{N}_n|}{\eta(n)}\nonumber \\&>2\times\frac{1}{2}=1.\nonumber
\end{align}
This contradicts the inequality $\mathcal{D}(\mathbb{P}_{\mathcal{C}(\infty,\mathbb{B})})\leq 1$ in Proposition \ref{propertydensity}. This shows that $\nu(n,\mathbb{P})>0$ for all sufficiently large values of $n\in 2\mathbb{N}$.
\end{proof}

\subsection{Binary Goldbach conjecture proof technique via circles of partition}

In this section, we propose a series of steps that could be taken to confirm the truth of the binary Goldbach conjecture. We enumerate the strategies chronologically as follows:

\begin{enumerate}
\item  We first construct a subset of the integers $\mathbb{B}$ that covers the primes with $|\mathcal{C}(n,\mathbb{B})|=|\mathbb{B}\cap \mathbb{N}_n|$ for all $n\in \mathbb{N}$ so that 
$$
\lim \limits_{n\longrightarrow \infty}\frac{|\mathbb{P}\cap \mathbb{N}_n|}{\eta(n)}>\frac{1}{2}
$$
where $\eta(n)$ is the counting function of all integers belonging to the set $\mathbb{B}\cap\mathbb{N}_n$.
\bigskip

\item  We remark that the following inequality also hold and this can be obtain by replacing the set $\mathbb{N}$  with the set $\mathbb{B}$
\begin{align}
\lim \limits_{n\longrightarrow \infty}\frac{\left\lfloor \frac{|\mathbb{P}\cap \mathbb{N}_n|}{2}\right\rfloor}{\left \lfloor\frac{\eta(n)-1}{2}\right\rfloor}\leq \mathcal{D}(\mathbb{P}_{\mathcal{C}(\infty,\mathbb{B})})\leq \lim \limits_{n\longrightarrow \infty}\frac{|\mathbb{P}\cap \mathbb{N}_n|}{\left\lfloor\frac{\eta(n)-1}{2}\right \rfloor}=2\lim \limits_{n\longrightarrow  \infty}\frac{|\mathbb{P}\cap\mathbb{N}_n|}{\eta(n)}.\nonumber
\end{align}
\bigskip

\item  By the uniqueness of the axes of CoPs, we can write  
\begin{align*}
    &\# \left \{\mathbb{L}_{[x],[y]}\inn \mathcal{C}(n,\mathbb{B})|~\{x,y\}\cap \mathbb{P}\neq \emptyset\right\}\\
    &=\nu(n,\mathbb{P})+\#\left\{\mathbb{L}_{[x],[y]}\inn \mathcal{C}(n,\mathbb{B})|~x\in\mathbb{P},~y\in \mathbb{B}\setminus\mathbb{P}\right\}.\nonumber
\end{align*}
\bigskip

\item We assume that $\nu(n,\mathbb{P})=0$ and get by the definition \ref{pointdensity}
\begin{align}
    \mathcal{D}(\mathbb{P}_{\mathcal{C}(\infty,\mathbb{B})})&=2\lim\limits_{n\longrightarrow  \infty}\frac{|\mathbb{P}\cap\mathbb{N}_n|}{\eta(n)}\nonumber \\&>2\times\frac{1}{2}=1.\nonumber
\end{align}
This contradicts the inequality $\mathcal{D}(\mathbb{P}_{\mathcal{C}(\infty,\mathbb{B})})\leq 1$ in Proposition \ref{propertydensity}. This proves that $\nu(n,\mathbb{P})>0$ for all sufficiently large values of $n\in 2\mathbb{N}$.
\end{enumerate}
\bigskip

\section{Open and Connected Circles of Partition}

In this section, we introduce the notion of \emph{open} CoP. We first introduce the notion of a \emph{path} connecting CoP. Also here and in the following sections only real axes are considered, the attribute \emph{real} is not used.

\begin{definition}
Let $\mathbb{M}\subseteq\mathbb{N}$ with $\mathcal{C}(n,\mathbb{M})\neq\emptyset$ and $\mathcal{C}(s,\mathbb{M})\neq \emptyset$ be two distinct CoPs. By the \emph{path} that joins CoP $\mathcal{C}(n,\mathbb{M})$ to CoP $\mathcal{C}(s,\mathbb{M})$, we mean the line that joins $[x]\in\mathcal{C}(n,\mathbb{M})$ to $[y]\in\mathcal{C}(s,\mathbb{M})$, denoted by $\mathcal{L}_{[x],[y]}$, such that $\mathcal{L}_{[x],[y]}$ is an axis of the CoP $\mathcal{C}(s,\mathbb{M})$
\[
\mathcal{L}_{[x],[y]}=\mathbb{L}_{[x],[y]}\inn\mathcal{C}(s,\mathbb{M}).
\]
We say that CoP $\mathcal{C}(n,\mathbb{M})$ is \emph{connected} to CoP $\mathcal{C}(s,\mathbb{M})$ if there exists such a path. We say that CoP $\mathcal{C}(n,\mathbb{M})$ is \emph{strongly connected} to some CoP $\mathcal{C}(m,\mathbb{M})$ if the connection exists for all possible dilations
\begin{align*}
\delta_r:\mathcal{C}(n,\mathbb{M})\longrightarrow \mathcal{C}(s,\mathbb{M})
\mbox{ by }s=n+r.
\end{align*}
with $\delta_r([x])=[y]$. We say that CoP $\mathcal{C}(n,\mathbb{\mathbb{M}})$ is \emph{fully connected} to CoP $\mathcal{C}(s,\mathbb{M})$ if there exists such a path for each $[x]\in\mathcal{C}(n,\mathbb{M})$.
\end{definition}

\begin{proposition}\label{P-common_points}
Let $\mathbb{M}\subseteq\mathbb{N}$ with $\mathcal{C}(n,\mathbb{M})\neq\emptyset$ and $\mathcal{C}(s,\mathbb{M})\neq \emptyset$ be two distinct CoPs with a common point $[x]$. The CoP $\mathcal{C}(n,\mathbb{M})$ is connected to CoP $\mathcal{C}(s,\mathbb{M})$.
\end{proposition}

\begin{proof}
Since $[x]\in\mathcal{C}(s,\mathbb{M})$ there must be an axis $\mathbb{L}_{[x],[s-x]}\inn\mathcal{C}(s,\mathbb{M})$. Since $[x]\in\mathcal{C}(n,\mathbb{M})$ there exists the path $\mathcal{L}_{[x],[s-x]}$. Hence, CoP $\mathcal{C}(n,\mathbb{M})$ is connected to CoP $\mathcal{C}(s,\mathbb{M})$. Otherwise, if there exists such a path $\mathcal{L}_{[x],[y]}$ with
a fixed $[x]\in\mathcal{C}(n,\mathbb{M})$ and any $[y]\in\mathcal{C}(s,\mathbb{M})$ such that $\mathbb{L}_{[x],[y]}\inn\mathcal{C}(s,\mathbb{M})$, then $[y]=[s-x]$ and $[x]$ is also a point of $\mathcal{C}(n,\mathbb{M})$. 
\end{proof}
\bigskip

\begin{proposition}
Let $\mathbb{M}\subseteq\mathbb{N}$ and $\mathcal{C}(n,\mathbb{M})$ be a CoP. If $\mathcal{C}(n,\mathbb{M})$ is fully connected to $\mathcal{C}(s,\mathbb{M})$, then 
\begin{align}
     \mathcal{C}(n,\mathbb{M})\subseteq\mathcal{C}(s,\mathbb{M}).\nonumber
\end{align}
\end{proposition}

\begin{proof}
We let $\mathbb{M}\subseteq\mathbb{N}$ and suppose that CoP $\mathcal{C}(n,\mathbb{M})$ is connected to CoP $\mathcal{C}(s,\mathbb{M})$. This implies that for each point $[x]\in \mathcal{C}(n,\mathbb{M})$ there exists an axis $\mathbb{L}_{[x],[y]}\inn \mathcal{C}(s,\mathbb{M})$ for some $[y]\in \mathcal{C}(s,\mathbb{M})$. It follows that $[x]\in\mathcal{C}(s,\mathbb{M})$. This completes the demonstration since the point $[x]$ is an arbitrary point in CoP $\mathcal{C}(n,\mathbb{M})$.
\end{proof}
\bigskip

\begin{theorem}\label{connectexistence}
Let $\mathbb{M}\subseteq\mathbb{N}$ and $\mathcal{C}(n,\mathbb{M})$ be any CoP that admits an aligned embedding. The CoP $\mathcal{C}(n,\mathbb{M})$ is strongly connected to some CoP $\mathcal{C}(m,\mathbb{M})$ that admits an aligned embedding.
\end{theorem}

\begin{proof}
We assume that $\mathcal{C}(n,\mathbb{M})$ is not strongly connected to any $\mathcal{C}(m,\mathbb{M})$. By the requirement that the CoPs admit aligned embedding, we may assume that $\mathcal{C}(n,\mathbb{M})\subset\mathcal{C}(s,\mathbb{M})$. The line $\mathbb{L}_{[x],[n-x]}$ is an axis of $\mathcal{C}(n,\mathbb{M})$ for any $[x]\in \mathcal{C}(n,\mathbb{M})$. It follows that $\mathbb{L}_{[x],[s-x]}$ is also an axis of the CoP $\mathcal{C}(s,\mathbb{M})$. Since no two CoPs are strongly connected there exists (by Theorem \ref{existence}) some dilation $\delta_{r_1}:\mathcal{C}(n,\mathbb{M})\longrightarrow \mathcal{C}(s,\mathbb{M})$ such that $[s-x]\neq \delta_{r_1}([x])$ for each $[x]\in \mathcal{C}(n,\mathbb{M})$. We produce a line $\mathcal{L}_{[x],[\delta_{r_1}([x])]}$ by joining $[x]$ to $\delta_{r_1}([x])$. Now, we can partition these lines as axes of large and small CoPs relative to the CoP $\mathcal{C}(s)$ as follows:
\begin{align}
\lbrace\mathbb{L}_{[x],\delta_{r_1}([x])}\inn \mathcal{C}(v,\mathbb{M})|~n<v\leq s-1\rbrace\bigcup  \lbrace\mathbb{L}_{[x],\delta_{r_1}([x])}\inn \mathcal{C}(k,\mathbb{M})|~k>s\rbrace.\nonumber
\end{align}
We now arbitrarily pick a small CoP relative to the CoP $\mathcal{C}(s,\mathbb{M})$ and large relative to the CoP $\mathcal{C}(n,\mathbb{M})$. In particular, we arbitrarily pick a CoP $\mathcal{C}(v,\mathbb{M})$ from the first set. We obtain the strict embedding
\begin{align}
    \mathcal{C}(n,\mathbb{M})\subset \mathcal{C}(v,\mathbb{M})\subset\mathcal{C}(s,\mathbb{M}).\nonumber
\end{align}
Otherwise, CoP $\mathcal{C}(n,\mathbb{M})$ will have the axis $\mathbb{L}_{[x],[\delta_{r_1}([x])]}$, which will violate the assumption. Under the assumption that no two CoPs are strongly connected, it follows that there exist some dilation
\begin{align}
    \delta_{r_2}:\mathcal{C}(n,\mathbb{M})\longrightarrow \mathcal{C}(v,\mathbb{M})\nonumber
\end{align}
such that for each $[x]\in \mathcal{C}(n,\mathbb{M})$ then $\delta_{r_2}([x])\neq [v-x]$. Repeating the argument in this manner under the assumption that no two CoPs are connected, we obtain the following infinite embedding (nesting) into the CoP $\mathcal{C}(s,\mathbb{M})$ as follows
\begin{align}
    \mathcal{C}(n,\mathbb{M})\subset \cdots \subset \mathcal{C}(t,\mathbb{M})\subset \mathcal{C}(v,\mathbb{M})\subset \mathcal{C}(s,\mathbb{M})\nonumber
\end{align}
and we have the following infinite descending sequence of generators toward the generator $n$
\begin{align}
    n<\cdots <t<v<s.\nonumber
\end{align}
This is impossible.
\end{proof}
\bigskip

\begin{corollary}
Let $\mathcal{C}(n,\mathbb{M})$ and $\mathcal{C}(n,\mathbb{M})$  be two CoPs that admit an aligned embedding. If $n<m$, then CoP $\mathcal{C}(n,\mathbb{M})$ is fully connected to CoP $\mathcal{C}(m,\mathbb{M})$.
\end{corollary}

\begin{proof}
By Theorem \ref{properties2}, we have
\[
\mathcal{C}(n,\mathbb{M})\subset\mathcal{C}(m,\mathbb{M}).
\]
Hence, each point of $\mathcal{C}(n,\mathbb{M})$ is also a point of $\mathcal{C}(m,\mathbb{M})$. By Proposition \ref{P-common_points}, CoP $\mathcal{C}(n,\mathbb{M})$ is connected to CoP $\mathcal{C}(m,\mathbb{M})$ for each point $[x]\in\mathcal{C}(n,\mathbb{M})$. Hence, CoP $\mathcal{C}(n,\mathbb{M})$ is fully connected to CoP $\mathcal{C}(m,\mathbb{M})$
\end{proof}

One could imagine an analog of the results of fully connected CoPs when we take the base set $\mathbb{M}$ to be the set $\mathbb{P}$ of prime numbers. Here, it becomes slightly complicated and will require a careful analysis of the situation. 

\begin{conjecture}\label{prime connectedness}
There are infinitely many pairs of fully connected CoPs of the form $\mathcal{C}(n,\mathbb{P})$.
\end{conjecture}

\begin{definition}
Let $\mathbb{M}\subseteq \mathbb{N}$ and $\mathcal{C}(n,\mathbb{M})$ be a CoP with $\mathbb{L}_{[x],[y]}\inn \mathcal{C}(n,\mathbb{M})$. By the \emph{open} CoP induced by the point $[x],[y]$, we mean the exclusion $\mathcal{C}(n,\mathbb{M})\setminus [x],[y]$. We call the points $[x],[y]$ the gates to the interior of the open CoP. We denote the induced open CoP by $\widehat{\mathcal{C}(n,\mathbb{M})_{[x],[y]}}\subset \mathcal{C}(n,\mathbb{M})$. We say that CoPs $\mathcal{C}(s,\mathbb{M})$ and $\mathcal{C}(n,\mathbb{M})$ form a two-member community if and only if there is a path that joins the gate $[x],[y]$ of $\widehat{\mathcal{C}(n,\mathbb{M})_{[x],[y]}}$ to CoP $\mathcal{C}(s,\mathbb{M})$.
\end{definition}

\section{Children, Offspring and Family Induced by Circles of Partition}

In this section, we introduce the notion of \emph{children}, the \emph{offspring} and \emph{family} induced by a typical CoP.

\begin{definition}
Let $\mathbb{M}\subseteq \mathbb{N}$ and $\mathcal{C}(n,\mathbb{M})\neq \emptyset$ and let $\left \{\mathbb{L}_{[u_i],[v_i]}\right \}_{i=1}^{N;N\geq 2}$ for some $N\geq 2$ be the set of all the axes. We say that CoP $\mathcal{C}(s,\mathbb{M})$ is a \emph{child} of CoP $\mathcal{C}(n,\mathbb{M})$ if there exist some axes $\mathbb{L}_{[u_k],[v_k]},\mathbb{L}_{[u_j],[v_j]}\in \left\{\mathbb{L}_{[u_i],[v_i]}\right\}_{i=1}^{N;N\geq 2}$ such that at least one of $\mathbb{L}_{[u_k],[u_j]},\mathbb{L}_{[u_k],[v_j]},\mathbb{L}_{[v_k],[u_j]},\mathbb{L}_{[v_k],[v_j]}$ is an axis of the child CoP $\mathcal{C}(s,\mathbb{M})$. This axis forms the \emph{principal axis} of the child CoP. We call the collection of all CoPs generated in this manner the \emph{offspring} of the \emph{parent} CoP $\mathcal{C}(n,\mathbb{M})$. The parent CoP $\mathcal{C}(n,\mathbb{M})$ together with its offspring forms a \emph{complete family} of CoPs. The \emph{size} of the family of CoPs is the number of CoPs in the family. A proper subset of a family is said to be an \emph{incomplete family} of CoPs.
\end{definition}

\begin{example}
We consider the CoP with $\Vert\mathcal{C}(20,\mathbb{P})\Vert=\left \{3,7,13,17\right \}$ with axes $\mathbb{L}_{[3],[17]}$ and $\mathbb{L}_{[7],[13]}$. We consider the following chords $\mathcal{L}_{[3],[7]}, \mathcal{L}_{[3],[13]},\mathcal{L}_{[7],[17]},\mathcal{L}_{[13],[17]}$. These chords correspond as principal axes to the following CoPs
\begin{align}
\mathcal{C}(10,\mathbb{P}),~\mathcal{C}(16,\mathbb{P}),~\mathcal{C}(24,\mathbb{P}),~\mathcal{C}(30,\mathbb{P}).\nonumber
\end{align}
Hence, we obtain a complete family of CoPs of size $5$.
\end{example}

\begin{proposition}\label{P-child}
Let $\mathcal{C}(n,\mathbb{M})$ be a non-empty CoP. Each axis point $[x]$ together with a point $[u]$ of another axis of $\mathcal{C}(n,\mathbb{M})$ generates a child $\mathcal{C}(s,\mathbb{M})$ of the parent $\mathcal{C}(n,\mathbb{M})$ with $s=\Vert[x]\Vert+\Vert[u]\Vert$.
\end{proposition}

\begin{proof}
We let $\mathbb{L}_{[x],[y]}$ and $\mathbb{L}_{[u],[v]}$ be two axes of $\mathcal{C}(n,\mathbb{M})$. By Proposition \ref{unique}, we have 
\[
\Vert[x]\Vert+\Vert[u]\Vert=s\neq n.
\]
Hence, the point $[x]$ and $[u]$ form the axis $\mathbb{L}_{[x],[v]}\inn\mathcal{C}(s,\mathbb{M})$ and $\mathcal{C}(s,\mathbb{M})$ is a child of $\mathcal{C}(n,\mathbb{M})$.
\end{proof}

\begin{proposition}\label{P-infinite_Mad}
Let $n\in \mathbb{N}$ with $\mathbb{M}\subseteq\mathbb{N}$ and $\mathcal{C}(n,\mathbb{M})$ be a CoP that admits an aligned embedding. If $\vert\mathcal{C}(n,\mathbb{M})\vert\geq 4$, then CoP $\mathcal{C}(n,\mathbb{M})$ admits an infinite chain of its descendants.
\end{proposition}

\begin{proof}
Due to $\vert\mathcal{C}(n,\mathbb{M})\vert\geq 4$, there is an axis point $[u]\in\mathcal{C}(n,\mathbb{M})$ with 
\[
u>\min\left(\Vert[w]\Vert\mid [w]\in\mathcal{C}(n,\mathbb{M})\right)
\] 
and a point $[v]\in\mathcal{C}(n,\mathbb{M})$ of another axis with $u+v=m>n$. It follows that there exists an axis $\mathbb{L}_{[u],[v]}\inn\mathcal{C}(m,\mathbb{M})$. Therefore $[u]\in\mathcal{C}(m,\mathbb{M})$. By Proposition \ref{P-child}, CoP $\mathcal{C}(m,\mathbb{M})$ is a child of CoP $\mathcal{C}(n,\mathbb{M})$. Since $m>n$ and $\mathcal{C}(n,\mathbb{M})$ admits an aligned embedding, we deduce 
\[
\mathcal{C}(n,\mathbb{M})\subset\mathcal{C}(m,\mathbb{M}).
\]
Now, we choose a point $[w]$ of $\mathcal{C}(m,\mathbb{M})$ and latter change its role to a parent. With the same procedure as above, we produce an axis $\mathbb{L}_{[u],[w]}\inn\mathcal{C}(r,\mathbb{M})$ with $[u],[w]\in\mathcal{C}(r,\mathbb{M})$ and
\[
\mathcal{C}(n,\mathbb{M})\subset\mathcal{C}(m,\mathbb{M})
\subset\mathcal{C}(r,\mathbb{M}).
\]
This procedure can be repeated infinitely often. We obtain an infinite chain of descendants of the CoP $\mathcal{C}(n,\mathbb{M})$ as its prime father.
\end{proof}

\begin{proposition}\label{breaker}
Let $\mathbb{M}\subseteq \mathbb{N}$ and $\mathcal{C}(n,\mathbb{M})$ be parent of a complete family. The CoP $\mathcal{C}(n,\mathbb{M})$ partitions the offspring into two incomplete families of equal sizes.
\end{proposition}

\begin{proof}
By Proposition \ref{P-child}, two points of distinct axes of CoP $\mathcal{C}(n,\mathbb{M})$ produces a child. Let 
\[
\mathbb{L}_{[u],[v]}, \mathbb{L}_{[x],[y]}\mid u<x
\]
be two arbitrary axes of $\mathcal{C}(n,\mathbb{M})$. Because $[u],[v]$ and $[x],[y]$ are axis points, we get
\begin{align*}
&n=u+v=x+y\mbox{ and therefore}\\
&v=x-u+y\mbox{ and because of }x>u\\
&v>y
\end{align*}
Hence, we get
\begin{align*}
&u<x<y<v\mbox{ and therefore}\\
&s_1:=u+x<s_2:=u+y<n=x+y\mbox{ and}\\
&t_1:=v+y>t_2:=v+x>n=v+u
\end{align*}
and a chain of children
\begin{align*}
\mathcal{C}(s_1,\mathbb{M}),\mathcal{C}(s_2,\mathbb{M})
&,\mathcal{C}(n,\mathbb{M}),
\mathcal{C}(t_2,\mathbb{M}),\mathcal{C}(t_1,\mathbb{M})
\mbox{ with}\\
s_1<s_2&<n<t_2<t_1.
\end{align*}
Therefore, we find that for all two axes $4$ children are produced, two on the left side of $\mathcal{C}(n,\mathbb{M})$ and two on the right side in a chain of children. Because CoP $\mathcal{C}(n,\mathbb{M})$ (for all two axes) is located in the middle of the chain, the parent CoP $\mathcal{C}(n,\mathbb{M})$ partitions its offspring in two halves, the incomplete families of equal sizes.   
\end{proof}

\begin{proposition}\label{P-child-embedding}
If the parent CoP admits an embedding, then their children admit an aligned embedding.
\end{proposition}

\begin{proof}
We look at the last proof and choose $[u]$ as the first point of the parent CoP $\mathcal{C}(n,\mathbb{M})$
\[
u:=\min \left(w\in\Vert\mathcal{C}(n,\mathbb{M})\Vert\right).
\]
We deduce
\begin{align*}
&[u]\in\mathcal{C}(s_1,\mathbb{M})\mbox{ and }
[u]\in\mathcal{C}(s_2,\mathbb{M})\mbox{ and}\\
&\max\left(w\in\Vert\mathcal{C}(s_1,\mathbb{M})\Vert\right)=x<y
=\max\left(w\in\Vert\mathcal{C}(s_2,\mathbb{M})\Vert\right)\\
&\mbox{and hence}\\
&\mathcal{C}(s_1,\mathbb{M})\subset\mathcal{C}(s_2,\mathbb{M})
\mbox{ under }s_1<s_2.
\end{align*}
Because $\mathcal{C}(n,\mathbb{M})$ admits an embedding, we deduce
\begin{align*}
\mathcal{C}(s_1,\mathbb{M})\subset\mathcal{C}(s_2,\mathbb{M})
&\subset \mathcal{C}(n,\mathbb{M})\subset\mathcal{C}(t_2,\mathbb{M})
\subset\mathcal{C}(t_1,\mathbb{M})\mbox{ under}\\
s_1<s_2&<n<t_2<t_1.
\end{align*}
\end{proof}

\begin{corollary}\label{C-children-chain}
If CoP $\mathcal{C}(n,\mathbb{M})$ admits an embedding and has $2k$ children, then we get (by Propositions \ref{breaker} and \ref{P-child-embedding}) for its complete family
\[
\mathcal{C}(s_1,\mathbb{M})\subset\ldots\subset\mathcal{C}(s_k,\mathbb{M})
\subset\mathcal{C}(n,\mathbb{M})\subset
\mathcal{C}(t_k,\mathbb{M})\subset\ldots\subset\mathcal{C}(t_1,\mathbb{M})
\]
and the following symmetry
\[
s_2-s_1=t_1-t_2,\ldots,s_k-s_{k-1}=t_{k-1}-t_k,n-s_k=t_k-n.
\]
\end{corollary}

\begin{proof}
The embedding chain is a direct consequence of the Propositions \ref{breaker} and \ref{P-child-embedding}. Now, we prove the symmetry of the differences of the children generators. We look at the proof of Proposition \ref{breaker} with
\begin{align*}
&u<x<y<v\mbox{ and therefore}\\
&s_1:=u+x<s_2:=u+y<n=x+y\mbox{ and}\\
&t_1:=v+y>t_2:=v+x>n=v+u
\end{align*}
for two arbitrary axes $\mathbb{L}_{[u],[v]},\mathbb{L}_{[x],[y]}$ of $\mathcal{C}(n,\mathbb{M})$. We deduce
\[
s_1<s_2<n<t_2<t_1
\]
and get
\begin{align*}
s_2-s_1=u+y-u-x&=\mathbf{y-x}\mbox{ and }n-s_2=x+y-u-y=\mathbf{x-u}\\
t_1-t_2=v+y-v-x&=\mathbf{y-x}\mbox{ and }t_2-n=v+x-v-u=\mathbf{x-u}.
\end{align*}
Because the two axes are arbitrary, this symmetry around the generator $n$ holds for all axes. This establishes the claim.
\end{proof}
\bigskip

\begin{theorem}\label{T-numberofchildren}
Let $\mathbb{M}\subseteq \mathbb{N}$ and $\mathcal{C}(n,\mathbb{M})$ be a CoP with $|\mathcal{C}(n,\mathbb{M})|=k$. The number of children in the family with parent $\mathcal{C}(n,\mathbb{M})$ satisfies the upper bound
\begin{align}
    \leq 2\left\lfloor\frac{k}{2}\right\rfloor\bigg(\left\lfloor \frac{k}{2}\right\rfloor-1\bigg)\nonumber
\end{align}
and the lower bound
\begin{align*}
\geq 2\left(n_a-2\right)=4\left(\left\lfloor\frac{k}{2}\right\rfloor-1\right)
\mbox{ with }n_a=2\left\lfloor\frac{k}{2}\right\rfloor.
\end{align*}
\end{theorem}

\begin{proof}
We first prove the upper bound. The CoP  $\mathcal{C}(n,\mathbb{M})$ with $|\mathcal{C}(n,\mathbb{M})|=k$ contains $\lfloor\frac{k}{2}\rfloor$ different axes. Each axis contains two points of the parent $\mathcal{C}(n,\mathbb{M})$ and determines children with at most $\lfloor\frac{k}{2}\rfloor-1$ number of axes. The upper bound follows from this counting argument. Now we prove the lower bound. By Corollary \ref{ordered_weights} the weights of the points of $\mathcal{C}(n,\mathbb{M})$ are strictly totally ordered. Now, we remove from this sequence the weight of the center if it exists. We are left with $n_a=2\left\lfloor\frac{k}{2}\right\rfloor$ weights. We enumerate them as
\[
x_1<x_2<\ldots<x_{n_a-1}<x_{n_a}
\]
and construct the following sequences
\begin{align*}
s_1:=x_1+x_2<s_2:=x_1+x_3<&\ldots<s_{n_a-2}:=x_1+x_{n_a-1}<x_1+x_{n_a}=n\\
&\mbox{and}\\
t_1:=x_{n_a}+x_{n_a-1}>t_2:=x_{n_a}+x_{n_a-2}>&\ldots
>t_{n_a-2}:=x_{na}+x_2>x_{n_a}+x_1=n.
\end{align*}
Hence, we obtain
\[
s_1<\ldots<s_{n_a-2}<n<t_{n_a-2}<\ldots<t_1
\]
and have at least $2(n_a-2)$ different generators for the children 
of $\mathcal{C}(n,\mathbb{M})$.
\end{proof}

\begin{remark}
If a CoP contains not more than $3$ points, then the CoP has no children. We call these CoPs \emph{childless}. If a CoP has two axes, then the CoP has $4$ children. Therefore, there are no CoPs with only one child or only two or three children.
\end{remark}

\begin{proposition}\label{P-notallMchildless}
Let $\mathbb{M}\subseteq\mathbb{N}$. There is no parent CoP $\mathcal{C}(n,\mathbb{M})$ that admits an embedding with $\vert\mathcal{C}(n,\mathbb{M})\vert\geq 4$ such that all the children are childless.
\end{proposition}

\begin{proof}
By Proposition \ref{P-child-embedding} the children admit an aligned embedding. Since the parent CoP has at least $4$ children, there are (by Proposition \ref{breaker}) at least $2$ children with generators $>n$. Because the children admit an aligned embedding, we get for a child $\mathcal{C}(s,\mathbb{M})$ with $s>n$
\[
\mathcal{C}(s,\mathbb{M})\supset\mathcal{C}(n,\mathbb{M})
\mbox{ and therefore }
\vert\mathcal{C}(s,\mathbb{M})\vert >\vert\mathcal{C}(n,\mathbb{M})\vert
\geq 4.
\]
Hence, there are at least two children with more than $4$ own children and hence not childless.
\end{proof}
\bigskip

Here, we show that one can always partition any complete family into incomplete families with equal dilation between the members.

\begin{lemma}[The regularity lemma]\label{regular}
The offspring of a CoP $\mathcal{C}(n,\mathbb{M})$ can be partitioned into incomplete families with equal scale dilation between their sequence of successive embeddings.
\end{lemma}

\begin{proof}
If there are no embedding among the children of the parent $\mathcal{C}(n,\mathbb{M})$, then we have a partition into one-member incomplete family and the dilation in each family is trivial. We now assume that $\mathcal{C}(s_1,\mathbb{M})\subset \mathcal{C}(s_2,\mathbb{M})\subset \cdots \subset \mathcal{C}(s_k,\mathbb{M})$ for $k\geq 2$ is a sequence of children of the  parent $\mathcal{C}(n,\mathbb{M})$ with equal scale dilation between successive embedding. If the sequence is all of the children of the parent $\mathcal{C}(n,\mathbb{M})$, then the parent must be inserted (by Corollary \ref{C-children-chain}) in the middle of the offset chain. Now, we remove from the chain the parent $\mathcal{C}(n,\mathbb{M})$ with the two closest children. Hence, we obtain a partition of collection of children in the embedding into two sub-chains of embedding with equal scale dilation between successive children, those to the left of the children closest to the parent $\mathcal{C}(n,\mathbb{M})$ and to the right of the children closest to the parent $\mathcal{C}(n,\mathbb{M})$. For the sequence removed from the sequence of children given below
\begin{align}
\mathcal{C}(s_i,\mathbb{M})\subset\mathcal{C}(n,\mathbb{M})\subset \mathcal{C}(s_{i+1},\mathbb{M})\nonumber
\end{align}
we remove the parent $\mathcal{C}(n,\mathbb{M})$ and obtain a third partition of offspring with equal scale dilation
\begin{align}
\mathcal{C}(s_i,\mathbb{M})\subset \mathcal{C}(s_{i+1},\mathbb{M}).\nonumber
\end{align}
For the case where not all children are contained in the a previous embedding, we have already obtained a partition of collection of children into an incomplete family with equal scale dilation between successive members. The remaining collection of children can also be partitioned into incomplete families by choosing an embedding with equal scale dilation between successive children.
\end{proof}

\begin{theorem}\label{connectoffspring}
The number of pairs of connected children in any complete family is satisfies
\[
\geq \frac{n_a(n_a-2)(n_a-3)}{2}=
2\left\lfloor\frac{k}{2}\right\rfloor\left(\left\lfloor\frac{k}{2}\right\rfloor-1\right)(n_a-3)
\geq 2\left\lfloor\frac{k}{2}\right\rfloor\left(\left\lfloor\frac{k}{2}\right\rfloor-1\right)
\]
if the parent CoP has $n_a$ axis points and $n_a=2\left\lfloor\frac{k}{2}\right\rfloor>3$.
\end{theorem}

\begin{proof}
By Proposition \ref{P-common_points}, two CoPs are connected if and only if they have a common point. The children are produced by pairs of points on different axes. Each such point $[x]$ of the parent CoP appears in at least $n_a-2$ children. Hence, there are $\frac{(n_a-2)(n_a-3)}{2}$ pairs of children containing the point $[x]$. Since there are $n_a$ axis points, we deduce the lower bound.
\end{proof}

In analogy with Theorem \ref{T-numberofchildren}, we observe that the number of pairs of connected children of a complete family is at least the number of its children. From the proof of Theorem \ref{connectoffspring}, we deduce that each child is connected to another child of the same family.

\begin{example}
We take as parent CoP
$$
\mathcal{C}(22,\mathbb{P})=\lbrace [3],[5],[11],[17],[19]\rbrace\rightarrow k=5,n_a=4.
$$
By Theorem \ref{T-numberofchildren}, it has maximal
\[
2\cdot 2\cdot 1=4
\]
children and by Theorem \ref{connectoffspring} at least
\[
\frac{4\cdot 2 \cdot 1}{2}=4
\]
pairs of connected children. For the children, we get
\begin{align*}
\mathcal{C}(8,\mathbb{P})&=\lbrace\textbf{[3]},\textbf{[5]}\rbrace\\
\mathcal{C}(20,\mathbb{P})&=\lbrace\textbf{[3]},[7],[13],\textbf{[17]}\rbrace\\
----&----------\\
\mathcal{C}(24,\mathbb{P})&=\lbrace\textbf{[5]},[7],[11],[13],[17],\textbf{[19]}\rbrace\\
\mathcal{C}(36,\mathbb{P})&=\lbrace[5],[7],[13],\textbf{[17]},\textbf{[19]},[23],[29],[31]\rbrace.
\end{align*}
We observe that the axis point $[3]$ appears twice among the children and there is a pair of children that contains the axis point $[3]$. The axis point $[5]$ appears thrice and there are at least three pairs of children that contains this axis point. The axis point $[17]$ appears thrice among children, while the axis point $[19]$ appears twice and is contained exactly one pair. Therefore, we have $8>6$ pairs of connected children with respect to the points of the parent CoP. However, we see that further more of these points are common points in the offset. Hence, there are six further pairs of connected children. In the preceding example, we note that each principal axes point is marked in boldface. Consequently, CoP $\mathcal{C}(24,\mathbb{P})$ contains six axis points and at most $12$ children with at least $36$ pairs of connected children.
\end{example}
\bigskip

\begin{theorem}\label{T-nonembedding}
Let $\mathbb{M}\subseteq\mathbb{N}$ and $\mathcal{C}(n,\mathbb{M})$ with $\mathcal{C}(m,\mathbb{M})$ be two CoPs. Also let $\hat{\mathcal{O}}_n$ and $\hat{\mathcal{O}}_m$ be their complete families. If $\vert\hat{\mathcal{O}}_n\vert<\vert\hat{\mathcal{O}}_m\vert$ and 
there exists a child $\mathcal{C}(s,\mathbb{M})\in\hat{\mathcal{O}}_n$ such that $\mathcal{C}(s,\mathbb{M})\not\in\hat{\mathcal{O}}_m$, then
\[
\mathcal{C}(n,\mathbb{M})\not\subset\mathcal{C}(m,\mathbb{M})
\mbox{ and }
\mathcal{C}(n,\mathbb{M})\not\supset\mathcal{C}(m,\mathbb{M}),
\]
which means that CoP $\mathcal{C}(n,\mathbb{M})$ and CoP $\mathcal{C}(m,\mathbb{M})$ do not admit an embedding.
\end{theorem}

\begin{proof}
Since CoP $\mathcal{C}(s,\mathbb{M})$ is a child of CoP $\mathcal{C}(n,\mathbb{M})$, there are two points $[x],[u]\in\mathcal{C}(n,\mathbb{M})$ with $x+u=s$. We get $\mathbb{L}_{[x],[u]}$ as the principal axis of $\mathcal{C}(s,\mathbb{M})$. Since $\mathcal{C}(s,\mathbb{M})$ is not a child of $\mathcal{C}(m,\mathbb{M})$, there are no points in $\mathcal{C}(m,\mathbb{M})$ with a weight sum equal to $s$. Therefore, the points $[x],[u]$ are not contained in $\mathcal{C}(m,\mathbb{M})$. We deduce $\mathcal{C}(n,\mathbb{M})\not\subset\mathcal{C}(m,\mathbb{M})$. Due to $\vert\hat{\mathcal{O}}_n\vert<\vert\hat{\mathcal{O}}_m\vert$, there is a child of $\mathcal{C}(m,\mathbb{M})$ which is not a child of $\mathcal{C}(n,\mathbb{M})$. Therefore $\mathcal{C}(n,\mathbb{M})\not\supset\mathcal{C}(m,\mathbb{M})$.  
\end{proof}
\bigskip

\begin{theorem}\label{T-embedding_child}
Let $\mathbb{M}\subseteq\mathbb{N}$ and let $\mathcal{C}(n,\mathbb{M})$ and $\mathcal{C}(m,\mathbb{M})$ be two CoPs that admit an aligned embedding. Assume that
\[
\mathcal{C}(n,\mathbb{M})\subset\mathcal{C}(m,\mathbb{M}).
\]
The CoP $\mathcal{C}(n,\mathbb{M})$ is a child of $\mathcal{C}(m,\mathbb{M})$. If there is a chord $\mathcal{L}_{[x],[y]}$ of $\mathcal{C}(n,\mathbb{M})$ with $x+y=m$, then $\mathcal{C}(m,\mathbb{M})$ is also a child of $\mathcal{C}(n,\mathbb{M})$. Furthermore, the complete family $\hat{\mathcal{O}}_n$ of $\mathcal{C}(n,\mathbb{M})$ is a subset of the complete family $\hat{\mathcal{O}}_m$ of $\mathcal{C}(m,\mathbb{M})$.   
\end{theorem}

\begin{proof}
Due to $\mathcal{C}(n,\mathbb{M})\subset\mathcal{C}(m,\mathbb{M})$, we deduce (by the definition \ref{aligned_embedding}) that $n<m$ and 
\[
\min\left(x\mid [x]\in\mathcal{C}(n,\mathbb{M})\right)=
\min\left(u\mid [u]\in\mathcal{C}(m,\mathbb{M})\right).
\]
All chords $\mathcal{L}_{[x],[y]}$ of $\mathcal{C}(n,\mathbb{M})$ are also chords of $\mathcal{C}(m,\mathbb{M})$ excluding the chords between points $[x],[y]\in\mathcal{C}(n,\mathbb{M})$ with $x+y=m$. Using the underlying embedding, we deduce that all chords of $\mathcal{C}(m,\mathbb{M})$ which are axes of $\mathcal{C}(n,\mathbb{M})$ produce the same child, the CoP $\mathcal{C}(n,\mathbb{M})$. Hence, CoP $\mathcal{C}(n,\mathbb{M})$ is a child of CoP $\mathcal{C}(m,\mathbb{M})$. If there is no chord $\mathcal{L}_{[x],[y]}$ of $\mathcal{C}(n,\mathbb{M})$ with $x+y=m$, then all children of $\mathcal{C}(n,\mathbb{M})$ are children of $\mathcal{C}(m,\mathbb{M})$. Hence, the complete family $\hat{\mathcal{O}}_n$ is a subset of the complete family $\hat{\mathcal{O}}_m$ in this case. \\

If such a chord of $\mathcal{C}(n,\mathbb{M})$ exists, then it is an axis of $\mathcal{C}(m,\mathbb{M})$, so that $\mathcal{C}(m,\mathbb{M})$ is a child of $\mathcal{C}(n,\mathbb{M})$. Because the parents belong to the complete family, it implies that the complete family $\hat{\mathcal{O}}_n$ is a subset of $\hat{\mathcal{O}}_m$ in this case.
\end{proof}

\subsection{Isomorphic Circles of Partition}

In this section, we introduce and study the notion of \emph{isomorphism} between CoPs. 

\begin{definition}\label{D-isomor}
Let $\mathbb{M}\subseteq \mathbb{N}$, and let $\mathcal{C}(n,\mathbb{M})$ and $\mathcal{C}(m,\mathbb{M})$ be parents with the complete families $\hat{\mathcal{O}}_n$ and $\hat{\mathcal{O}}_m$, respectively. We say that the parents $\mathcal{C}(n,\mathbb{M})$ and $\mathcal{C}(m,\mathbb{M})$ are isomorphic if 
\begin{align}
    \hat{\mathcal{O}}_m\cap\hat{\mathcal{O}}_n\neq\emptyset. \nonumber
\end{align}
We call the number $|\hat{\mathcal{O}}_m\cap \hat{\mathcal{O}}_n|$ the \emph{degree} of isomorphism. We denote this isomorphism by $\mathcal{C}(n,\mathbb{M})\cong \mathcal{C}(m,\mathbb{M})$. We say that the degree of isomorphism is \emph{high} if at least one of the following equality holds
\begin{align}
    \frac{|\hat{\mathcal{O}}_m\cap\hat{\mathcal{O}}_n|}{|\hat{\mathcal{O}}_n|}=1\nonumber
\end{align}
or 
\begin{align}
     \frac{|\hat{\mathcal{O}}_m\cap\hat{\mathcal{O}}_n|}{|\hat{\mathcal{O}}_m|}=1.\nonumber
\end{align}
Otherwise, we say that the degree of isomorphism is \emph{low}.
\end{definition}
\bigskip

\begin{proposition}\label{connect isomorphism}
Let $\mathcal{C}(n,\mathbb{M})$ and $\mathcal{C}(m,\mathbb{M})$ be two CoPs. If $\mathcal{C}(n,\mathbb{M})$ is connected to $\mathcal{C}(m,\mathbb{M})$ by at least three distinct paths, then $\mathcal{C}(n,\mathbb{M})\cong\mathcal{C}(m,\mathbb{M})$.
\end{proposition}

\begin{proof}
We assume that CoP $\mathcal{C}(n,\mathbb{M})$ is connected to CoP $\mathcal{C}(m,\mathbb{M})$ by at least three distinct paths. We deduce that there exist some distinct points $[x],[y],[z]\in \mathcal{C}(m,\mathbb{M})$ such that $\mathbb{L}_{[x],[u]},\mathbb{L}_{[y],[v]},\mathbb{L}_{[z],[w]}\inn \mathcal{C}(m,\mathbb{M})$. It follows that there exist at least the following chords $\mathcal{L}_{[x],[y]},\mathcal{L}_{[x],[z]},\mathcal{L}_{[y],[z]}\inn \mathcal{C}(m,\mathbb{M})$. It follows from the pigeonhole principle that at least one of the following lines $\mathcal{L}_{[x],[y]},\mathcal{L}_{[x],[z]},\mathcal{L}_{[y],[z]}$ must be a chord of CoP $\mathcal{C}(n,\mathbb{M})$. We deduce $\mathcal{O}_n\cap\mathcal{O}_m\neq \emptyset$.
\end{proof}

\begin{theorem}\label{isomor1}
Let $\mathcal{C}(n,\mathbb{M})$ and $\mathcal{C}(m,\mathbb{M})$ be two parent CoPs that admit an aligned embedding. We have
\[
\mathcal{C}(n,\mathbb{M})\cong \mathcal{C}(m,\mathbb{M})
\]
with a high degree.
\end{theorem}

\begin{proof}
Without loss of generality, we assume that $\mathcal{C}(n,\mathbb{M})\subset \mathcal{C}(m,\mathbb{M})$. By Theorem \ref{T-embedding_child}, all children of CoP $\mathcal{C}(n,\mathbb{M})$ are children of CoP $\mathcal{C}(m,\mathbb{M})$. Hence $\hat{\mathcal{O}}_n\subset\hat{\mathcal{O}}_m$ and therefore
\[
\frac{\vert\hat{\mathcal{O}}_n\cap\hat{\mathcal{O}}_m\vert}{\vert\hat{\mathcal{O}}_n\vert}=1.
\]
\end{proof}

\subsection{Compatible and Incompatible Circles of Partition}

In this section, we introduce the notion of \emph{compatibility} and \emph{incompatibility} of circle of partitions.

\begin{definition}
Let $\mathbb{M}\subseteq \mathbb{N}$, and let $\mathcal{C}(n,\mathbb{M})$ and $\mathcal{C}(m,\mathbb{M})$ be two CoPs. We say that CoP $\mathcal{C}(n,\mathbb{M})$ and CoP $\mathcal{C}(m,\mathbb{M})$ are \emph{compatible} if there exists some CoP $\mathcal{C}(r,\mathbb{M})$ satisfying
\begin{align}
    \mathcal{C}(n,\mathbb{M})\cup\mathcal{C}(m,\mathbb{M})\subseteq\mathcal{C}(r,\mathbb{M})\nonumber
\end{align}
such that for each $[x]\in \mathcal{C}(n,\mathbb{M})\cup \mathcal{C}(m,\mathbb{M})$ with $2x\neq n$, there exists some $[y]\in\mathcal{C}(n,\mathbb{M})\cup\mathcal{C}(m,\mathbb{M})$ so that 
\begin{align}
    \mathbb{L}_{[x],[y]}\inn \mathcal{C}(r,\mathbb{M}).\nonumber
\end{align}
We denote the compatibility by $\mathcal{C}(n,\mathbb{M})\diamond \mathcal{C}(m,\mathbb{M})$. We call CoP $\mathcal{C}(r,\mathbb{M})$ the \emph{cover} of this compatibility.
\end{definition}
\bigskip

\begin{proposition}\label{P-compatible_embedding}
Let $\mathcal{C}(n,\mathbb{M})$ and $\mathcal{C}(m,\mathbb{M})$ be any two CoPs that admit an aligned embedding. We have $\mathcal{C}(n,\mathbb{M})\diamond\mathcal{C}(m,\mathbb{M})$.
\end{proposition}

\begin{proof}
 Without loss of generality, we assume that
\[
\mathcal{C}(n,\mathbb{M})\subset\mathcal{C}(m,\mathbb{M}).
\]
This implies
\begin{align*}
\mathcal{C}(n,\mathbb{M})\cup\mathcal{C}(m,\mathbb{M})
&=\mathcal{C}(m,\mathbb{M})\mbox{ and therefore}\\
\mathcal{C}(n,\mathbb{M})&\diamond\mathcal{C}(m,\mathbb{M}).
\end{align*}
\end{proof}
\bigskip

\begin{theorem}\label{rigornocompatible}
Let $\mathbb{M}\subseteq \mathbb{N}$. There are no CoPs of the forms $\mathcal{C}(n,\mathbb{M})$ and $\mathcal{C}(m,\mathbb{M})$ with all axes points concentrated at their center and satisfies $\mathcal{C}(n,\mathbb{M})\cap\mathcal{C}(m,\mathbb{M})=\emptyset$ for $|\mathcal{C}(n,\mathbb{M})|>2$ and $|\mathcal{C}(m,\mathbb{M})|>2$ with 
\begin{align*}
\nu(n,\mathbb{M})\neq\nu(m,\mathbb{M})
\end{align*}
such that
\begin{align}
    \mathcal{C}(n,\mathbb{M})\diamond\mathcal{C}(m,\mathbb{M})\nonumber
\end{align}
with a cover whose axes points are away from the center.
\end{theorem}

\begin{proof}
We suppose that there exists at least a pair of CoPs of the form $\mathcal{C}(n,\mathbb{M})$ and $\mathcal{C}(m,\mathbb{M})$ with $m\neq n$ such that $\mathcal{C}(n,\mathbb{M})\cap \mathcal{C}(m,\mathbb{M})=\emptyset$ for $|\mathcal{C}(n,\mathbb{M})|,|\mathcal{C}(m,\mathbb{M})|>2$ and satisfies
\begin{align*}
\nu(n,\mathbb{M})\neq\nu(m,\mathbb{M}) 
\end{align*}
so that $\mathcal{C}(n,\mathbb{M})\diamond \mathcal{C}(m,\mathbb{M})$. It follows that there exists some CoP $\mathcal{C}(s,\mathbb{M})$ such that 
\begin{align}
    \mathcal{C}(n,\mathbb{M})\cup \mathcal{C}(m,\mathbb{M})\subseteq \mathcal{C}(s,\mathbb{M})\nonumber
\end{align}
so that for each $[x]\in\mathcal{C}(n,\mathbb{M})\cup\mathcal{C}(m,\mathbb{M})$ there exists some $[y]\in \mathcal{C}(n,\mathbb{M})\cup \mathcal{C}(m,\mathbb{M})$ such that 
\begin{align}
    \mathbb{L}_{[x],[y]}\inn\mathcal{C}(s,\mathbb{M}).\nonumber
\end{align}
Under the conditions 
\begin{align*}
\nu(n,\mathbb{M})\neq \nu(m,\mathbb{M}) 
\end{align*}
and 
\begin{align}
    \mathcal{C}(n,\mathbb{M})\cap \mathcal{C}(m,\mathbb{M})=\emptyset \nonumber
\end{align}
we deduce from the pigeonhole principle and the uniqueness of axes CoPs that there exists some $\mathbb{L}_{[x],[y]}\inn\mathcal{C}(s,\mathbb{M})$ such that $[x],[y]\in \mathcal{C}(n,\mathbb{M})$ or $[x],[y]\in \mathcal{C}(m,\mathbb{M})$. Without loss of generality, we assume that $[x],[y]\in \mathcal{C}(n,\mathbb{M})$. By the embedding 
\begin{align}
    \mathcal{C}(n,\mathbb{M})\subset\mathcal{C}(s,\mathbb{M})\nonumber
\end{align}
the line $\mathcal{L}_{[x],[y]}\inn\mathcal{C}(n,\mathbb{M})$ is such that $\mathcal{L}_{[x],[y]}\neq \mathbb{L}_{[x],[y]}\inn \mathcal{C}(n,\mathbb{M})$. It follows that the line $\mathcal{L}_{[x],[y]}$ must be a chord in $\mathcal{C}(n,\mathbb{M})$ and $\mathcal{C}(s,\mathbb{M})$ must be a \emph{child} of the parent $\mathcal{C}(n,\mathbb{M})$. Now, we locate all the remaining \emph{chords} $\mathcal{L}_{[u],[v]}\neq \mathcal{L}_{[x],[y]}$ in the parent $\mathcal{C}(n,\mathbb{M})$. We claim that each chord $\mathcal{L}_{[u],[v]}$ must be an axis of the child $\mathcal{C}(s,\mathbb{M})$. We assume to the contrary that some chord $\mathcal{L}_{[u],[v]}\inn \mathcal{C}(n,\mathbb{M})$ is also a chord in the child $\mathcal{C}(s,\mathbb{M})$. This implies that there exist some axes
\begin{align}
    \mathbb{L}_{[u],[a]},\mathbb{L}_{[v],[b]}\inn \mathcal{C}(n,\mathbb{M}).\nonumber
\end{align}
Using underlying embedding, we deduce that the lines
\begin{align}
    \mathcal{L}_{[u],[a]},\mathcal{L}_{[v],[b]}\nonumber
\end{align}
cannot be axes of the CoP $\mathcal{C}(s,\mathbb{M})$ so that $\mathcal{L}_{[u],[a]}$ and $\mathcal{L} _{[v],[b]}$ are chords in $\mathcal{C}(s,\mathbb{M})$ with 
\begin{align}
\Gamma([u],[b])=\Gamma([v],[a]).\label{chord1}
\end{align}
It follows that at least one of $\mathcal{L}_{[u],[b]}$ and $\mathcal{L}_{[v],[a]}$ must be chords in CoP $\mathcal{C}(s,\mathbb{M})$. Otherwise, it would mean both lines $\mathcal{L}_{[u],[b]}=\mathbb{L}_{[u],[b]}\inn \mathcal{C}(s,\mathbb{M})$ and $\mathcal{L}_{[v],[a]}=\mathbb{L}_{[v],[a]}\inn \mathcal{C}(s,\mathbb{M})$, which in relation to \eqref{chord1} cannot hold for axes points of CoPs. Without loss of generality, we assume that $\mathcal{L}_{[u],[b]}$ is a chord then so is $\mathcal{L}_{[v],[a]}$ under the condition $\mathbb{L}_{[u],[a]},\mathbb{L}_{[v],[b]}\inn \mathcal{C}(n,\mathbb{M})$. Otherwise, it would imply that the chord $\mathcal{L}_{[a],[v]}$ must be an axis of $\mathcal{C}(s,\mathbb{M})$. Since all the axes points of $\mathcal{C}(n,\mathbb{M})$ are concentrated around the center, it implies that 
\begin{align}
\frac{n}{2}=\frac{a+u}{2}\approx a\quad \mathrm{and} \quad\frac{n}{2}=\frac{a+u}{2}\approx u\label{conc1}
\end{align}
and 
\begin{align}
\frac{n}{2}=\frac{b+v}{2}\approx b \quad \mathrm{and} \quad \frac{n}{2}=\frac{b+v}{2}\approx v\label{conc2}
\end{align}
so that we have $a\approx b \approx u\approx v$ and deduce that the co-axis point $[a],[v]$ of the cover CoP $\mathcal{C}(s,\mathbb{M})$ is close to the center by the relation 
\begin{align}
\frac{s}{2}=\frac{a+v}{2}\approx a \approx v.\nonumber
\end{align} 
This contradicts the requirement of proximity of the axes points of the cover $\mathcal{C}(s,\mathbb{M})$. It follows that $\mathcal{L}_{[u],[v]}$ and $\mathcal{L}_{[a],[b]}$ are also chords in $\mathcal{C}(s,\mathbb{M})$ with 
\begin{align}
    \Gamma([u],[v])=\Gamma([a],[b])\label{chord2}
\end{align}
since the lines $\mathbb{L}_{[u],[a]},\mathbb{L}_{[v],[b]}\inn \mathcal{C}(n,\mathbb{M})$ with the embedding $\mathcal{C}(n,\mathbb{M})\subset\mathcal{C}(s,\mathbb{M})$. It follows from \eqref{chord1} and \eqref{chord2} that
\begin{align}
    \mathbb{L}_{[u],[a]},\mathbb{L}_{[v],[b]}\inn \mathcal{C}(s,\mathbb{M})\nonumber
\end{align}
so that $n=u+a=v+b=s$ and $\mathcal{C}(n,\mathbb{M})=\mathcal{C}(s,\mathbb{M})$, contradicting the embedding 
\begin{align}
    \mathcal{C}(n,\mathbb{M})\subset\mathcal{C}(s,\mathbb{M}).\nonumber
\end{align}
Hence, each \emph{chord} in $\mathcal{C}(n,\mathbb{M})$ must be an axis of the \emph{child} CoP $\mathcal{C}(s,\mathbb{M})$. Therefore, the parent has only one child $\mathcal{C}(s,\mathbb{M})$, which is impossible since $|\mathcal{C}(n,\mathbb{M})|>2$.
\end{proof}
\bigskip

\begin{conjecture}\label{compatible conjecture}
Let $\mathcal{C}(n,\mathbb{M})$ and $\mathcal{C}(m,\mathbb{M})$ be parents CoPs with the offspring $\mathcal{O}_n$ and $\mathcal{O}_m$, respectively. We have $\mathcal{C}(n,\mathbb{M})\diamond\mathcal{C}(m,\mathbb{M})$ if and only if there exists some $\mathcal{C}(s,\mathbb{M})\in \mathcal{O}_m$ and $\mathcal{C}(t,\mathbb{M})\in \mathcal{O}_n$ such that 
\begin{align}
    \mathcal{C}(s,\mathbb{M})\diamond \mathcal{C}(t,\mathbb{M}).\nonumber
\end{align}
\end{conjecture}

\noindent
For a CoP there are two possibilities:
\begin{itemize}
\item
The CoP admits an embedding. The Proposition \ref{P-compatible_embedding} holds for the parents and for their children and the conjecture \ref{compatible conjecture} is valid.
\bigskip

\item
The CoP does not admit an embedding. An example for such CoPs is $\mathcal{C}(n,\mathbb{P})$. The following example demonstrates that Conjecture \ref{compatible conjecture} does not hold for such CoPs.
\end{itemize}

\begin{example}
We consider the weights of CoPs 
\begin{align*}
&\Vert\mathcal{C}(16,\mathbb{P})\Vert=\lbrace 3,5,11,13\rbrace\mbox{ and}\\
&\Vert\mathcal{C}(18,\mathbb{P})\Vert=\lbrace 5,7,11,13\rbrace\mbox{ and}\\
&\Vert\mathcal{C}(24,\mathbb{P})\Vert=\lbrace 5,7,11,13,17,19\rbrace
\mbox{ as child of }\mathcal{C}(16,\mathbb{P})\mbox{ and}\\
&\Vert\mathcal{C}(12,\mathbb{P})\Vert=\lbrace 5,7\rbrace
\mbox{ as child of }\mathcal{C}(18,\mathbb{P}).
\end{align*}
We obtain
\begin{align*}
\mathcal{C}(24,\mathbb{P})\cup\mathcal{C}(12,\mathbb{P})
&=\mathcal{C}(24,\mathbb{P})\mbox{ and therefore}\\
\mathcal{C}(24,\mathbb{P})&\diamond\mathcal{C}(12,\mathbb{P}).
\end{align*}
On the other hand, we have
\begin{align*}
&\Vert\mathcal{C}(16,\mathbb{P})\cup\mathcal{C}(18,\mathbb{P})\Vert=
\lbrace 3,5,7,11,13\rbrace.
\end{align*}
Because $3,5,7$ are the only primes in the CoP with distance two apart, there exists no CoP $\mathcal{C}(n,\mathbb{P})$ for which the weights of the remaining points have a distance of $2$ between each other. Hence $\mathcal{C}(16,\mathbb{P})$ and $\mathcal{C}(18,\mathbb{P})$ are not compatible although they have children that are compatible.
\end{example}

One can obtain a CoP from the union of CoP $\mathcal{C}(16,\mathbb{P})$ and CoP $\mathcal{C}(18,\mathbb{P})$ when we remove the point $[3]$ or $[7]$ from the union. In the first case, we get CoP $\mathcal{C}(18,\mathbb{P})$ and in the second case CoP $\mathcal{C}(16,\mathbb{P})$. In both cases, we have the \emph{weak compatibility} $\mathcal{C}(16,\mathbb{P})\circ\mathcal{C}(18,\mathbb{P})$.   

\begin{definition}
Let $\mathbb{M}\subseteq \mathbb{N}$, and let $\mathcal{C}(n,\mathbb{M})$ and $\mathcal{C}(m,\mathbb{M})$ be two CoPs. We say that CoP $\mathcal{C}(n,\mathbb{M})$ and CoP $\mathcal{C}(m,\mathbb{M})$ are \emph{weakly compatible} if there exist some CoP $\mathcal{C}(r,\mathbb{M})$ and a point $[z]\in\mathcal{C}(n,\mathbb{M})\cup\mathcal{C}(m,\mathbb{M})$ satisfying
\begin{align*}
    \mathcal{C}(n,\mathbb{M})\cup\mathcal{C}(m,\mathbb{M})\setminus\lbrace[z]\rbrace\subseteq\mathcal{C}(r,\mathbb{M})
\end{align*}
such that for each $[x]\in\mathcal{C}(n,\mathbb{M})\cup\mathcal{C}(m,\mathbb{M})$ with $2x\neq n$, there exists some $[y]\in \mathcal{C}(n,\mathbb{M})\cup\mathcal{C}(m,\mathbb{M})$ such that
\begin{align}
    \mathbb{L}_{[x],[y]}\inn \mathcal{C}(r,\mathbb{M}).\nonumber
\end{align}
We denote the weak compatibility by $\mathcal{C}(n,\mathbb{M})\circ\mathcal{C}(m,\mathbb{M})$. We call CoP $\mathcal{C}(r,\mathbb{M})$ the \emph{cover} of this compatibility.
\end{definition}

The following is an example of weakly compatible CoPs:
\begin{align*}
\Vert\mathcal{C}(28,\mathbb{P})\Vert&=\lbrace 5,11,17,23\rbrace\mbox{ and}\\
\Vert\mathcal{C}(30,\mathbb{P})\Vert&=\lbrace 7,11,13,17,19,23\rbrace
\mbox{ and}\\
\mathcal{C}(28,\mathbb{P})\cup\mathcal{C}(30,\mathbb{P})&\setminus\lbrace[5]\rbrace=\mathcal{C}(30,\mathbb{P})\\
\mbox{and }&\mbox{therefore}\\
\mathcal{C}(28,\mathbb{P})&\circ\mathcal{C}(30,\mathbb{P}).
\end{align*}
\bigskip

The conjecture \ref{compatible conjecture} could have several consequences if it is valid. Progress on this conjecture would require an expansion of the notion of compatibility.

\section{Extended Circles of Partition}

In this section, we introduce the notion of \emph{extended} CoPs. 

\begin{definition}\label{D-extended_CoP}
Let $\mathbb{M}\subset\mathbb{N}$, and let $\mathcal{C}(n)$ be a CoP with $\mathbb{N}$ as a base set. We call the CoP 
\[
\mathcal{C}^*(n,\mathbb{M}):=\lbrace [x]\in\mathcal{C}(n)\mid\lbrace x,n-x\rbrace\cap\mathbb{M}\neq\emptyset, x>2\rbrace
\]
the \emph{extended circle of partition}. We abbreviate this as \textbf{xCoP}. We denote by $\mathcal{O}_n^*(\mathbb{M})$ the extended family of xCoP $\mathcal{C}^*(n,\mathbb{M})$---the collection of all children $\mathcal{C}^*(s,\mathbb{M})$ with the principal axes $\mathbb{L}_{[x],[y]}\inn\mathcal{C}^*(s,\mathbb{M}),s=x+y\neq n$ where $[x]$ or $[y]$ is not the center of the parent xCoP. We denote by $\hat{\mathcal{O}}_n^*(\mathbb{M})$ the union of the parent xCoP with its extended family as a complete extended family. We denote the set of \emph{generators} of a complete extended family by $\mathbb{F}_n(\mathbb{M})$. We call an axis whose axis points are members of $\mathbb{M}$ the full--$\mathbb{M}$ axis and in the other case the half--$\mathbb{M}$ axis. If $\mathbb{M}=\mathbb{P}$, then we get the full prime axis (resp. the half prime axis).
\end{definition}

\begin{proposition}\label{P-ext_subset}
We have $\mathcal{C}(n,\mathbb{M})\subseteq\mathcal{C}^*(n,\mathbb{M})$ for all $n\in\mathbb{N}$ with $\mathcal{C}(n,\mathbb{M})\neq\emptyset$.
\end{proposition}

\begin{proof}
If $x\in\mathbb{M}$ and $n-x\in\mathbb{M}$, then $[x]\in\mathcal{C}(n,\mathbb{M})$ and $[x]\in\mathcal{C}^*(n,\mathbb{M})$. If $x\in\mathbb{M}$ and $[n-x]\not\in\mathbb{M}$, then $[x]\in\mathcal{C}^*(n,\mathbb{M})$ and $[x]\not\in\mathcal{C}(n,\mathbb{M})$.
\end{proof}
\bigskip

\begin{corollary}\label{C-ext_Psubset}
Let $\mathbb{P}$ be the set of all prime numbers. By Proposition \ref{P-ext_subset}, for all $n\in2\mathbb{N}$ with $\mathcal{C}(n,\mathbb{P})\neq\emptyset$
\[
\mathcal{C}(n,\mathbb{P})\subset\mathcal{C}^*(n,\mathbb{P}).
\]
\end{corollary}
\bigskip

\begin{proposition}\label{P-ext_progression}
Let $\mathbb{M}_{a,d}\subset\mathbb{N}$ be as defined in \eqref{Mad}. For all $n\in\mathbb{M}_{2a,d}$, we have
\[
\mathcal{C}(n,\mathbb{M}_{a,d})\equiv\mathcal{C}^*(n,\mathbb{M}_{a,d}).
\]
\end{proposition}

\begin{proof}
Due to $n\in\mathbb{M}_{2a,d}$, it implies that for each member $x$ of $\mathbb{M}_{a,d}$, the summand $n-x$ is a member of $\mathbb{M}_{a,d}$. Therefore, there are no $x\in\mathbb{M}_{a,d}$ with $n-x\notin\mathbb{M}_{a,d}$ and vise versa. These are extensions of $\mathcal{C}(n,\mathbb{M}_{a,d})$. 
\end{proof}
\bigskip

\begin{proposition}\label{P-ext_primes}
Let $n\in 2\mathbb{N}$ be with $ n\geq 8$. The xCoP $\mathcal{C}^*(n,\mathbb{P})$ contains all odd primes $\leq n-3$ and 
\[
\vert\mathcal{C}^*(n,\mathbb{P})\vert\geq\pi(n-3)-1\geq 2.
\]
\end{proposition}

\begin{proof}
The xCoP $\mathcal{C}^*(n,\mathbb{P})$ contains only odd numbers because if $x$ even with $>2$, then $n-x$ is even, and therefore both are not prime. The first member of each such xCoP is $[3]$ with the co-axis point $[n-3]$. Therefore, CoP $\mathcal{C}^*(n,\mathbb{P})$ contains points with weight as prime numbers $p$ that satisfy $3\leq p\leq n-3$. All axis points $[u]\in\mathcal{C}(n)$ and their axis partner $[n-u]$ have weights that are not both primes. This implies that the points $[u]$ and $[n-u]$ are not contained in $\mathcal{C}^*(n,\mathbb{P})$. Because $ \Vert\mathcal{C}(n)\Vert$ contains all positive integers $1\leq x\leq n-1$, we find that $ \Vert\mathcal{C}(n)\Vert $ contains all primes $3\leq p\leq n-3$. Therefore, CoP $\mathcal{C}^*(n,\mathbb{P})$ contains all primes $p$ satisfying $3\leq p\leq n-3$ and therefore
\[
\vert\mathcal{C}^*(n,\mathbb{P})\vert\geq\pi(n-3)-1\geq 2.
\]
\end{proof}

\begin{lemma}\label{L-ext_symmetrCoP}
Let $\mathbb{M}\subseteq\mathbb{N}$, and let $\mathcal{C}^*(n,\mathbb{M})$ be an xCoP. The weights of the xCoP are symmetrically distributed around $\frac{n}{2}$. 
\end{lemma}

\begin{proof}
This is true because
\[
x_{i+1}-x_i = n-n+x_{i+1}-x_i=(n-x_i)-(n-x_{i+1})
\mbox{ for } i=1,2,\ldots,\left\lfloor\frac{k}{2}\right\rfloor
\]
if $x_1,x_2,\ldots,x_k$ are the weights of the points of xCoP.
\end{proof}

\begin{lemma}\label{L-ext_symmetrFam}
Let $\mathbb{M}\subseteq\mathbb{N}$, and let $\mathcal{C}^*(n,\mathbb{M})$ be an xCoP. The generators of its children are symmetrically distributed around the generator of the parent xCoP. 
\end{lemma}

\begin{proof}
In analogy with the proof of Proposition \ref{breaker}, we get for two arbitrary axes 
\[
\mathbb{L}_{[u],[v]},\mathbb{L}_{[x],[y]}\inn\mathcal{C}^*(n,\mathbb{M}), u<x
\]
of the xCoP the following children generators
\begin{align*}
s_1=u+x<s_2=u+y<~ &n<t_2=v+x<t_1=v+y\\ 
\mbox{and the }&\mbox{following distances}\\
s_2-s_1=u+y-u-x &\mbox{ and }t_1-t_2=v+y-v-x\\
\mbox{and }&\mbox{hence}\\
s_2-s_1=y &-x=t_1-t_2.
\end{align*}
Because the chosen axes are arbitrary and the parent generator is always located in the middle of such inequalities the distances between the children generators are symmetrically distributed around $n$.
\end{proof}

\begin{corollary}\label{C-ext_evenChildren}
From the proof of Lemma \ref{L-ext_symmetrFam}, we conclude that the number of children of any xCoP is even and hence every complete extended family has an odd number of members.
\end{corollary}

\begin{theorem}\label{T-ext_embFams}
Let $m,n\in 2\mathbb{N}$, and let $\mathcal{C}^*(m,\mathbb{P})$ and $\mathcal{C}^*(n,\mathbb{P})$ be two xCoPs with $m<n$ and $m\geq 16$. If $\mathcal{C}(k,\mathbb{P})\neq\emptyset$ for $k\in 2\mathbb{N}\mid 8\leq k <n$, then $\mathbb{F}_n(\mathbb{P})$ contains all even integers $x$ satisfying $8\leq x\leq 2n-8$ and
\[
\mathbb{F}_m(\mathbb{P})\subset\mathbb{F}_n(\mathbb{P}).
\]
Furthermore, we have $\vert\mathbb{F}_n(\mathbb{P})\vert=n-7$. This means that xCoP $\mathcal{C}^*(n,\mathbb{P})$ has $n-8$ children.
\end{theorem}

\begin{proof}
Because the first three odd numbers $\geq 3$ all are primes, we get the first three points that are not on a common axes or on a degenerated axis of all xCoPs to be $[3],[5],[7]$ for $m\geq 16$. 
Hence, the first three generators of children of all such xCoPs are
\[
f_1=3+5=8,
f_2=3+7=10\mbox{ and }
f_3=5+7=12.
\]
By the lemma \ref{L-ext_symmetrFam}, we only have to prove that all even integers $8\leq x<n$ are members of $\mathbb{F}_n(\mathbb{P})$. We let $t$ be the cardinality of $\mathbb{F}_n(\mathbb{P})$ and
\[
s:=\frac{t+1}{2}.
\]
By Corollary \ref{C-ext_evenChildren}, we find that $t$ is an odd number and $\mathbb{F}_n(\mathbb{P})$ has the parent generator $n$ in the middle and that there are equal even number of children generators left and right of the parent generator. Therefore $n$ is the $s^{th}$ element of $\mathbb{F}_n(\mathbb{P})$
\[
f_s=n.
\]
The point $[n-5]$ is a member of xCoP $\mathcal{C}^*(n,\mathbb{P})$ and therefore $\mathcal{L}_{[3][n-5]}$ is not an axis. Hence $\mathcal{L}_{[3][n-5]}$ must be the principal axis for the child $\mathcal{C}^*(n-2,\mathbb{P})$. Hence the greatest member of $\mathbb{F}_n(\mathbb{P})$ with $\leq n$ is   
\[
f_{s-1}=n-2.
\]
Because of $\mathcal{C}(k,\mathbb{P})\neq\emptyset$ for $k\in 2\mathbb{N}$ satisfying $8\leq k<n$, all such even numbers $k$ have at least one representation as sum of two primes. By Proposition \ref{P-ext_primes}, the xCoP $\mathcal{C}^*(n,\mathbb{P})$ contains all primes $p$ satisfying $3\leq p\leq n-3$. Therefore every even number $r$ satisfying $8\leq r\leq n-2$ has a representation as sum of two weights of xCoP $\mathcal{C}^*(n,\mathbb{P})$ whose points do not form an axis of xCoP $\mathcal{C}^*(n,\mathbb{P})$. Therefore these even numbers are generators of children of the parent xCoP $\mathcal{C}^*(n,\mathbb{P})$, and hence the elements $f_1=8,f_2=10,\ldots,f_{s-1}=n-2\in\mathbb{F}_n(\mathbb{P})$.
There are
\[
s-1=\frac{n-2-8}{2}+1=\frac{n}{2}-4
\]
members. We deduce
\[
t=\vert\mathbb{F}_n(\mathbb{P})\vert=2s-1=n-6-1=n-7.
\]
By the lemma \ref{L-ext_symmetrFam}, the members from $f_{s+1}$ up to $f_t$ are the $s-1$ even numbers $>n$. For the greatest children generator, we obtain
\[
f_t =n-3+n-5=2n-8
\]
Therefore all even numbers $r$ satisfying $8\leq r\leq 2n-8$ are the members of $\mathbb{F}_n(\mathbb{P})$. We deduce $\mathbb{F}_m(\mathbb{P})\subset\mathbb{F}_n(\mathbb{P})$ since $\vert\mathbb{F}_m(\mathbb{P})\vert<\vert\mathbb{F}_n(\mathbb{P})\vert$. We deduce that xCoP $\mathcal{C}^*(n,\mathbb{P})$ has $n-8$ children. We remark that there are no even $r$ with $8\leq r\leq 2n-8$ that is a child generator.
\end{proof}

\begin{corollary}\label{C-ext_isomorphic}
Due to $\mathbb{F}_m(\mathbb{P})\subset\mathbb{F}_n(\mathbb{P})$ for $m<n$ all xCoPs $\mathcal{C}^*(n,\mathbb{P})$ with $n\geq 12$
and $\mathcal{C}(k,\mathbb{P})\neq\emptyset$ satisfying $8\leq k<n$ implies
\[
\hat{\mathcal{O}}^*_m\subset\hat{\mathcal{O}}^*_n\mbox{ if }m<n,
\]
which means that all such two xCoPs are isomorphic with a high degree (see the definition \ref{D-isomor}). 
\end{corollary}
\bigskip

\begin{corollary}\label{C-ext_connected}
 All xCoPs $\mathcal{C}^*(n,\mathbb{P})$ for $n\geq 6$ are connected (see Proposition \ref{P-common_points}).
\end{corollary}
\bigskip

\begin{proposition}\label{T-ext_axis}
Let $n\in 2\mathbb{N}$, and let $\mathcal{C}^*(n,\mathbb{P})$ be an xCoP with $n\geq 8$ and $\mathcal{C}(n,\mathbb{P})\neq\emptyset$. There exists an axis $\mathbb{L}_{[x],[n-x]}\inn\mathcal{C}^*(n,\mathbb{P})$ such that $\mathbb{L}_{[x],[n-x]}$ is also an axis of $\mathcal{C}(n,\mathbb{P})$.
\end{proposition}

\begin{proof}
This is obvious and can be verified from the development of this section.
\end{proof}
\bigskip

\section{The Asymptotic Binary Goldbach and Lemoine Conjectures}

\subsection{The Squeeze Principle}

In this section, we introduce the \emph{squeeze principle} and its consequences if the set of all odd prime numbers is the base set of CoPs.

\begin{theorem}[The squeeze principle]\label{L_squeeze principle}
Let $\B\subset\M\subseteq\N$, and suppose that $\C(m,\B)$ and $\C(m+t,\B)\neq\emptyset$ for some $t\geq 4$. Suppose that there exist $\mathbb{L}_{[x],[y]}\inn \mathcal{C}(m+t,\mathbb{M})$ with $x\in \mathbb{B}$ and $x<y$ such that 
\begin{align}\label{E_maxcovered}
y>w= \max\set{u\in||\C(m,\M)||\mid u\in \B}>m-x.
\end{align}
There exists a generator $s$ with $m<s<m+t$ such that $\C(s,\B)\neq\emptyset$.
\end{theorem}

\begin{proof}
By \eqref{E_maxcovered}, we have $w\in\B$. By hypothesis, the axis $\mathbb{L}_{[x],[y]}\inn \mathcal{C}(m+t,\mathbb{M})$ exists with $x\in \mathbb{B}$ such that $m-w<x<y$. We deduce
\begin{align}
m=w+(m-w)<\ul{w+x}&=w+(m+t-y)=m+t+(w-y)\nonumber\\
&<m+t \mbox{, since }y>w \label{E_inequality}
\end{align}
and $m-w<x=m+t-y$ implies $y-w<t$. With $s=\ul{w+x}$, there is an axis $\LL_{[x],[w]}\inn\C(s,\B)$ and it follows that $\mathcal{C}(s,\mathbb{B})\neq \emptyset$ with $m<s<m+t$. 
\end{proof}
\bigskip

Theorem \ref{L_squeeze principle} may be viewed as a basic apparatus for studying the possibility of partitioning numbers of a particular parity into components that belong to a special subset of positive integers. It works by choosing two non-empty CoPs with the same base set and finding further non-empty CoPs with generators trapped in between these two generators. This principle can be used with care to study the broader question that concerns the feasibility of partitioning numbers with each summand belonging to the same subset of positive integers. We state the following proposition as an outgrowth of Theorem \ref{L_squeeze principle}.

\begin{proposition}[The interval binary Goldbach partition detector]\label{P_The interval binary Goldbach partition detector}
Let $\mathbb{P}$ be the set of all prime numbers and let $\mathcal{C}(m,\mathbb{P}),\mathcal{C}(m+t,\mathbb{P})\neq \emptyset$ for some $t\geq 4$. Suppose that there exist $\mathbb{L}_{[x],[y]}\inn \mathcal{C}(m+t,\mathbb{N})$ with $x\in\p$ and $x<y$ such that
\begin{align}\label{E_maxcoveredP}
y>w= \max\set{u\in||\C(m,\N)||\mid u\in \p}>m-x.
\end{align}
There exists a generator $s$ with $m<s<m+t$ such that $\C(s,\p)\neq \emptyset$.
\end{proposition}

\begin{proof}
This is a simple consequence of Theorem \ref{L_squeeze principle} by taking $\M=\N$ and $\B=\p$.
\end{proof}

\begin{proposition}[Interval Goldbach partition]\label{P_intervalGoldbach}
Let $\p$ be the set of all prime numbers and $\C(m,\p), \C(m+t,\p)\neq\emptyset$ be for some $t\geq 4$. If $m-1\in \p$, then there exists some $s\equiv 0\pmod 2$ with $m<s<m+t$ such that $\C(s,\p)\neq\emptyset$.
\end{proposition}

\begin{proof}
By the requirements $\C(m,\p),\C(m + t,\p)\neq\emptyset$ for $t \geq 4$ and with $w$ by virtue of \eqref{E_maxcoveredP}, we choose
$\LL_{[3],[y]}\inn \C(m + t,\N)$ so that $w=m-1$ and $y>w$
since $y=m+t-3>m$ for $t\geq 4$ and $m-1\in\p$. We deduce
\begin{align*}
y-w=y-(m-1)\leq(m+t-3)-(m-1)<t
\end{align*}
and the conditions in Proposition \ref{P_The interval binary Goldbach partition detector} are satisfied, so that there exist some $s\equiv 0\pmod 2$ with $m<s<m+t$ such that $\C(s, \p)\neq\emptyset$, f.i. $s=3+m-1=m+2$ with $\LL_{[3],[m-1]}\inn\C(m+2,\p)$.
\end{proof}

\begin{proposition}\label{P_finite Goldbach in an interval}
Let $\mathbb{P}$ be the set of all prime numbers and $\mathcal{C}(m,\mathbb{P}),\mathcal{C}(m+t,\mathbb{P})\neq \emptyset$ be for some $t\geq 4$ such that $m-1\in \mathbb{P}$. There are finitely many $s\equiv 0\pmod 2$ with $m<s<m+t$ such that $\C(s,\p)\neq \emptyset$.
\end{proposition}

\begin{proof}
The result is obtained by iterating on the generators $s\equiv 0\pmod 2$ with $m<s<m+t$ such that $\C(s,\p)\neq \emptyset$.
\end{proof}

\begin{theorem}[Conditional Goldbach]\label{T_Conditional Goldbach}
Let $\mathbb{P}$ be the set of all prime numbers and $m\in 2\mathbb{N}$ such that $\mathcal{C}(m,\mathbb{P})\neq \emptyset$ for a \textbf{sufficiently} large $m$ . If for \textbf{all} $t\geq 4$ there exists $\mathbb{L}_{[x],[y]}\inn \mathcal{C}(m+t,\mathbb{N})$ with $x\in \mathbb{P}$ and $x<y$ such that
\[
y> w= \max\set{u\in||\C(m,\N)||\mid u\in \p}>m-x,
\]
then there are CoPs $\mathcal{C}(s,\mathbb{P})\neq \emptyset$ for all (sufficiently large) $s\in 2\mathbb{N}$ with $s>m$.
\end{theorem}

\begin{proof}
It is known that there are infinitely many even numbers that can be written as the sum of two primes, so that for \textbf{sufficiently} large $m\in 2\mathbb{N}$ with $\mathcal{C}(m,\mathbb{P})\neq \emptyset$, then $t\geq 4$ can be \textbf{arbitrarily} chosen large such that $\mathcal{C}(m+t,\p) \neq \emptyset$. By the requirements and Proposition \ref{P_The interval binary Goldbach partition detector}, there must exist some $s\equiv 0\pmod 2$ with $m<s<m+t$ such that $\mathcal{C}(s,\mathbb{P})\neq \emptyset$. Now, we continue our arguments on the intervals of the generators $[m,s]$ and $[s,s+r]$. If there exist some $u,v\in 2\mathbb{N}$ such that $m<u<s$ and $s<v<s+r$, then we repeat the argument under the requirements (for arbitrary $t$) to deduce that $\mathcal{C}(u,\mathbb{P})\neq \emptyset$ and $\mathcal{C}(v,\mathbb{P})\neq \emptyset$. We can iterate the process so long as there exist some even generators trapped in the following sub-intervals of generators $[m,u],[u,s],[s,v],[v,v+r]$ where $v+r=m+t$ for $t\geq 4$. Since $t$ can be arbitrarily chosen so that $\mathcal{C}(m+t,\p) \neq \emptyset$, the assertion follows immediately.
\end{proof}
\bigskip

Now, we use the \emph{squeeze principle} to solve the Lemoine conjecture in analogy to its use for the binary Goldbach conjecture above.

\begin{proposition}[The first interval Lemoine partition detector]\label{P_The first interval Lemoine partition detector}
Let $\mathbb{P}$ and $2\mathbb{P}$ be the set of all prime numbers and their doubles, respectively, and let $\mathcal{C}(m,\mathbb{P}\cup 2\mathbb{P}),\mathcal{C}(m+t,\mathbb{P}\cup 2\mathbb{P})\neq \emptyset$ be for some $t\geq 4$. Suppose that there exist $\mathbb{L}_{[x],[y]}\inn \mathcal{C}(m+t,\mathbb{N})$ with $x\in\p$ and $x<y$ such that
\begin{align}\label{E_maxcoveredP2}
y>w=\max\set{u\in||\C(m,\N)||\mid u\in \p\cup 2\p}\in 2\mathbb{P}>m-x.
\end{align}
There exists a generator $s$ with $m<s<m+t$ such that $\C(s,\p\cup 2\p)\neq \emptyset$.
\end{proposition}

\begin{proof}
This is a consequence of Theorem \ref{L_squeeze principle} when we take $\M=\N$ and $\B=\p\cup 2\p$.
\end{proof}

\begin{proposition}[The second interval Lemoine partition detector]\label{P_The second interval Lemoine partition detector}
Let $\mathbb{P}$ and $2\mathbb{P}$ be the set of all prime numbers and their doubles, respectively, and $\mathcal{C}(m,\mathbb{P}\cup 2\mathbb{P}),\mathcal{C}(m+t,\mathbb{P}\cup 2\mathbb{P}) \neq \emptyset$ by $t\geq 4$. Suppose that there exist $\mathbb{L}_{[x],[y]}\inn \mathcal{C}(m+t,\mathbb{N})$ with $x\in 2\p$ and $x<y$ such that
\begin{align}\label{E_maxcoveredP3}
y>w= \max\set{u\in||\C(m,\N)||\mid u\in \p\cup 2\p}\in \mathbb{P}>m-x.
\end{align}
There exists a generator $s$ with $m<s<m+t$ such that $\C(s,\p\cup 2\p)\neq \emptyset$.
\end{proposition}

\begin{proof}
The proof is the same as in Proposition \ref{P_The first interval Lemoine partition detector}.
\end{proof}

\begin{theorem}[Conditional Lemoine]\label{T_Conditional Lemoine}
Let $\mathbb{P}$ and $2\mathbb{P}$ be the set of all prime numbers and their doubles, respectively, and $m\in 2\mathbb{N}+1$ be such that $\mathcal{C}(m,\mathbb{P})\neq \emptyset$ for a \textbf{sufficiently} large $m$. If for \textbf{all} $t\geq 4$ there exists $\mathbb{L}_{[x],[y]}\inn \mathcal{C}(m+t,\mathbb{N})$ with $x\in \mathbb{P}$ and $x<y$ such that
\[
y>w=\max\set{u\in||\C(m,\N)||\mid u\in \p\cup 2\p}\in 2\p>m-x,
\]
or there exist $\mathbb{L}_{[x],[y]}\inn \mathcal{C}(m+t,\mathbb{N})$ with $x\in 2\mathbb{P}$ and $x<y$ such that
\[
y> w= \max\set{u\in||\C(m,\N)||\mid u\in \p\cup 2\p}\in \p>m-x
\]
then there are CoPs $\mathcal{C}(s,\mathbb{P}\cup 2\mathbb{P})\neq \emptyset$ for all (sufficiently large) $s\in 2\mathbb{N}+1$ with $s>m$.
\end{theorem}

\begin{proof}
It is known that there are infinitely many odd numbers that can be written as the sum of a prime and a double of a prime, so that for a \textbf{sufficiently} large $m\in 2\mathbb{N}+1$ with $\mathcal{C}(m,\mathbb{P}\cup 2\mathbb{P})\neq \emptyset$, then $t\geq 4$ can be \textbf{arbitrarily} chosen large such that $\mathcal{C}(m+t,\p)\neq \emptyset$. By the requirements and Propositions \ref{P_The first interval Lemoine partition detector} and \ref{P_The second interval Lemoine partition detector}, there must exist some $s\equiv 1\pmod 2$ with $m<s<m+t$ such that $\mathcal{C}(s,\mathbb{P}\cup 2\mathbb{P})\neq \emptyset$. Now, we continue our arguments on the intervals of the generators $[m,s]$ and $[s,s+r]$. If there exist some $u,v\in 2\mathbb{N}+1$ such that $m<u<s$ and $s<v<s+r$, then we repeat the argument under the requirements (for arbitrary $t$) to deduce that $\mathcal{C}(u,\mathbb{P}\cup 2\mathbb{P})\neq \emptyset$ and $\mathcal{C}(v,\mathbb{P}\cup 2\mathbb{P})\neq \emptyset$. We can iterate the process so long as there exist some odd generators (numbers) trapped in the following sub-intervals of generators $[m,u],[u,s],[s,v],[v,v+r]$ where $v+r=m+t$ for $t\geq 4$. Since $t$ can be arbitrarily chosen so that $\mathcal{C}(m+t,\p\cup 2\p)\neq \emptyset$, the assertion follows immediately.
\end{proof}

\subsection{Application to the Binary Goldbach Conjecture}

In this section, we apply the \emph{squeeze principle} to study the asymptotic version of the binary Goldbach conjecture. Despite the Estermann proof from 1938 (see, e.g, \cite{estermann1938goldbach}, \cite{chudakov1938goldbach}) that the binary Goldbach conjecture is true for almost all positive integers, and the work of Chen that every sufficiently large even number can be written as a prime and almost prime (see, e.g, \cite{chen2002representation}) we can use our elementary tool to independently establish and prove the binary Goldbach conjecture in an asymptotic sense. We gather the following established facts that will feature prominently in our arguments.

\begin{lemma}[The prime number theorem]\label{L_the prime number theorem}
Let $\pi(m)$ denote the number of prime numbers less than or equal to $m$ and $p_{\pi(m)}$ denote the $\pi(m)^{th}$ prime number. We have the asymptotic relation
\begin{align}
p_{\pi(m)}\sim \pi(m)\log \pi(m) \quad \text{and} \quad \pi(m)\sim \frac{m}{\log m}.\nonumber
\end{align}
\end{lemma}

In keeping with the notation of the previous section, we write
\begin{align}\label{E_w-max-p}
w=\max\set{u\in||\C(m,\N)|| \mid u\in\p}=p_{\pi(m)}.
\end{align}
and set
$$
w''=2p_{\pi(\frac{m}{2})}.
$$

\bigskip

\begin{lemma}[Bertrand's postulate]\label{L_Bertrand's postulate}
There exists a prime number in the interval $(k,2k)$ for all $k>1$.
\end{lemma}
\bigskip

\begin{lemma}[The little lemma]\label{L_the little lemma}
Let $\mathbb{P}$ be the set of all prime numbers and $m\in \mathbb{N}$ be \textbf{sufficiently} large such that $\mathcal{C}(m,\mathbb{P})\neq \emptyset$. There exists some $h:=h(m)=o(\log m)$ such that
\begin{align*}
m-w''\sim m-w\sim h(m)\frac{m}{\log m}.
\end{align*}
For all $x\in \mathbb{P}$ satisfying 
$$
m-w<x\leq 3(m-w)
$$ 
the inequality
\begin{align*}
0\lesssim |w-(m+t-x)|\lesssim t
\end{align*}
hold for $t\geq 4$.
\end{lemma}

\begin{proof}
The first assertion can be deduced from the prime number theorem. By the Bertrand postulate, there exists a prime number $x\in \mathbb{P}$ satisfying $m-w<x\leq 3(m-w)$, since 
\begin{align}
w=\max\set{u\in||\C(m,\N)|| \mid u\in\p}=p_{\pi(m)}.
\end{align}
imply $m-w\geq 1$.
We deduce
\begin{align*}
m+t-x&\gtrsim m+t-3h(m)\frac{m}{\log m}\\
&=m\left(1-3\frac{h(m)}{\log m}\right)+t\\
&\sim m+t>p_{\pi(m)}=w
\end{align*}
and 
\begin{align*}
|w-(m+t-x)|&=|m+t-x-w|\\
&<|m+t-(m-w)-w|\\
&=t\nonumber
\end{align*}
for $t\geq 4$.
\end{proof}
\bigskip

We are now ready to prove the binary Goldbach conjecture for all \textbf{sufficiently} large even numbers. The following result is a combination of ideas developed in this paper. 

\begin{theorem}[Asymptotic Goldbach theorem]\label{T_Asymoptotic Goldbach theorem}
Every \textbf{sufficiently} large even number can be written as the sum of two prime numbers.
\end{theorem}

\begin{proof}
The claim is equivalent to the following statement:
\begin{center}
For every sufficiently large even number $n$, we have $\C(n,\p)\neq\emptyset$.
\end{center}
It is known that there are infinitely many even numbers $m>0$ with $\mathcal{C}(m,\mathbb{P})\neq \emptyset$. We choose a \textbf{sufficiently} large $m\in 2\N$ such that $\mathcal{C}(m,\mathbb{P})\neq \emptyset$ and choose $t\geq 4$ such that $\mathcal{C}(m+t,\mathbb{P})\neq \emptyset$. We set
\begin{align}
w=\max\set{u\in||\C(m,\N)|| \mid u\in\p}=p_{\pi(m)}\nonumber
\end{align}
and choose a prime number $x\leq 3(m-w)$ such that $x>m-w$, since by the Bertrand postulate (Lemma \ref{L_Bertrand's postulate}) there exists a prime number $x$ such that $x\in (k,2k)$ for every $k>1$. We obtain for the axis partner $[y]$ of the axis point $[x]$ of $\LL_{[x],[y]}\inn\C(m+t,\N)$ the inequality 
\begin{align*}
y=m+t-x&\gtrsim m+t-3h(m)\frac{m}{\log m}\\
&=m\left(1-3\frac{h(m)}{\log m}\right)+t\\
&\sim m+t>p_{\pi(m)}=w
\end{align*}
for $t\geq 4$ and by the lemma \ref{L_the little lemma}, we deduce
\begin{align*}
h(m)\frac{m}{\log m}\sim m-w<x
\end{align*}
with $h(m)=o(\log m)$ and 
\begin{align*}
\abs{y-w}=\abs{(m+t-x)-w}=\abs{m-w+t-x}\lesssim\abs{x+t-x}= t.
\end{align*}
The requirements in Theorem \ref{T_Conditional Goldbach} are \textbf{asymptotically} fulfilled 
with
\[
y\gtrsim w\mbox{ and }x\gtrsim m-w\mbox{ and }0\lesssim \abs{y-w}\lesssim t.
\]
The result follows by arbitrarily choosing $t\geq 4$ so that $\mathcal{C}(m+t,\mathbb{P})\neq \emptyset$ and adapting the proof in Theorem \ref{T_Conditional Goldbach}.
\end{proof}
\bigskip

\subsection{Application to the Lemoine Conjecture}

In this section, we apply the \emph{squeeze principle} to study Lemoine conjecture (see, e.g, \cite{levi1963}, \cite{lemoine1894}) in the asymptotic.

\begin{lemma}\label{L_the prime number theorem2}
Let $\pi(m)$ denote the number of prime numbers less than or equal to $m$ and $p_{\pi(m)}$ denote the $\pi(m)^{th}$ prime number. We have the asymptotic relation
\begin{align*}
p_{\pi(m)}\sim 2p_{\pi(\frac{m}{2})}.
\end{align*}
\end{lemma}

\begin{proof}
This is the prime number theorem.
\end{proof}

In keeping with the notation of the previous section, we write
\begin{align}\label{E_w-max-p1}
w=\max\set{u\in||\C(m,\N)||~\mid u\in\p\cup 2\p}=p_{\pi(m)}
\end{align}
provided that $w\in \mathbb{P}$ and 
\begin{align}\label{E_w-max-p2}
w'=\max\set{u\in||\C(m,\N)||~\mid u\in\p\cup 2\p}=2p_{\pi(\frac{m}{2})}
\end{align}
provided that $w'\in 2\mathbb{P}$.
\bigskip

\begin{lemma}[The first little lemma]\label{L_the first little lemma}
Let $\mathbb{P}$ and $2\mathbb{P}$ be the set of all prime numbers and their doubles, respectively, and let $m\in \mathbb{N}$ be \textbf{sufficiently} large such that $\mathcal{C}(m,\mathbb{P}\cup 2\mathbb{P})\neq \emptyset$. Suppose that
\begin{align}
w=\max\set{u\in||\C(m,\N)||~\mid u\in\p\cup 2\p}=p_{\pi(m)}.
\end{align}
There exists some $h:=h(m)=o(\log m)$ such that the asymptotic relation
\begin{align*}
m-w\sim m-w'\sim h(m)\frac{m}{\log m}
\end{align*}
holds. Furthermore, for all $x'\in 2\mathbb{P}$ with $x'=2x$ for $x\in \mathbb{P}$ satisfying 
\[
\frac{1}{2}(m-w')<x<m-w'
\]
the inequality
\begin{align*}
0\lesssim |w-(m+t-x')|\lesssim t
\end{align*}
holds for $t\geq 4$.
\end{lemma}

\begin{proof}
The first assertion can be deduced from the prime number theorem. By the lemma \ref{L_Bertrand's postulate}, there is a prime between $\frac{1}{2}(m-w')$ and $m-w'$, since 
\begin{align}
w=\max\set{u\in||\C(m,\N)||~\mid u\in\p\cup 2\p}=p_{\pi(m)}.
\end{align}
imply $m-w'\geq 3$.
We deduce
\begin{align*}
m+t-x'=m+t-2x&\gtrsim m+t-2h(m)\frac{m}{\log m}\\
&=m\left(1-2\frac{h(m)}{\log m}\right)+t\\
&\sim m+t\geq p_{\pi(m)}=w
\end{align*}
and 
\begin{align*}
|w-(m+t-x')|&\sim |m+t-x'-w'|\\
&\lesssim |m+t-(m-w')-w'|\\
\\
&=t\nonumber
\end{align*}
for $t\geq 4$.
\end{proof}

\begin{lemma}[The second little lemma]\label{L_the second little lemma}
Let $\mathbb{P}$ be the set of all prime numbers and their doubles, respectively, and $m\in \mathbb{N}$ be \textbf{sufficiently} large such that $\mathcal{C}(m,\mathbb{P}\cup 2\mathbb{P})\neq \emptyset$. Suppose that 
\begin{align}
w'=\max\set{u\in||\C(m,\N)||~\mid u\in\p\cup 2\p}=2p_{\pi(\frac{m}{2})}.
\end{align}
There exists some $h:=h(m)=o(\log m)$ such that 
$$
m-w'\sim m-w\sim h(m)\frac{m}{\log m}.
$$
For all $x\in \mathbb{P}$ satisfying 
\[
m-w<x<2(m-w)
\]
we have
\begin{align*}
0\lesssim |w'-(m+t-x)|\lesssim t
\end{align*}
for $t\geq 4$.
\end{lemma}

\begin{proof}
By the lemma \ref{L_Bertrand's postulate}, there is a prime between $m-w$ and $2(m-w)$, since 
\begin{align}
w'=\max\set{u\in||\C(m,\N)||~\mid u\in\p\cup 2\p}=2p_{\pi(\frac{m}{2})}
\end{align}
imply $m-w\geq 2$. We deduce
\begin{align*}
m+t-x&>m+t-2h(m)\frac{m}{\log m}\\
&=m\left(1-2\frac{h(m)}{\log m}\right)+t\\
&\sim m+t\geq p_{\pi(m)}\sim 2p_{\pi(\frac{m}{2})}=w'
\end{align*}
and 
\begin{align*}
|w'-(m+t-x)|&=|m+t-x-w'|\\
&\lesssim |m+t-(m-w')-w'|\\
&=t\nonumber
\end{align*}
for $t\geq 4$.
\end{proof}
\bigskip

We are now ready to prove the Lemoine conjecture for all \textbf{sufficiently} large odd numbers. It is a case-by-case argument and a culmination of ideas developed in this paper. 

\begin{theorem}[Asymptotic Lemoine theorem]\label{T_Asymoptotic Lemoine theorem}
Every \textbf{sufficiently} large odd number can be written as the sum of a prime number and a double of a prime number.
\end{theorem}

\begin{proof}
The claim is equivalent to the following statement:
\begin{center}
For every sufficiently large odd number $n\in 2\mathbb{N}+1$, we have $\C(n,\p\cup 2\p)\neq\emptyset$, since only the sum of an odd and an even number provides an odd number and therefore each axis $\LL_{[x],[y]}\inn\C(m,\p\cup 2\p)$ has an odd and an even axis point.
\end{center}
It is known that there are infinitely many odd numbers $m>0$ with $\mathcal{C}(m,\mathbb{P}\cup 2\mathbb{P})\neq \emptyset$. We choose a \textbf{sufficiently} large $m\in 2\N+1$ such that $\mathcal{C}(m,\mathbb{P}\cup 2\mathbb{P})\neq \emptyset$ and choose $t\geq 4$ such that $\mathcal{C}(m+t,\mathbb{P}\cup 2\mathbb{P})\neq \emptyset$. Now, we distinguish and examine two special cases as follows:
\begin{enumerate}
\item  The case $$w=\max\set{u\in||\C(m,\N)||~\mid u\in\p\cup 2\p}=p_{\pi(m)}$$

\item The case $$w'=\max\set{u\in||\C(m,\N)||~\mid u\in\p\cup 2\p}=2p_{\pi(\frac{m}{2})}$$
\end{enumerate} 
In the case 
$$
w=\max\set{u\in||\C(m,\N)||~\mid u\in\p\cup 2\p}=p_{\pi(m)}
$$ 
we choose a prime number $x<m-w'$ such that $x>\frac{1}{2}(m-w')$, since by the Bertrand postulate (Lemma \ref{L_Bertrand's postulate}) there exists a prime number $x$ such that $x\in (k,2k)$ for every $k>1$ and set $2x=x'\in 2\mathbb{P}$. We obtain for the axis partner $[y']$ of the axis point $[x']$ of $\LL_{[x'],[y']}\inn\C(m+t,\N)$ the inequality 
\begin{align*}
y'=m+t-x'&\gtrsim m+t-2h(m)\frac{m}{\log m}\\
&=m\left(1-2\frac{h(m)}{\log m}\right)+t\\
&\sim m+t\geq p_{\pi(m)}=w
\end{align*}
for $t\geq 4$ and by the lemma \ref{L_the first little lemma} the following asymptotic inequalities
\begin{align*}
m-w\sim m-w'<x'
\end{align*}
and 
\begin{align*}
\abs{y'-w}=\abs{(m+t-x')-w}=\abs{m-w+t-x'}\lesssim\abs{x'+t-x'}= t.
\end{align*}
The requirements in Theorem \ref{T_Conditional Lemoine} are \textbf{asymptotically} fulfilled in this case
with
\[
y'\gtrsim w\mbox{ and }x'\gtrsim m-w\mbox{ and }0\lesssim \abs{y'-w}\lesssim t.
\]
In the case 
$$
w'=\max\set{u\in||\C(m,\N)||~\mid u\in\p\cup 2\p}=2p_{\pi(\frac{m}{2})}
$$ 
we choose a prime number $x<2(m-w)$ such that $x>(m-w)$, since by Bertrand's postulate (Lemma \ref{L_Bertrand's postulate}) there exists a prime number $x$ such that $x\in (k,2k)$ for every $k>1$. We obtain for the axis partner $[y]$ of the axis point $[x]$ of $\LL_{[x],[y]}\inn\C(m+t,\N)$ the inequality 
\begin{align*}
y=m+t-x&\gtrsim m+t-2h(m)\frac{m}{\log m}\\
&=m\left(1-2\frac{h(m)}{\log m}\right)+t\\
&\sim m+t\geq p_{\pi(m)}\sim 2p_{\pi(\frac{m}{2})}=w'
\end{align*}
for $t\geq 4$ and by appealing to Lemma \ref{L_the second little lemma} the following asymptotic inequalities
\begin{align*}
h(m)\frac{m}{\log m}\sim m-w'\sim m-w<x
\end{align*}
where $h(m):=o(\log m)$ and 
\begin{align*}
\abs{y-w'}=\abs{(m+t-x)-w'}=\abs{m-w'+t-x}\lesssim\abs{x+t-x}= t.
\end{align*}
The requirements in Theorem \ref{T_Conditional Lemoine} are \textbf{asymptotically} fulfilled in this second case
with
\[
y\gtrsim w'\mbox{ and }x\gtrsim m-w'\mbox{ and }0\lesssim \abs{y-w'}\lesssim t.
\] 
The result follows by arbitrarily choosing $t\geq 4$ so that $\mathcal{C}(m+t,\mathbb{P}\cup 2\mathbb{P})\neq \emptyset$ and adapting the proof in Theorem \ref{T_Conditional Lemoine}.
\end{proof}

\section{Further remarks}

Theorem \ref{T_Asymoptotic Goldbach theorem} and Theorem \ref{T_Asymoptotic Lemoine theorem} are both equivalent to the statement: there must exist some $N>0$ such that for all $m\geq N$ every even number $m$ can be written as the sum of two prime numbers (resp. every odd number is the sum of a prime and a double of a prime). This result--although constructive to some extent--loses its constructive flavour in a way that one cannot perform this construction to cover all odd numbers, since we are unable to obtain any quantitative (lower) bound for the threshold $N$. Regardless, we are able to \textbf{asymptotically} get a handle on the conjecture. The asymptotic version of the Lemoine conjecture implies an asymptotic version ternary of the Goldbach conjecture \cite{helfgott2013ternary}

%
%%%%%%%%%%%%%%%%%%%%%%%%%%%%%%%%%%%%%%%%%%%%%%%%%%%%%%%%%%%%%%%%%%%%%%%%
%\footnote{
%\par
%.}%
%%%%%%%%%%%%%%%%%%%%%%%%%%%%%%%%%%%%%%%%%%%%%%%%%%%%%%%%%%%%%%%%%%%%%%%%
\rule{100pt}{1pt}

\bibliographystyle{amsplain}

\end{document}